\documentclass[aps,prd,longbibliography,a4paper,floatfix,nofootinbib,12pt,report]{revtex4-1}

\usepackage{graphicx}
\usepackage[dvipsnames]{xcolor}
\usepackage[hidelinks]{hyperref}
\usepackage{amsmath,amssymb,amsfonts,amsthm,amscd,bm}
\usepackage{booktabs}
\usepackage{cancel}
\usepackage{url}
\usepackage{listings}
\usepackage{bookmark}

\usepackage{setspace}

\usepackage[hang]{footmisc}
\newcommand\blfootnote[1]{%
  \begingroup
  \renewcommand\thefootnote{}\footnote{#1}%
  \addtocounter{footnote}{-1}%
  \endgroup
}

\def\tilde{\widetilde}
\def\bar{\overline}
\def\hat{\widehat}
\def\*{\star}
\def\[{\left[}
\def\]{\right]}
\def\({\left(}      

\def\){\right)}

\def\zbar{{\bar{z} }}
\def\frac#1#2{{#1 \over #2}}
\def\frac#1#2{\dfrac{#1}{#2}}
\def\inv#1{{1 \over #1}}
\def\inv#1{\dfrac{1}{#1}}
\def\half{{1 \over 2}}
\def\half{\tfrac{1}{2}}
\def\d{\partial}

\def\2pi{\hbox{$2\pi i$}}

\def\dsl{\raise.15ex\hbox{/}\kern-.57em\partial}
\def\Dsl{\,\raise.15ex\hbox{/}\mkern-.13.5mu D}
\def\Zmath{\mathbb{Z}}

\def\sig{\sigma}

\def\vphi{\varphi}

      \def\CC{{\cal C}}
\def\CD{{\cal D}}   \def\CE{{\cal E}}   \def\CF{{\cal F}}
      \def\CI{{\cal J}}
      
\def\CM{{\cal M}}

\def\2pi{\hbox{$2\pi i$}}

\def\dsl{\raise.15ex\hbox{/}\kern-.57em\partial}
\def\Dsl{\,\raise.15ex\hbox{/}\mkern-.13.5mu D}

\def\barray{\begin{eqnarray}}
\def\earray{\end{eqnarray}}
\def\beq{\begin{equation}}
\def\eeq{\end{equation}}

\def\kvec{{\bf{k}}}
\def\gradvec{\vec{\nabla}}
\def\Evec{\vec{E}}

\def\vphi{\varphi}
\def\vphibar{\bar{\varphi}}
\def\smallhalf{\tfrac{1}{2}}
\def\phi{\Phi}

\def\Li{{\rm Li}}

\def\AA{\leavevmode\setbox0=\hbox{h}
\dimen0=\ht0 \advance\dimen0 by-1ex\rlap{\raise.67\dimen0\hbox{\char'27}}A}

\def\Arg{{\mathrm{Arg}\,}}

\makeatletter
\def\iddots{\mathinner{\mkern1mu\raise\p@
\vbox{\kern7\p@\hbox{.}}\mkern2mu
\raise4\p@\hbox{.}\mkern2mu\raise7\p@\hbox{.}\mkern1mu}}
\makeatother

\def\Li{{\rm Li}}
\def\yzero{y_\bullet}
\def\SL{SL_2\(\mathbb{Z}\)}
\def\curve{\CC} 
\def\yplus{y^{(+)}}
\def\yminus{y^{(-)}}

\makeatletter
\g@addto@macro\bfseries{\boldmath}
\makeatother

\newtheorem{theorem}{Theorem}

\theoremstyle{remark}

\newtheoremstyle{note}
{}
{}
{}
{}
{\bfseries}
{.}
{.5em}
{}
\theoremstyle{note}
\newtheorem{noteadded}[theorem]{Note added}
\newenvironment{added}[1][]{\begin{noteadded}[#1]}{\hfill$\blacktriangle$\end{noteadded}}

\begin{document}

\title{A theory for the zeros of Riemann \texorpdfstring{$\zeta$}{Zeta} 
and  other \texorpdfstring{$L$}{L}-functions (updated)}

\author{Guilherme Fran\c ca}
\email{guifranca@gmail.com}
\affiliation{Cornell University, Physics Department, Ithaca, New York 14850 USA} 

\author{Andr\' e  LeClair}
\email{andre.leclair@cornell.edu}
\affiliation{Cornell University, Physics Department, Ithaca, New York 14850 USA}%
\blfootnote{Original version: July 2014;  Updated version: August 2024} 


\begin{abstract}
In these lectures we first review the important 
properties of the Riemann $\zeta$-function that are necessary to understand the 
nature and importance of the Riemann hypothesis (RH). 
In particular this  first part describes 
the analytic continuation, the functional equation,
trivial zeros, the Euler product formula,
Riemann's main result relating the zeros on the critical strip to the 
distribution of primes, the exact counting formula for the number 
of zeros on the strip $N(T)$, and the GUE statistics
of the zeros on the critical line.
We then turn to presenting some new results obtained in the past year 
and describe several strategies towards proving the RH.    
First we describe an
electrostatic analogy and argue that if the electric potential 
along the line $\Re (z) =1$ 
is a regular alternating function,  the RH would follow.    
The main new result is that 
the zeros on the critical line are in one-to-one correspondence 
with the zeros of the cosine function,
and this leads to a transcendental equation for the $n$-th zero on 
the critical line that depends only
on $n$.    
If there is a unique solution to this equation for every $n$,  
denoting $N_0 (T)$ the number of zeros on the critical line,  
then $N_0(T) = N(T)$,  i.e. all zeros are on the critical 
line.    These results are generalized to two infinite classes of 
functions,  Dirichlet $L$-functions
and $L$-functions based on modular forms.   
We present extensive numerical analysis of 
the solutions of these equations.     
We apply these methods to the Davenport-Heilbronn 
$L$-function,  which is known to have zeros off of the line,   
and explain why the RH fails in this case.    
We also present a new approximation to the $\zeta$-function that 
is analogous to the Stirling approximation to the $\Gamma$-function.

In this updated version of these lectures,  we removed some   incorrect   comments   concerning 
$\arg \zeta \big(\tfrac12 + i y\big)$,  and added references to  subsequent publications,  the most recent in 2024 which made some progress on some open questions posed  in the previous version from 10 years ago.

\vspace{0.6cm}

\centerline{{\footnotesize 
Lectures delivered by A. LeClair at  Riemann Master School on $\zeta$eta Functions}}

\centerline{{\footnotesize
Riemann Center, Hanover, Germany,  June 10-14, 2014}}
\end{abstract}

\keywords{Riemann Hypothesis, Modular forms, $L$-functions, 
transcendental equations}

\maketitle

\tableofcontents

\newpage

\section{Introduction}

Riemann's zeta function $\zeta (z)$  plays a central role in many areas
of mathematics and physics.    It was present at the birth of Quantum Mechanics,
where the energy density of photons at finite temperature is proportional to
$\zeta (4)$ \cite{Planck}.    In mathematics,  it is at the foundations of 
analytic number theory.   

Riemann's major contribution to number theory was an explicit formula 
for the arithmetic function $\pi (x)$, which counts the number of primes 
less than $x$, in terms of an infinite sum over the non-trivial zeros of 
the $\zeta (z) $ function,  i.e. roots  $\rho$ of the equation $\zeta (z ) =0$ 
on the \emph{critical strip} $0\leq \Re (z)  \leq 1$ \cite{Edwards}.
It was later proven by Hadamard and de la Vall\'ee Poussin that there are 
no zeros on the line $\Re (z) =1$,  which in turn  proved the prime number 
theorem $\pi (x)  \sim  \Li (x)$;  see section \ref{sec:prime} for a review.  Hardy proved that there are an infinite number of zeros on the 
\emph{critical line} $\Re (z) = \tfrac{1}{2}$. The \emph{Riemann hypothesis} 
(RH) was his  statement, in 1859,  that all zeros on the critical strip 
have $\Re(\rho)=\tfrac{1}{2}$.   
Despite strong numerical evidence of its validity,  it remains unproven 
to this day.  Some excellent 
introductions to the RH are \cite{Conrey,Sarnak,Bombieri}.
In this context, the common convention is that the argument of $\zeta$ is
$s=\sigma + i t$.   Throughout these lectures, the argument of $\zeta$ will
instead be denoted as $z=x+i y$,  zeros will be denoted as $\rho$,   and zeros on the positive
critical line will be denoted as $\rho_n =  \smallhalf + i y_n$ for $n=1,2,\dotsc$.   

The aim of these lectures is clear from the Table of Contents.    
Not all the material was actually lectured on,  since the topics in 
some of the first sections
were covered by other lecturers in the school.   
There are two main components. Section~\ref{QuantStat}   
reviews the most important applications of $\zeta$ to quantum statistical 
mechanics, and discusses other physics-based approaches,   but 
is not necessary for understanding the rest of the lectures.     
 The following three sections  review 
the most important properties of Riemann's zeta function.   This introduction
to $\zeta$ is  self-contained, and describes all the  necessary ingredients to 
understand the importance of the Riemann zeros and the nature of the 
Riemann hypothesis, and may perhaps be more accessible to mathematical physicists than the more advanced mathematics literature.     
The remaining sections present new work carried out over the past year.    
Sections~\ref{Electric}--\ref{sec:dirichlet} 
are based on the published works \cite{RHLeclair,FL1}.    
The main new result 
is a transcendental equation satisfied by individual 
zeros of $\zeta$ and other $L$-functions.
These equations provide a novel characterization of the zeros,  
hence the title 
``A theory for the zeros of Riemann $\zeta$ and other $L$-functions''.          
Sections~\ref{Soft},  \ref{sec:Davenport} and \ref{sec:saddle}  
contain additional results that have not been previously published.  

Let us summarize some of the main points of the material presented  in 
the second half of these 
lectures.      
It is difficult to visualize a complex function since it is a hypersurface
in a $4$-dimensional space.
We present one way to visualize
the RH based on an electrostatic analogy.   
From the real and imaginary parts of the function 
$\xi (z)$ defined in \eqref{xidef}  we construct a 
2-dimensional vector field 
$\Evec$ in section \ref{Electric}.    
By virtue of the Cauchy-Riemann equations, 
this field satisfies the conditions of an electrostatic field with no 
sources.  It can be
written as the gradient of an electric potential $\Phi (x,y)$,  
and visualization is reduced to 
one real function over the 2-dimensional complex plane.
We argue that if the real  potential $\Phi$  along the line
$\Re (z) =1$ alternates between positive and negative values in the most regular 
manner possible,
then the RH would appear to be true.   
Asymptotically one can analytically understand this
regular alternating behavior,  however the asymptotic expansion 
is not controlled enough to 
rule out zeros off the critical line.

The primary  new result is described in section \ref{sec:zeta_function}.    
There we present
the equation \eqref{exact_eq2} which is a transcendental equation for 
the ordinate of 
the $n$-th zero on the critical line.\footnote{\emph{Note added:} We will sometimes refer to this equation as the Fran\c ca-LeClair equation \cite{FL1}.}  These zeros are in one-to-one 
correspondence with
the zeros of the cosine function.    
This equation involves two terms,  a smooth one
from $\log \Gamma$ and a small fluctuating term 
involving $\arg \zeta \(\smallhalf + i y \)$,  i.e.
the function commonly referred to as $S(y)$ \eqref{SofyDef}. 
If the small term $S(y)$ is neglected,
then there is a  solution to the equation for every $n$ which  
can be expressed explicitly in terms of the  Lambert $W$-function  
(see section~\ref{sec:lambert}).
The equation with the $S(y)$ term
can be solved numerically to calculate zeros to very high accuracy,  
$1000$ digits or more. 

More importantly,   there is  a clear strategy 
for proving the RH based on the Fran\c ca-LeClair equation.    
It is easily stated: if there is a unique solution 
to \eqref{exact_eq2} for each $n$,   then the zeros can be counted 
along the critical line
since they are enumerated by the integer $n$.    
Specifically, one can determine  
$N_0 (T)$ which is the number of zeros on the critical line 
with $0<y<T$ (see eq.~\eqref{counting2_exact}).
On the other hand,   there is a known expression for $N(T)$ which counts the 
zeros in the entire critical strip due to Backlund; the asymptotic expansion of 
$N(T)$ was known to Riemann.\footnote{Riemann only presented the asymptotic limit of $N(T)$ as $T \to \infty$,  but
must have known the exact formula eventually written by Backlund.}    
We find that $N_0 (T) = N(T)$,  which indicates that
all zeros are on the critical line.     
The proof is not complete mainly because we cannot 
establish that there is a unique solution to \eqref{exact_eq2} for 
every positive integer $n$ 
due to the fluctuating  behavior of the function $S(y)$.
Nevertheless,  we argue that the 
$\delta \to 0^+$ prescription smooths out the discontinuities of 
$S(y)$ sufficiently 
such that there are unique solutions for every $n$.
There is extensive numerical evidence 
that this is indeed the case,  presented in sections \ref{sec:numerical}, 
\ref{sec:GUE}, and \ref{sec:numerical_exact}.

Understanding the  properties of $\arg \zeta\big(\half + i y\big)$ is crucial to our construction,
and section \ref{Soft} is devoted to it.   

$L$-functions are generalizations of the 
Riemann $\zeta$-function, the latter being the trivial case~\cite{Apostol}.
It is straightforward to extend the results for $\zeta$ to two 
infinite classes of important $L$-functions, the Dirichlet $L$-functions and
$L$-functions associated with cusp (modular) forms.
The former have applications primarily in multiplicative number theory,
whereas the latter in additive number theory.
These functions can be analytically continued to the entire 
(upper half) complex plane. The \emph{Generalized Riemann Hypothesis}  
and \emph{Grand Riemann Hypothesis} refers to these classes of functions.  
These would  also would follow if there were a unique solution for each $n$
of the appropriate transcendental equation,  which we present below.    

There is a well-known counterexample to the RH which is based on 
the Davenport-Heilbronn function.   
It is an $L$-series that is  a linear combination
of two Dirichlet $L$-series,  that also satisfies a functional equation.    
It has an infinite number of zeros on the critical line,  
but also zeros off  of it, thus violating the RH.
It is therefore interesting to apply our construction  to this function.
One finds that the zeros on the line do satisfy a transcendental equation, 
where again  the solutions are enumerated by an integer $n$.    However 
for some $n$ there is no solution, and this happens precisely where there
are zeros off of the line.  This can be traced to the behavior of the 
analog of $S(y)$,  which changes branch at these points.
Nevertheless the zeros off the line continue to satisfy our general 
equation  \eqref{Blimit}.

The final topic concerns a new approximation for $\zeta$.
Stirling's approximation to $n!$ is extremely useful; much of 
statistical mechanics would be impossible without it.
Stirling's approximation to $\Gamma (n) = (n-1)!$  is a saddle point
approximation,  or steepest decent approximation,  to the 
integral representation for $\Gamma$ when $n$ is real and positive.  
However since $\Gamma(z)$ is an analytic function of the complex variable $z$,
Stirling's approximation extends to the entire complex plane for $z$ large.    
We describe essentially the same kind of analysis for the $\zeta$ function.  
Solutions to the saddle point equation are explicit in terms of the Lambert 
$W$-function.     The main complication compared to $\Gamma$ is that 
one must sum over more than one branch of the $W$-function,  
however the approximation is quite  good.

\begin{added} 
Since these are unpublished  lectures,  we took the liberty to update them  now  10 years later.
We removed some incorrect comments about the behavior of  the argument of $\zeta \big(\tfrac12 + i y\big)$,   the notorious  ``$S(t)$",  about  which we were somewhat naive at the time due to our limited background in Analytic Number Theory.\footnote{We originally  thought that our 
$\arg \zeta \big(\tfrac12 + i y\big)$ \emph{on the zeros}  was not identical to  the standard definition of $S(t=y)$ from piecewise integration.   Furthermore, we were misled by numerical studies that that did not go high enough up the critical line to observe the extremely slow growth of   $\arg \zeta \big(\tfrac12 + i y\big)$,   which is roughtly $\log \log y$.}
This is the only significant modification of the original content,  which is in  section~\ref{Soft},  apart from some corrected typos.
The more important changes are additional remarks concerning work done by us subsequent to the 2014 lectures \cite{FrancaLeclairRW}, and some followup work  
done in 
collaboration with Giuseppe Mussardo \cite{Googleth,MLRW1,LMRW,MLRW2,LecMussScattering,SpectralFlow}.   The  latest articles  \cite{LecMussScattering,SpectralFlow} only appeared this year of 2024.     The latter works are progresses on some of the open questions raised in this original lectures from 2014.
 However,  we do not review these new works in any detail in this update, but rather comment on them mostly in the final discussion section~\ref{sec:discussion}.  These additional remarks are 
highlighted  throughout  
such that the reader can easily spot them.\footnote{Additional references can be found in the published versions of these articles and others mentioned in this updated version of the lectures.}  
\end{added}

\section{Riemann \texorpdfstring{$\zeta$}{Zeta} in quantum 
statistical mechanics and physics based approaches}
\label{QuantStat}

In this section we describe some significant occurrences of Riemann's 
zeta function 
in the quantum statistical mechanics of free gases.   
We also discuss approaches to the RH based on Brownian Motion in one spatial dimension,  or Random Walks. 
 For  other  interesting connections of the RH
to physics see \cite{Connes2,Schumayer} (and references therein).
One prominent idea is the Hilbert-P\'olya conjecture,  where they suggested
that for zeros expressed in the form $\rho_n =  \smallhalf + i y_n$,  
the $y_n$ are the real  eigenvalues of
a hermitian operator, perhaps an unknown  quantum mechanical Hamiltonian. 
Most of this literature are variations of Berry-Keating's Hamiltonian $H= xp$ \cite{BerryKeating,Sierra3},  and we do not include these ideas in these 
 review lectures.   We refer the reader to various reviews such as \cite{Schumayer}.     But it is fair to say that such a Hamiltonian has remained elusive in its specifics,   and it is not at all clear if it exists as a bound state problem in Quantum Mechanics.

\begin{added}    
Recently we proposed that 
the Riemann zeros arise as quantized energies of a scattering problem in Quantum Mechanics,  rather than a bound state problem.  
Although this is not what was presumably envisioned by Hilbert-P\'olya,    it is in a similar spirit and we believe is the 
 proper version of the Hilbert-P\'olya 
conjecture.       
We defined a model of particles scattering about a circle of impurities \cite{LecMussScattering},   where each impurity is associated with a scattering phase such that the S-matrix is unitary.    This relates to the    Hilbert-P\'olya conjecture since roughly speaking $S=e^{-i H}$,  where $H$ is a hermitian Hamiltonian, 
thus $S^\dagger S = 1$.\footnote{This is an over-simplification since one actually needs to consider a time-ordered exponential.}           
In this model the Riemann zeros on the critical line arise from a Bethe-ansatz equation which is readily solved for the \emph{exact}   Riemann zeros.    The spectral flow for these eigen-energies was studied in \cite{SpectralFlow},   where it was argued that the scattering problem necessarily has real eigenvalues since the 
S-matrix is unitary,  and this would imply the Riemann Hypothesis is true.  Further remarks are in the last  discussion section~\ref{sec:discussion}.
\end{added}

\begin{added}
 Another approach to proving the RH based on physical ideas involves  Brownian Motion or Random Walks \cite{FrancaLeclairRW,LMRW,MLRW1,MLRW2}. 
There it was argued that certain sums of non-principal Dirichlet characters over prime numbers  behave like 
Brownian motion in one spatial dimension.  Namely they behave as independent,   identically distributed random variables,  hence the central limit theorem applies.    In this approach,  the significance of the critical line $\Re (z) = \half$ is identical to the universal critical exponent  $1/2$ 
for the standard deviation $N^{1/2}$ for $N$ steps.\footnote{$\half = 1/2$!   However this by itself is not very meaningful since the $L$ functions can easily be redefined to shift the critical line,   whereas the random walk exponent $1/2$ is universal.    For instance,  for $L$-functions based on cusp forms of weight $k$,   the critical line is $\Re (z)= k/2$.}   
If this is the case,   then the Euler product  converges to the right of the critical line and this would also establish the validity of the RH for these 
$L$-functions.    In this approach the random matrix (GUE) statistics reviewed below is intimately related to this Brownian motion.  
These  Brownian motion  arguments also extend to $L$-functions based on cusp forms \cite{FrancaLeclairRW}. 
For the principal Dirichlet $L$ functions,  such as $\zeta (z)$ itself,   these random walk ideas also apply \cite{MLRW2} when one considers M\"obius inversion for $1/\zeta (z)$.   An extensive battery of statistical tests on such series were performed in \cite{MLRW2} and rather convincing evidence was presented that non-principal  Dirichlet characters $\chi (p)$ where $p$ is prime, or M\"obius $\mu$ over the integers,  can out-perform current random number generators.    The main obstruction  for proving the RH is that  pseudo-randomness of the primes is difficult to formulate
rigorously,  except perhaps along the lines just stated based on the Dirichlet characters. 
     Assuming the particular series involved are indeed random walks,   then 
the Riemann Hypothesis can only be proven \emph{almost surely}   in probability theory, 
namely with probability equal to 1,    and this is somewhat of a brick wall toward proving   that it is  
\emph{surely true}.
\end{added}

\subsection{The quantum theory of light:  physical derivation of the functional equation}

Perhaps the most important application of  Riemann's zeta function in physics 
occurs in quantum statistical mechanics.    One may argue that 
Quantum Mechanics was born
in Planck's seminal paper on black body radiation \cite{Planck},  
and $\zeta$  was present at this birth; 
in fact his study led to the first determination
of the fundamental Planck constant $h$,  or more commonly $\hbar =h/2\pi$.   

In order to explain discrepancies between the spectrum of radiation 
in a cavity at temperature $T$
and the prediction of classical statistical mechanics,  
Planck proposed that 
the energy of radiation was ``quantized'',  i.e. took on only the 
discrete values 
$\hbar \, \omega_\kvec$,  where $\kvec$ is the wave-vector,  
or momentum,  and $\omega_\kvec = |\kvec| $.    
It is now well understood that these quantized energies are those 
of  real light particles called photons.      The following integral,  
related to the energy density of photons,  appears in Planck's paper:
\beq
\label{Planck1}
\frac{8 \pi h}{c^3}  \int_0^\infty   \frac{\nu^3  d\nu}{e^{h\nu/k} -1 }  
  =   \frac{8 \pi h}{c^3} \int_0^\infty  \nu^3  
  \(e^{-h\nu/k}  +  e^{-2h\nu/k}  +  e^{-3h\nu/k} + \dotsm \) d\nu .
\eeq
He then proceeds to evaluate this expression by termwise integration and writes
\beq
\label{Planck2} 
\frac{8 \pi h}{c^3}  \int_0^\infty   \frac{\nu^3  d\nu}{e^{h\nu/k} -1 }  =  
\frac{8 \pi h}{c^3}  \cdot 6 \( \frac{k}{h} \)^4 \( 1 + \inv{2^4}  +  
\inv{3^4} + \inv{4^4}  \ldots \)  .
\eeq
It is not clear whether Planck knew that the above sum 
was $\zeta (4)$,   since he simply
writes that it is approximately $1.0823$ since it converges rapidly.    
Due to Euler,  it  was already known at the time 
that $\zeta (4) = \pi^4/90$.  

The quantum statistical mechanics of photons leads to a physical 
demonstration  of the most important  functional equation satisfied 
by $\zeta (z)$.   
Consider a free quantum field theory of massless bosonic particles
in $d+1$ spacetime dimensions with Euclidean action
$S = \int d^{d+1} x  ~ (\d \phi)^2 $.   
The geometry of Euclidean spacetime is taken to be $S^1 \times \mathbb{R}^d$
where the circumference of the circle $S^1$ is $\beta$.   We will
refer to the $S^1$ direction as $\hat{x}$.  
Endow the flat space
$\mathbb{R}^d$ with a large but finite volume as follows.  Let us refer
to one of the directions perpendicular to $\hat{x}$ as $\hat{y}$
with length $L$ and let the remaining $d-1$ directions have volume $A$. 

\begin{figure}
\centering
\includegraphics{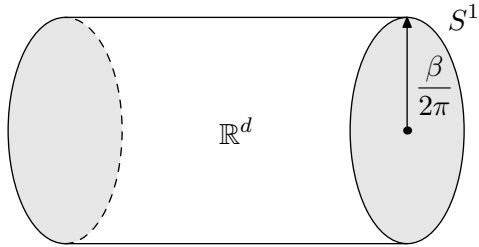}
\caption{Spacetime geometry for the partition function in \eqref{ap.1}.}
\label{fig_Casimir} 
\end{figure}

Let us first view the $\hat{x}$ direction as compactified Euclidean time,
so that we are dealing with finite temperature $T=1/\beta$
(see Figure~\ref{fig_Casimir}).    
As a quantum statistical mechanical system, the partition function 
in the limit $L, \, A \to \infty$ is
\beq
\label{ap.1}
Z = e^{-\beta V \, \CF (\beta )}
\eeq
where $V=L\cdot A$ and $\CF$ is the free energy density.  
Standard results give:
\beq
\label{ap.2}
\CF (\beta ) = \inv{\beta} \int \frac{d^d \kvec}{(2\pi)^d} \,
\log \( 1 - e^{-\beta k} \).
\eeq

The Euclidean rotational symmetry allows one to view the above
system with time now along the $\hat{y}$ direction.  In $1d$,
interchanging the role of $\hat{x}$ and $\hat{y}$ is a special case
of a modular transformation of the torus.   In this version, 
the problem is a zero temperature quantum mechanical system with
a finite size $\beta$ in one direction,  and the total volume of
the system is $V'= \beta \cdot A$.   The quantum mechanical path
integral leads to 
\beq
\label{ap.3}
Z = e^{- L E_0 (A,\beta)}
\eeq
where $E_0$ is the ground state energy.  Let $\CE_0 = E_0/V'$ 
denote the ground state energy per volume.  Comparing the two
``channels'',  their equivalence requires $\CE_0 (\beta) = \CF (\beta)$.  
In this  finite-size channel,  the modes of the field in the
$\hat{x}$ direction are quantized with wave-vector $k_x = 2 \pi n/\beta$,
and the calculation of $\CE_0$ is as in the Casimir effect (see below):
\beq
\label{ap.4}
\CE_0 = \inv{2\beta} \sum_{n \in \Zmath} 
\int \frac{ d^{d-1} \kvec}{ (2\pi)^{d-1} } 
\( \kvec^2 + \(2\pi n /\beta\)^2 \)^{1/2} .
\eeq

The free energy density $\CF$ can be calculated using 
\beq \int d^d \kvec = 
\frac{ 2\pi^{d/2} }{\Gamma(d/2)} \int  d|\kvec|  \, |\kvec|^{d-1}.
\eeq
For $d>0$ the integral is convergent
and one finds
\beq
\label{ap.5}
\CF = -  \inv{\beta^{d+1}}  \,
\dfrac{ \Gamma (d+1) \, \zeta (d+1)}{2^{d-1} \pi^{d/2} d \, \Gamma(d/2) }. 
\eeq

For the Casimir energy,  after performing the $\kvec$ integral,   
$\CE_0$ involves 
$\sum_{n \in \Zmath} |n|^d $ which must be regularized.  
As is usually done, let us regularize this as $2 \, \zeta (-d)$.  
Then
\beq
\label{ap.6}
\CE_0 = - \inv{\beta^{d+1}} ~ \pi^{d/2} \Gamma(-d/2) \zeta (-d) .
\eeq

Let us analytically continue in $d$ and define the function 
\beq
\label{ap.7}
\chi (z) \equiv  \pi^{-z /2} \Gamma(z/2) \zeta (z) .
\eeq
Then the equality $\CE_0 = \CF$ requires the  identity
\beq
\label{ap.8}
\chi (z ) = \chi (1-z) .
\eeq
The above relation is a known functional identity that 
was proven by Riemann  using complex analysis.   
It will play an essential role in the rest of these lectures.   

The above calculations show   that $\zeta$-function regularization
of the ground state energy $\CE_0$  is consistent with a modular 
transformation to the finite-temperature channel.   
Our calculations can  thus be viewed as a
proof of the identity \eqref{ap.8} based on physical consistency.

\begin{figure}
\centering
\includegraphics{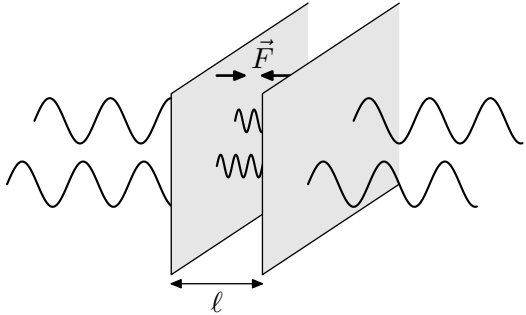}
\caption{Geometry of the Casimir effect.}
\label{fig_Casimir2} 
\end{figure}

For spatial dimension $d=3$,   the ground state energy $\CE_0$ is 
closely related
to the measurable Casimir effect,  the difference only being in the 
boundary conditions.
In the Casimir effect one measures the ground state energy of the 
electromagnetic field
between two plates by measuring the force as one varies their separation,  
as illustrated in Figure~\ref{fig_Casimir2}.

There is a simple relation between the vacuum energy densities 
$\rho$ in the cylindrical
geometry above,  and that of the usual Casimir effect:
\beq
\label{Casimir2} 
\rho_{\rm vac}^{\rm casimir} (\ell )  =  
2 \rho_{\rm vac}^{\rm cylinder}  (\beta = 2 \ell)  =
- \frac{\pi^2}{720 \ell^4 }.
\eeq
It is remarkable that since the Casimir effect has been measured in 
the laboratory, 
such a measurement verifies 
\beq
\label{Casimir3}
\zeta (-3)  =  1 + 2^3 + 3^3 + 4^3  + \dotsm  =  \inv{120}.
\eeq
Of course the above equation is non-sense on its own.   
It only makes sense after
analytically continuing $\zeta (z)$ from $\Re (z)>1$ to the 
rest of the complex plane.

\subsection{Bose-Einstein condensation and the pole of 
\texorpdfstring{$\zeta$}{Zeta}}

It is known that $\zeta (z)$  has only one pole,  a simple one at $z=1$.   
This property is also related to some basic physics.     
In the phenomenon of
Bose-Einstein condensation,  below a critical temperature $T$ all of 
the bosonic particles 
occupy the lowest energy single particle state.   
This critical temperature depends on 
the density $n$.   In $d$ spatial dimensions the formula reads 
\beq
\label{BEC}
n = \int \frac{d^d \kvec}{(2 \pi)^d } \inv{e^{\omega_\kvec /T}  -1}  
=  \( \frac{mT}{2 \pi} \)^{d/2} \zeta (d/2)  .
\eeq 
The Coleman-Mermin-Wagner theorem in statistical physics says that 
Bose-Einstein condensation is impossible in $d=2$   spatial dimensions.   
This is manifest in the above formula since 
$\zeta (1) = \infty$.

\section{Important properties of \texorpdfstring{$\zeta$}{Zeta}}

After having discussed some applications of $\zeta$ in physics, 
let us now  focus on its most basic mathematical properties as a complex
analytic function.

\subsection{Series representation}

The $\zeta$-function is defined for $\Re(z) > 1$ through the series
\beq\label{zeta_series}
\zeta(z) = \sum_{n=1}^{\infty}\dfrac{1}{n^z}.
\eeq
The first appearance of $\zeta(z)$ was in 
the so called ``Basel problem'' posed by Pietro Mengoli in
$1644$. This problem consists in finding the precise sum
of the infinite series of squares of natural numbers:
\beq\label{basel}
\zeta(2) = \sum_{n=1}^{\infty} \dfrac{1}{n^2} = \dfrac{1}{1^2} + 
\dfrac{1}{2^2} + \dfrac{1}{3^2} + \dotsb = \ ?
\eeq
The leading mathematicians of that time, like the Bernoulli 
family, attempted  the problem unsuccessfully. It was only in
$1735$ that Euler, only $28$ years old, solved the problem 
claiming that such sum is equal to $\pi^2/6$, and was
brought to fame. Nevertheless, his arguments were not fully 
justified,  as he manipulated infinite series with abandon,  
 and it was only in $1741$ that Euler could provide  a formal proof.
Even such a proof had to wait $100$ years until all the steps
could be rigorously justified by the Weierstrass factorization 
theorem.
Just out of curiosity let us reproduce Euler's steps. From the Taylor series
of $\sin z$ we have 
\beq\label{taylor_sin}
\dfrac{\sin\pi z}{\pi z} = 
1 - \dfrac{\(\pi z\)^2}{3!} + \dfrac{(\pi z)^4}{5!} + \dotsb.
\eeq
From the Weierstrass product formula for the $\Gamma(z)$ function it is
possible to obtain the product formula for $\sin z$ which reads
\beq\label{prod_sin}
\dfrac{\sin\pi z}{\pi z} = 
\prod_{n=1}^{\infty}\(1-\(\dfrac{z}{n}\)^2\).
\eeq
Now collecting terms of powers of $z$ in \eqref{prod_sin} one finds
\beq\label{prod_sin_sum}
\dfrac{\sin\pi z}{\pi z} = 1 - z^2 \sum_{n=1}^{\infty}\dfrac{1}{n^2}
+ z^4 \sum_{\substack{n,m=1 \\ m > n}}^{\infty}\dfrac{1}{n^2 m^2} + 
\dotsb.
\eeq
Comparing the $z^2$ coefficient of \eqref{prod_sin_sum} and 
\eqref{taylor_sin} we immediately obtain Euler's result,
\beq\label{zeta_of_2}
\zeta(2) = \sum_{n=1}^{\infty}\dfrac{1}{n^2} = \dfrac{\pi^2}{6}.
\eeq
Comparing the $z^4$ coefficients we have
\beq\label{coeff4}
\sum_{\substack{n,m=1 \\ m > n}}^{\infty}\dfrac{1}{n^2 m^2} = 
\dfrac{\pi^4}{5!}.
\eeq
Now consider the full range sum
\beq
\sum_{n,m=1}^{\infty}\dfrac{1}{n^2m^2} =
\sum_{n=1}^{\infty}\dfrac{1}{n^4} + 
\sum_{\substack{n,m=1 \\ m > n}}^{\infty}\dfrac{1}{n^2m^2} + 
\sum_{\substack{n,m=1 \\ m < n}}^{\infty} \dfrac{1}{n^2m^2}.
\eeq
The LHS is nothing but $\zeta^2(2)$. The first term in the RHS is
$\zeta(4)$ and the two other terms are both equal to \eqref{coeff4}. 
Therefore we have
\beq\label{zeta_of_4}
\zeta(4) = \sum_{n=1}^{\infty}\dfrac{1}{n^4} = \dfrac{\pi^4}{90}.
\eeq
Considering the other powers $z^{2n}$ we can obtain $\zeta(2n)$ for even
integers. We will see later a better approach to compute these values
through an interesting relation with the Bernoulli numbers.
For positive odd integer $n$,  $\zeta (n)$ has no such simple expression, and this is related to its pole at $n=1$.

\subsection{Convergence}

Let us now analyze the domain of convergence of 
the series \eqref{zeta_series}. Absolute convergence implies convergence,
therefore it is enough to consider
\beq
\sum_{n=1}^{\infty}\dfrac{1}{|n^z|} = \sum_{n=1}^{\infty}\dfrac{1}{n^x},
\eeq
where $z=x+iy$. Let $x=1+ \delta $ with $\delta > 0$. By the integral
test we have
\beq\label{integral_test}
\int_{1}^{\infty}\dfrac{1}{u^{1+\delta}} \, du = 
-\dfrac{1}{\delta u^\delta }\bigg\vert_{1}^{\infty} = \dfrac{1}{\delta}.
\eeq
If $\delta > 0$ the above integral is finite, therefore the series is
absolutely convergent for $\Re(z) > 1$ and $\zeta(z)$ is an analytic function
on this region. Note that if $\delta = 0$ the above
integral diverges as $\log(\infty)$. If $\delta \le -1$ the integral
also diverges. Therefore, the series representation given by 
\eqref{zeta_series} is defined only for $\Re(z) > 1$.

\subsection{The golden key: Euler product formula}

The importance of the $\zeta$-function in number theory is mainly because 
of its relation to prime numbers. 
The first one to realize this connection was again Euler.
To see this, let us consider the following product
\beq\label{prod1}
\prod_{i=1}^{\infty}\(1-\dfrac{1}{p_i^{z}}\)^{-1}
\eeq
where $p_i$ denotes a prime number, i.e. $p_1=2$, $p_2=3$ and so on.
We know that $\(1-z\)^{-1} = 
\sum_{n=0}^\infty z^{n}$ for $|z| < 1$, thus \eqref{prod1} is equal to
\beq\label{prod2}
\prod_{i=1}^{\infty}\sum_{n=0}^{\infty}\dfrac{1}{p_i^{nz}} = 
\(1+\dfrac{1}{p_1^z} + \dfrac{1}{p_1^{2z}}+\dotsb\)\times
\(1+\dfrac{1}{p_2^z} + \dfrac{1}{p_2^{2z}}+\dotsb\)\times \dotsb .
\eeq
If we open the product in \eqref{prod2} we have an infinite sum of terms,
each one having the form 
\beq \(p_1^{\alpha_1}p_2^{\alpha_2}
\dotsc p_j^{\alpha_j}\)^{-z} \quad \mbox{where} \ j=1,2,\dotsc \ \mbox{and} \
\alpha_j=0,1,2,\dotsc,
\eeq
i.e. we have a finite product of primes raised to every possible power. 
From the fundamental theorem of
arithmetic, also known as the unique prime factorization theorem, we know 
that every natural number can be expressed in exactly one way through a
product of powers of primes; 
$n=p_1^{\alpha_1}p_2^{\alpha_2}\dotsc
p_j^{\alpha_j}$ where $p_1 < p_2 < \dotsb < p_j$. 
Therefore \eqref{prod2} involves a 
sum of all natural numbers, which is exactly the definition
\eqref{zeta_series}. Thus we have the very important
result known as the \emph{Euler product formula}
\beq\label{euler_prod}
\zeta(z) = \sum_{n=1}^{\infty}\dfrac{1}{n^z} = 
\prod_{i=1}^{\infty}\(1-\dfrac{1}{p_i^z}\)^{-1}.
\eeq
This result is of course only valid for $\Re(z) > 1$. 
A simpler derivation of it is also given in section~\ref{sec:Gauss}.   
From this formula we can easily see that there are no zeros in this
region,
\beq
\zeta(z) \ne 0 \qquad \mbox{for $\Re(z) > 1$,}
\eeq
since each factor $\(1-p_i^{-z}\)$ never diverges.

\subsection{The Dirichlet \texorpdfstring{$\eta$}{Eta}-function}

Instead of the $\zeta$-function defined in \eqref{zeta_series}, let us
consider its alternating version
\beq\label{dirichlet_eta}
\eta(z) = \sum_{n=1}^{\infty}\dfrac{(-1)^{n-1}}{n^z}.
\eeq
This series is known as the \emph{Dirichlet $\eta$-function}. 
For an alternating
series $\sum (-1)^{n-1}a_n$, if $\lim_{n\to\infty}a_n = 0$ and 
$a_n > a_{n+1}$ then the series converges. We can see that both conditions
are satisfied for \eqref{dirichlet_eta} provided that $\Re(z) > 0$.

The $\eta$-function \eqref{dirichlet_eta} has a negative sign at even
naturals:
\beq\label{sum1}
\eta(z) = 1 - \dfrac{1}{2^z} + \dfrac{1}{3^z} - \dfrac{1}{4^z} + \dotsb.
\eeq
On the other hand the $\zeta$-function \eqref{zeta_series} has
only positive signs, thus if we multiply it by $2/2^z$ we double all
the even naturals appearing in \eqref{sum1} but with a positive sign
\beq\label{sum2}
\dfrac{2}{2^z}\zeta(z) = \dfrac{2}{2^z} + \dfrac{2}{4^z} + 
\dfrac{2}{6^z} + \dotsb.
\eeq
Summing \eqref{sum1} and \eqref{sum2} we obtain $\zeta(z)$ again, therefore
\beq\label{zeta_eta}
\zeta(z) = \dfrac{1}{1-2^{1-z}} \, \eta(z).
\eeq
In obtaining this equality we had to assume $\Re(z) > 1$, nevertheless
\eqref{dirichlet_eta} is defined for $\Re(z) > 0$, thus \eqref{zeta_eta}
yields the analytic continuation of $\zeta(z)$ in the region $\Re(z) > 0$.

\subsection{Analytic continuation}

One of Riemann's main contributions is the analytic continuation of
the $\zeta$-function to  the entire complex plane (except for a simple
pole at $z=1$). 
Actually, Riemann was the first one to consider
the function \eqref{zeta_series} over a complex field. Following his
steps, let us start from the integral definition of the $\Gamma$-function
\beq\label{gamma_int}
\Gamma(z) = \int_{0}^{\infty} u^{z-1} e^{-u} du.
\eeq
Replacing $u \to n u$ and summing over $n$ we have\footnote{
According to Edwards book \cite{Edwards} this formula also appeared 
in one of Abel's paper and a similar one in a paper of Chebyshev. 
Riemann should probably be aware of this.}
\beq\label{sum_int1}
\Gamma(z)\sum_{n=1}^{\infty}\dfrac{1}{n^z} = 
\sum_{n=1}^{\infty}\int_{0}^{\infty} u^{z-1} e^{-nu} du.
\eeq
If we are allowed to interchange the order of the sum with the integral, and
noting that  
$\sum_{n=1}^{\infty}e^{-nz} = \(e^z - 1\)^{-1}$, then formally we have
\beq\label{sum_int2}
\Gamma(z)\zeta(z) = \int_{0}^{\infty} \dfrac{u^{z-1}}{e^u - 1}du.
\eeq
The integral in \eqref{sum_int1} converges for $\Re(z) > 1$. Under
this condition it is possible to show that the integral \eqref{sum_int2}
also converges. Thus the step in going 
from \eqref{sum_int1} to \eqref{sum_int2} is justified as long as we assume
$\Re(z) > 1$.

\begin{figure}[h]
\begin{center}
\includegraphics{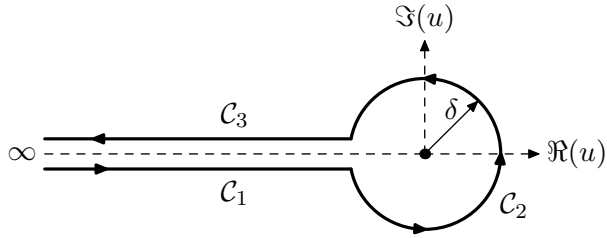}
\caption{Contour $\CC=\CC_1 + \CC_2 + \CC_3$ excluding the origin. 
$\CC_2$ has radius $\delta < 2\pi$.}
\label{fig_contour} 
\end{center}
\end{figure}

Since the integrand in \eqref{sum_int2} diverges at $u=0$, 
let us promote $u$ to a complex variable and consider the following 
integral
\beq\label{I_def}
\CI(z) = \dfrac{1}{2\pi i}\int_{\CC} \dfrac{u^{z}}{e^{-u} - 1}\dfrac{du}{u}
\eeq
over the path $\CC$ illustrated in
Figure~\ref{fig_contour}, which excludes the origin through a circle
of radius $\delta < 2\pi$, avoiding the pole of the integrand in \eqref{I_def}.

On $\CC_1$ we choose $u=r e^{-i\pi}$ then
\beq
\dfrac{1}{2\pi i}\int_{\CC_1} \dfrac{u^{z}}{e^{-u} - 1}\dfrac{du}{u}
= - \dfrac{e^{-i\pi z}}{2\pi i}
\int_{\delta}^{\infty} \dfrac{r^{z-1}}{e^r - 1}dr.
\eeq
Analogously, on $\CC_3$ we choose $w=r e^{i\pi}$, then we have
\beq\label{int_c1_c3}
\lim_{\delta\to 0}\dfrac{1}{2\pi i}\(\int_{\CC_1} + \int_{\CC_3}\)
\dfrac{u^{z}}{e^{-u} - 1}\dfrac{du}{u}
= \dfrac{ \sin \pi z }{\pi} \int_{0}^{\infty} \dfrac{r^{z-1}}{e^r - 1}dr.
\eeq
For the integral over $\CC_2$ we have $u = \delta e^{i\theta}$,
and since we are interested in $\delta\to 0$ we obtain
\beq\label{intc2}
\dfrac{1}{2\pi i}\int_{\CC_2} \dfrac{u^{z}}{e^{-u} - 1}\dfrac{du}{u}
= \dfrac{\delta^z}{2\pi} \int_{-\pi}^{\pi}\dfrac{e^{i\theta z}}
{e^{-\delta e^{i\theta}}-1} d\theta 
\approx -\dfrac{\delta^{z-1}}{2\pi}\int_{-\pi}^{\pi}
e^{i\theta(z-1)}d\theta.
\eeq
Thus we have
\beq
\left| 
\dfrac{1}{2\pi i}\int_{\CC_2} \dfrac{u^{z}}{e^{-u} - 1}\dfrac{du}{u}
\right| \le 
\dfrac{\delta^{x-1}}{2\pi}\int_{-\pi}^{\pi} e^{- y \theta} d\theta
= \delta^{x-1}\dfrac{\sinh \pi y}{\pi y} . 
\eeq
Since we are assuming $x > 1$, the above result vanishes when  
$\delta \to 0$. Therefore from \eqref{sum_int2} 
and \eqref{int_c1_c3} we conclude that 
\beq\label{zeta_anal1}
\zeta(z) = \dfrac{\pi \, \CI(z)}{ \Gamma(z)\sin\(\pi z\)}.
\eeq 
Using the well known identity
\beq\label{gamma_ident1}
\Gamma(z)\Gamma(1-z) = \dfrac{\pi}{\sin\(\pi z\)}
\eeq
we can further simplify \eqref{zeta_anal1} obtaining
\beq\label{zeta_analytic}
\zeta(z) = \Gamma(1-z)\, \CI(z).
\eeq

Although in obtaining \eqref{zeta_analytic} we assumed $\Re(z) > 1$, the
function $\Gamma(z)$ and the integral \eqref{I_def} remains valid 
for every complex $z$. Thus we can define $\zeta(z)$ through 
\eqref{zeta_analytic} on the whole complex plane (except at $z=1$ as will
be shown below). Since \eqref{zeta_analytic}
and \eqref{zeta_series} are the same on the half plane $\Re(z) > 1$, then
\eqref{zeta_analytic} yields the analytic continuation of the $\zeta$-function
to  the entire complex plane. The function defined by 
\eqref{zeta_analytic} is known as the \emph{Riemann $\zeta$-function}.

The integral \eqref{I_def} is an entire function of $z$, and
$\Gamma(1-z)$ has simple poles at $z=1,2,\dots$. Since $\zeta(z)$ does
not have zeros for $\Re(z) > 1$ then $\CI(z)$ must vanish at
$z=2,3,\dots$. In fact it is easy to see this explicitly. Note that
for $z= n \in \mathbb{Z}$ the integral \eqref{I_def} over $\CC_1$ and
$\CC_3$ cancel each other, thus it remains the integral over $\CC_2$ only.
From Cauchy's integral formula we have
\beq\label{res_in}
\CI(n) = \lim_{u \to 0} \dfrac{u^n}{e^{-u}-1} = 0 \qquad (n=2,3,\dotsc).
\eeq
Moreover, since $\CI(n) = 0$ for $n\ge2$, cancelling the poles of 
$\Gamma(1-n)$, the only possible pole of 
$\zeta(z)$ occurs at $z=1$. From \eqref{res_in} we have $\CI(1) = -1$, 
showing that $\zeta(z)$ indeed has a \emph{simple pole} inherited from 
the simple pole of $\Gamma(0)$. 
Since this pole is simple, we can compute the residue of $\zeta(1)$ through
\beq
\mbox{Res} \, \zeta(1) = 
\lim_{z \to 1}\(z-1\)\zeta(z) =
\lim_{z \to 1}\(z-1\)\Gamma(1-z)\CI(z) =
- \Gamma(1)\CI(1) = 1,
\eeq
where we have used \eqref{zeta_analytic} and the well known 
property $\Gamma(z+1)=z\, \Gamma(z)$.

\subsection{Functional equation}

The functional equation was already motivated by physical arguments 
in \eqref{ap.8}. Now we proceed to derive it mathematically. There
are at least $7$ different ways to do this \cite{Titchmarsh}. 
We are going to present only one way, which we think is the simplest.

\begin{figure}[h]
\centering
\includegraphics{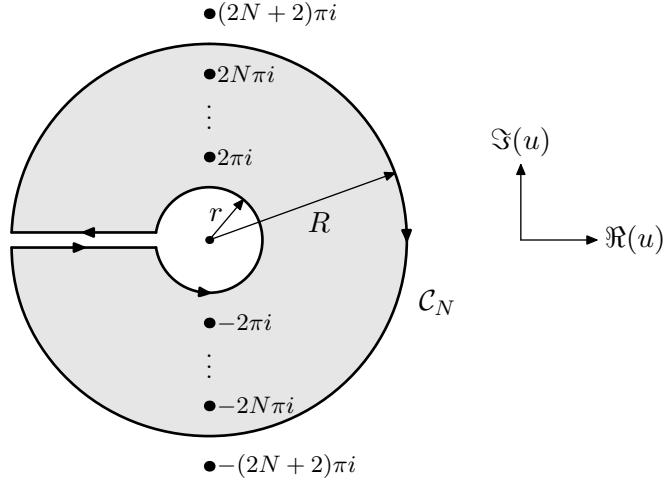}
\caption{Contour $\CC_N$ where $r < 2\pi$ and $2N \pi < R < (2N+2)\pi$.
This contour encloses $2N$ poles.}
\label{fig_path_cn} 
\end{figure}

Let us now consider the following integral
\beq\label{In_def}
\CI_N(z) = \dfrac{1}{2\pi i}\oint_{\CC_N} \dfrac{u^{z-1}}{e^{-u}-1}du
\eeq
with the path of integration $\CC_N$ as shown in Figure~\ref{fig_path_cn}.
Note that the integrand has the poles $u_n=2\pi n i $ ($n\in \mathbb{Z}$) 
along the imaginary axis.
On this contour $r < 2 \pi$ and
$2\pi N i < R < (2N+2)\pi i$, thus the annular 
region encloses $2N$ poles corresponding to $n=\pm1,\pm2,\dotsc,\pm N$.

Let us first consider the integral \eqref{In_def} over the outer circle with
radius $R$, i.e. $u=Re^{i\theta}$. Let $g(u) = \(e^{-u}-1\)^{-1}$. The
function $g(u)$ is meromorphic, i.e. has only isolated singularities
at $u_n=2\pi  n i$. Around the circle of radius $R$ there are no singularities
thus $g(u)$ is bounded in this region, i.e. $|g(u)| \le M$.
We also have 
$\left| u^{z-1} \right| = R^{x-1} e^{-\theta y} \le R^{x-1} e^{|y|\pi}$.
Therefore
\beq
\left|\int_{R} \dfrac{u^{z-1}}{e^{-u}-1}du\right| \le
\int_{-\pi}^{\pi} M R^{x-1} e^{|y|\pi} R d \theta = 2\pi M R^{x} e^{|y|\pi}.
\eeq
Now if $N\to \infty$ then $R \to \infty$ and the above result goes
to zero if $x < 0$. In this way the only contribution from the integral 
\eqref{In_def} is due to the smaller circle $r$ and this is equal to
the integral \eqref{I_def}. Thus we conclude that 
\beq\label{lim_in}
\lim_{N\to\infty} \CI_{N}(z) = \CI(z) \qquad (\Re(z) < 0).
\eeq

We can now evaluate the integral \eqref{In_def} through the residue theorem.
The integrand is $f(u) = u^{z-1}/\(e^{-u}-1\)$, then we have
\beq
\CI_N(z) = - \sum_{n=1}^{N+1} \( \mbox{Res} f\(u_n\) + 
\mbox{Res} f\(u_{-n}\) \),
\eeq
where the minus sign is because the contour is traversed in the clockwise 
direction.
The residue can be computed as follows:
\beq
\mbox{Res}f\(u_n\) = \lim_{u\to u_n}\(u-u_n\)f\(u\) = 
u_n^{z-1}\lim_{u\to u_n }\dfrac{u-u_n}{e^{-u}-1} = -(2\pi n i )^{z-1}.
\eeq
Therefore,
\beq\label{sum_in}
\CI_N(z) = \sum_{n=1}^{N+1} (2\pi n)^{z-1} \( i^{z-1}+(-i)^{z-1} \)
= 2(2\pi)^{z-1}\cos\(\pi(z-1)/2\) \sum_{n=1}^{N+1}n^{z-1}.
\eeq
Now the sum in \eqref{sum_in} becomes particularly interesting if we
replace $z\to 1-z$, yielding the same term $n^{-z}$ that appears in the 
series of the $\zeta$-function. In view of \eqref{lim_in}, now valid 
for $\Re(z) > 1$, we thus have
\beq
\lim_{N\to\infty}\CI_N(1-z) = \CI(1-z) = 
2 (2\pi)^{-z}\cos\(\pi z/2\)\zeta(z).
\eeq
Using \eqref{zeta_analytic} we obtain
\beq\label{zeta_rel}
\zeta(1-z) = 2 (2\pi)^{-z} \cos\(\pi z/2\) \Gamma(z) \zeta(z),
\eeq
which relates $\zeta(z)$ to $\zeta(1-z)$. This equality can be written
in a much more symmetric form using the following two well known properties
of the $\Gamma$-function:
\begin{align}
\Gamma(z)\Gamma(1-z) &= \dfrac{\pi}{\sin(\pi z)}, 
\label{gamma_reflection} \\
\Gamma(z)\Gamma(z+1/2) &= 2^{1-2z} \sqrt{\pi}\, \Gamma(2z). 
\label{gamma_dupplication}
\end{align}
Replacing $z\to (z+1)/2$ in \eqref{gamma_reflection} and substituting into
\eqref{zeta_rel} we have
\beq
\zeta(1-z) = \dfrac{2 (2\pi)^{-z}\pi \, \Gamma(z) \,\zeta(z)}{
\Gamma\((z+1)/2\)\Gamma\((1-z)/2\)}.
\eeq
Now replacing $z\to z/2$ in \eqref{gamma_dupplication} and substituting
the term $\Gamma(z) / \Gamma\((z+1)/2\)$ into the above equation we
obtain
\beq\label{func_equation}
\chi(z) = \chi(1-z), \qquad 
\chi(z) \equiv \pi^{-z/2} \, \Gamma\(z/2\) \zeta(z).
\eeq
This equality is known as the \emph{functional equation} for 
the $\zeta$-function. This is an amazing  relation, first discovered by
Riemann.
In deriving \eqref{func_equation} we had to assume $\Re(z) > 1$, but 
through analytic continuation it is valid on the whole complex plane, 
except at $z=1$ where $\zeta(z)$ has a simple pole. Note also that
replacing $z\to 1+2n$ for $n=1,2,\dotsc$ into \eqref{zeta_rel} we see that
$\zeta(-2n) = 0$, corresponding to the zeros of $\cos\(\pi/2 + n \pi\)$.

\subsection{Trivial zeros and specific values}

We already have seen that due to the Euler product formula \eqref{euler_prod}
$\zeta(z)$ has  no zeros for $\Re(z) > 1$. We have also seen in connection
with \eqref{zeta_analytic} that $\zeta(z)$ has a simple pole at $z=1$.
It follows from this pole that there are an infinite number of prime numbers.   
Moreover, due to \eqref{zeta_rel} we have $\zeta(-2n)=0$ for $n=1,2,\dotsc$.
These are the so called \emph{trivial zeros}. Since there are
no zeros for $\Re(z) > 1$, the functional equation
\eqref{func_equation} implies that there are no other zeros for $\Re(z) < 0$.

The other possible zeros of $\zeta(z)$ must therefore be inside the 
so-called \emph{critical strip}, $0 \le \Re(z) \le 0$. These are called the
\emph{non-trivial zeros}, since they can be complex contrary to the 
trivial ones. Note that from \eqref{func_equation}, since
$\Gamma(z/2)$ has no zeros on the critical strip, then $\zeta(z)$ and
$\chi(z)$ have the same zeros on this region. Moreover, these zeros
are symmetric between the line $\Re(z) = 1/2$. Since 
$\(\chi(z)\)^* = \chi(z^*)$ if $\rho$ is a zero so is its complex conjugate
$\rho^*$. Thus zeros occur in a quadruple: $\rho$, $\rho^*$, 
$1-\rho$ and $1-\rho^*$. The exception are for zeros on the so
called \emph{critical line} $\Re(z) = 1/2$, where $\rho$ and $1-\rho^*$
coincide.  It is known that there are an infinite number of zeros on 
the critical line. 
It remains unknown whether they are all simple zeros.   
We will discuss these non-trivial zeros in more detail later.

Now let us consider special values of the $\zeta$-function. Let us start
by considering the negative integers $\zeta(-n)$. From
\eqref{zeta_analytic} we have $\zeta(-n) = n! \, \CI(-n)$. The integral
\eqref{I_def} at this point can be computed through the residue
theorem. Since the only pole is at $z=0$ we have
\beq\label{zeta_residue}
\zeta(-n) = n! \, 
\mbox{Res}\left.\( \dfrac{u^{-n-1}}{e^{-u}-1} \)\right|_{u=0}.
\eeq
The Bernoulli polynomials $B_n(x)$ are 
defined through the generating function
\beq\label{bernoulli_pol}
\dfrac{u e^{x u}}{e^u-1} = \sum_{n=0}^{\infty} \dfrac{B_n(x)}{n!} u^{n}
\qquad (|u| < 2\pi).
\eeq
The values $B_n\equiv B_n(0)$ are called Bernoulli numbers and defined
by setting $x=0$ in \eqref{bernoulli_pol}. Then using \eqref{bernoulli_pol}
into \eqref{zeta_residue} we have
\beq\label{zeta_bernoulli1}
\zeta(-n) = n! \, \mbox{Res}\left.
\( -u^{-n-2}\sum_{m=0}^{\infty} \dfrac{B_m(1)}{m!} u^m \)\right|_{u=0}
= - \dfrac{B_{n+1}(1)}{n+1}.
\eeq
From the well known relation $B_n(x+1) -B_n(x) = nx^{n-1}$ for $n\ge 1$
we have $B_n(1) = B_n(0)$ for
$n\ge 2$. Thus we have
\beq\label{zeta_bernoulli}
\zeta(-n) = - \dfrac{B_{n+1}}{n+1} \qquad (n\ge 1).
\eeq
The case $\zeta(0)$ can be obtained from \eqref{zeta_bernoulli1} since
$B_1(1) = 1/2$, then $\zeta(0) = -1/2$. Also, since $B_{2n+1}=0$ we see
once again that $\zeta(-2n)=0$ for $n\ge1$.

Setting $z=2k$ in \eqref{zeta_rel} 
we have $\zeta(1-2k) = 2 (2\pi)^{-2k}\cos\(\pi k\)\Gamma(2k)\zeta(2k)$,
but from \eqref{zeta_bernoulli} we have $\zeta(1-2k) = -B_{2k}/2k$,
and thus we have $\zeta(z)$ over the positive even integers,
\beq\label{zeta_bernoulli2}
\zeta(2k) = (-1)^{k+1} \dfrac{(2\pi)^{2k} B_{2k}}{2 (2k)!}
\qquad (k\ge1).
\eeq
The previous results \eqref{zeta_of_2} and \eqref{zeta_of_4} are particular
cases of the above formula.

A very interesting question concerns  the values of $\zeta(2k+1)$.
Setting $z=2k+1$ into \eqref{zeta_rel} we get zero on both sides, and for
$z=-2k$ the poles of $\Gamma(-2k)$ cancel with the zeros of $\zeta(-2k)$
but this product is still undetermined, so we get no information about
$\zeta(2k+1)$. No simple formula analogous to \eqref{zeta_bernoulli2} is
known for these values. It is known that $\zeta(3)$ is irrational, this
number is called \emph{Ap\' ery's constant} 
\cite{Apery1,Apery2,Beukers,Zudilin}. It is also known that
there are infinite numbers in the form $\zeta(2k+1)$ which are 
irrational \cite{Rivoal}. Due to the unexpected nature of the result,
when Ap\' ery first showed his proof many
mathematicians considered it as flawed, however, H. Cohen, H. Lenstra
and A. van der Poorten confirmed that in fact Ap\' ery was correct.
It has been conjectured that $\zeta(2k+1)/ \pi^{2k+1}$ is 
transcendental \cite{Kohnen}. A very interesting physical connection, 
providing a link between number theory and statistical mechanics, is that
the most fundamental correlation function of the $XXX$ spin-$1/2$ chain 
can be expressed in terms of $\zeta(2k+1)$ \cite{Korepin}.

\section{Gauss and the prime number theorem}
\label{sec:Gauss}

\begin{figure}
\centering
\begin{minipage}{.49\textwidth}
\centering
\includegraphics[width=1\linewidth]{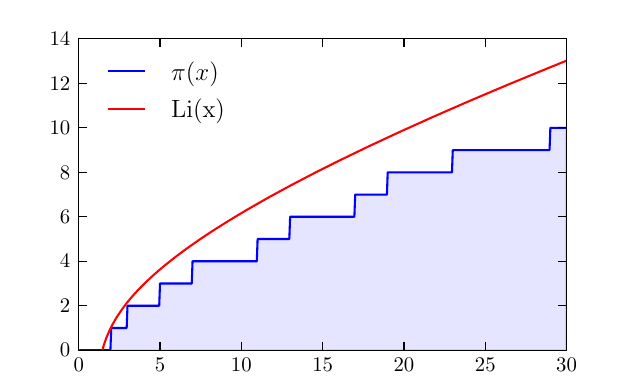}
\end{minipage}
\begin{minipage}{.49\textwidth}
\centering
\includegraphics[width=1\linewidth]{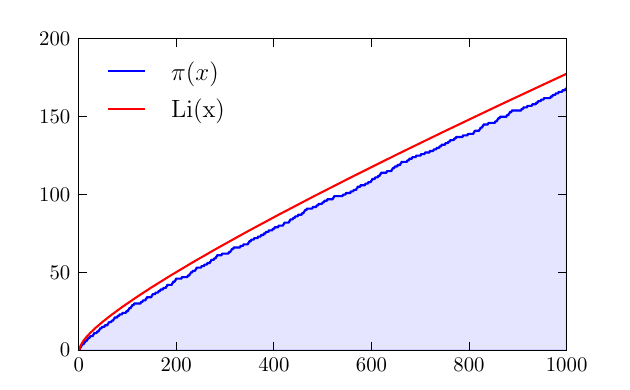}
\end{minipage}
\caption{
The $\Li (x)$ approximation to the prime number counting 
function $\pi (x)$.}
\label{fig_LiPi1}
\end{figure}

Let $\pi (x)$ denote the number of primes less than the positive 
real number $x$. 
It is a staircase function that jumps by one at each prime.      
In 1792,  Gauss, when only
15 years old,  based on examining data on the known primes,  guessed that 
their density 
went as $1/\log(x)$.     This leads to the approximation
\beq
\label{Gauss1}  
\pi (x)  =  \sum_{p\leq x}  1  \approx  \Li (x) 
\eeq
where $\Li (x)$ is the log-integral function,
\beq
\label{Li}
\Li (x)  =  \int_0^x  \frac{dt}{\log t} .
\eeq
$\Li (x)$ is a smooth function,  and does indeed provide a smooth 
approximation
to $\pi (x)$ as Figure~\ref{fig_LiPi1}  shows.     

The \emph{Prime number theorem} (PNT) is the statement 
that $\Li (x)$ is the leading approximation to
$\pi (x)$.   It was only proven 100 years later 
using the main result of Riemann described in the next section.   
As we will explain later,   
the PNT follows if there are no Riemann zeros with $\Re(z) = 1$.

The key ingredient in Riemann's derivation of his result is the Euler 
product formula relating
$\zeta (z)$ to the prime numbers.    
A simple derivation of Euler's formula is based to the 
ancient ``sieve'' method for locating primes.   
One begins with a list of integers.   First one removes
all even integers,  then all multiples of $3$, then all multiples of $5$,  
and so on.   Eventually one ends up
with the primes.     
We can describe this procedure analytically as follows.  
Begin with 
\beq
\zeta (z)  =  1  +  \inv{2^z}  +  \inv{3^z} + \inv{4^z} + \dotsm .
\eeq
One has 
\beq
\inv{2^z} \zeta (z) =  \inv{2^z}  + \inv{4^z}  + \inv{6^z}  + \dotsm
\eeq
thus
\beq
\( 1- \inv{2^z} \) \zeta (z) =  1 + \inv{3^z}  +  \inv{5^z}  + \dotsm. 
\eeq
Repeating this process with powers of $3$ we have
\beq
\(1- \inv{3^z} \)  \(1- \inv{2^z}\)  \zeta (z)  
=   1 + \inv{5^z}  +  \inv{7^z} +  \dotsm. 
\eeq
Continuing this process to infinity,  the right hand side equals $1$.   
Thus 
\beq
\label{EulerProd}
\zeta (z)  =  \prod_{p}  \inv{1-p^{-z}} .
\eeq

Chebyshev tried to prove the PNT using $\zeta (z)$ in 1850.  
It was finally proven in 1896  by Hadamard and de la Vall\'e  Poussin by
demonstrating that $\zeta(z)$ indeed has no zeros with $\Re (z) =1$.

\section{Riemann zeros and the primes}

\subsection{Riemann's main result}

Riemann obtained an explicit and exact expression for $\pi (x)$ in terms of 
the non-trivial zeros  $\rho$ of $\zeta(z)$.  
There are simpler but equivalent versions 
of the main result,  based  on the function $\psi (x) $ below.   
However, let us present the main formula for $\pi (x)$ itself,
since it is historically more important.   
The derivation is given in the next subsection. 

The function $\pi (x)$ is related to another number-theoretic 
function $J(x)$, defined as 
\beq
\label{Jx}
J(x)  =  \sum_{2\leq n \leq x}    \frac{\Lambda (n)}{\log n} 
\eeq
where  $\Lambda (n)$,  the  von Mangoldt function, is defined by
\beq
\Lambda(n) = \begin{cases} 
\log p & \mbox{if $n=p^m$ for some prime $p$ and integer $m\ge 1$,} \\
0 & \mbox{otherwise.}
\end{cases}
\eeq
For instance $\Lambda (3)= \Lambda (9) = \log 3$.  
The two functions 
$\pi (x)$ and $J(x)$ are related by M\" obius inversion as follows:
\beq
\label{mob1}
\pi (x) = \sum_{n\geq 1} \frac{\mu (n)}{n}  J(x^{1/n}).
\eeq
Here  $\mu (n)$ is the M\" obius function defined as follows.
For $n>1$, through the prime decomposition theorem we can write 
$n=p_1^{\alpha_1}\dotsm p_k^{\alpha_k}$. Then
\beq
\mu(n) = \begin{cases}
(-1)^k & \mbox{if $\alpha_1=\alpha_2=\dotsm=\alpha_k=1$,} \\
0 & \mbox{otherwise.}
\end{cases}
\eeq
We also have $\mu(1) = 1$. Note that $\mu(n)=0$ if and only 
if $n$ has a square
factor $> 1$.
The above expression \eqref{mob1} is actually a finite sum,  
since for large enough $n$,  $x^{1/n} <2$ and $J=0$.

\begin{figure}
\centering
\begin{minipage}{.49\textwidth}
\centering
\includegraphics[width=1\linewidth]{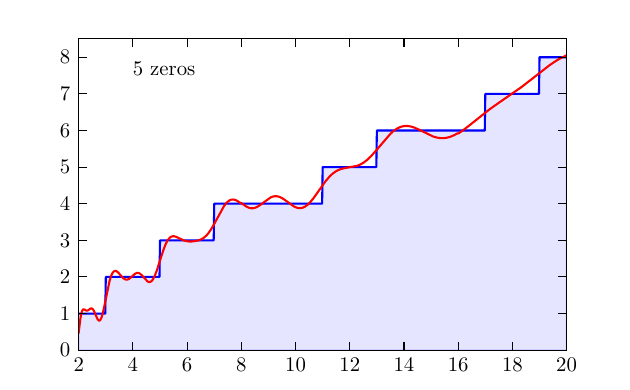}
\includegraphics[width=1\linewidth]{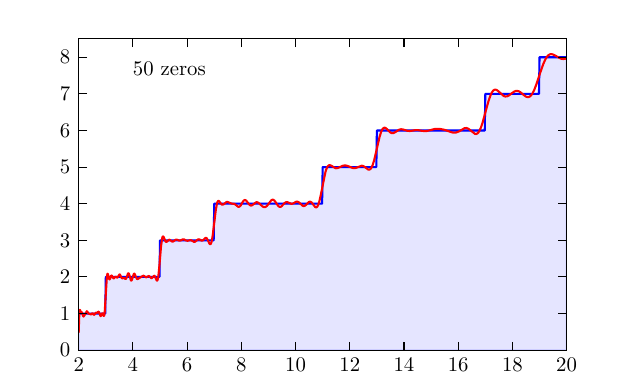}
\end{minipage}
\begin{minipage}{.49\textwidth}
\centering
\includegraphics[width=1\linewidth]{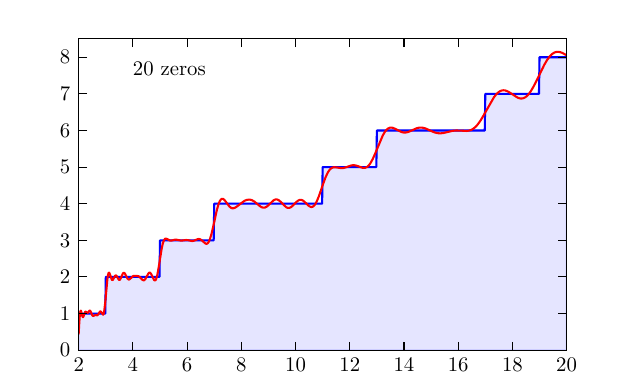}
\includegraphics[width=1\linewidth]{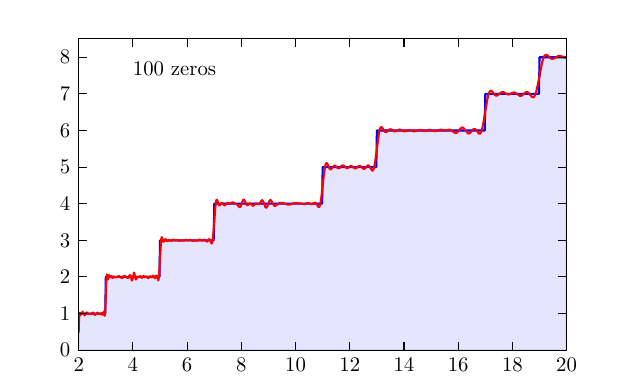}
\end{minipage}
\caption{The function $\pi (x)$ as a sum of Riemann zeros from equations  
\eqref{mob1} and \eqref{Jzeros}. The dashed (blue) region represents the
exact $\pi(x)$ while the (red) oscillating curve is \eqref{mob1}, obtained
with different number of zeros as indicated in the figures.}
\label{fig_Pizeros}
\end{figure}

The main result of Riemann is a formula for $J(x)$, expressed 
as an infinite sum over non-trivial zeros $\rho$,
\beq
\label{Jzeros}
J(x) =  \Li (x) - \sum_\rho  \Li\(x^\rho\)  +  
\int_x^\infty  \dfrac{dt}{\log t} ~  \inv{ \(t^2 -1\) t} - \log 2.
\eeq
Riemann derived the result \eqref{Jzeros} starting from the Euler product formula 
and utilizing some insightful complex analysis that was sophisticated for the time.   
Some care must be taken in  numerically 
evaluating $\Li (x^\rho )$  since it has a branch point. 
It is more properly defined through the exponential integral function
\beq
\Li(x) = {\rm Ei} (\rho \log x ), \qquad 
{\rm Ei} (z) = - \int_{-z}^\infty  dt \, \dfrac{e^{-t}}{t}.
\eeq
The sum in \eqref{Jzeros} is real because the $\rho$'s come 
in conjugate pairs.
If there are no zeros on the line $\Re(z) = 1$, then the dominant 
term is the first one, i.e. $J(x) \sim  \Li (x)$,  
and this proves the PNT.
The sum over $\rho$ corrections to $\Li(x)$ 
deform it to the staircase function $\pi(x)$ as Figure \ref{fig_Pizeros}
shows.    
Thus, the complete knowledge of the primes is contained in the Riemann zeros.

Von Mangoldt  provided a simpler formulation based on the function 
\beq
\label{PsiDef}
\psi (x)  =   \sum_{n \leq x}   \Lambda (n).
\eeq
The  function $\psi (x)$ has a simpler expression in terms of Riemann zeros
which reads
\beq
\label{psizeros}
\psi (x)  =   
x - \sum_\rho\dfrac{x^\rho}{\rho} - \log (2\pi)  - \inv{2} 
\log \( 1 - \inv{x^2} \).
\eeq
In this formulation,  the PNT follows from the fact 
that the leading term is $\psi (x) \sim x$.

\subsection{\texorpdfstring{$\psi (x)$}{Psi} and the Riemann zeros}

We first derive the formula \eqref{psizeros}. 
From the Euler product formula one has 
\beq
\d_z \log \zeta (z) =  -  \sum_p \d_z  \log\(1-p^{-z}\)  
= - \sum_p  \log p \,  \dfrac{p^{-z}}{1-p^{-z}}.
\eeq
Taylor expanding the factor $1/(1-p^{-z})$ one obtains 
\beq
\label{dlogZ}
\d_z \log \zeta (z)  =  -  \sum_p  \sum_{m=1}^\infty  \frac{\log p}{p^{mz}}.
\eeq
For any arithmetic function $a(n)$,   the Perron formula relates  
\beq\label{Aa}
A(x) = \sum_{n \leq x}{}^{'} a(n)
\eeq
to the poles of the Dirichlet series 
\beq
g(z) =  \sum_{n=1}^\infty  \frac{a(n)}{n^z}.
\eeq
In \eqref{Aa} the restriction on the sum is such that if $x$ is an integer,  
then the last term of the sum must be multiplied by $1/2$.
Now  $\zeta (z)$ can be factored
in terms of its zeros,  $\zeta (z)  \propto  \prod_\rho (z- \rho)$,   thus  
$\d_z \log \zeta (z)$ has poles at each zero $\rho$.    This implies that the 
Perron formula can be used to relate $\psi (x)$ to the Riemann zeros.   

The Perron formula is essentially an inverse Mellin transform.    
If the series for $g(z)$ converges 
for $\Re (z) > z_1$,  then 
\beq
\label{Perron}
A(x) =  \inv{2 \pi i}  \int_{c - i\infty}^{c+ i \infty}  \dfrac{dz}{z}  
\, g(z) \, x^z 
\eeq
where the $z$-contour of integration is a straight vertical 
line from $-\infty$ to $+\infty$ 
with $c> z_1$.   
For completeness,  we present a derivation of this formula in 
Appendix~\ref{sec:perron}.   

Let us apply the Perron formula to $\psi (x)$,
\beq
\label{PerronPsi}
\psi (x) = \inv{2 \pi i}  \int_{c - i\infty}^{c+ i \infty}    
\dfrac{dz}{z}  ~g(z) \,   x^z  
\eeq
where $g(z) = - \d_z \log \zeta (z) $ and $c>1$.    The line of integration can be 
made into a closed contour by closing at infinity with $\Re(z) \leq c$.    
Now
\beq
g(z) =  - \d_z \(  \sum_\rho \log (z-\rho)  + \sum_{\rho'} \log (z - \rho')  -  \log(z-1) \)  
\eeq
where $\rho$ are  zeros of $\zeta(z)$ on the critical strip 
and $\rho'$ are the trivial zeros
on the negative real axis at $\rho' = -2n$.   
The $- \log (z-1)$ is due to the pole at $z=1$.
The sum of the residues gives 
\beq
\psi (x) = x  - \sum_{\rho}  \dfrac{x^\rho}{\rho}  -  
\sum_{\rho'} \dfrac{x^{\rho'} }{\rho'}  + g(0) .
\eeq
The first term comes from the $z=1$ pole and $g(0) = -\log(2 \pi)$ 
comes from the 
$z=0$ pole.    
Finally
\beq
\sum_{\rho'} \frac{x^{\rho'} }{\rho'} = 
- \sum_{n=1}^\infty   \frac{x^{-2n}}{2n}  
 = \inv{2}  \log (1 - 1/x^2)
\eeq
and this gives the result \eqref{psizeros}.

\subsection{\texorpdfstring{$\pi (x)$}{Pi} and  the Riemann zeros} 

\begin{figure}
\centering
\includegraphics[width=0.6\linewidth]{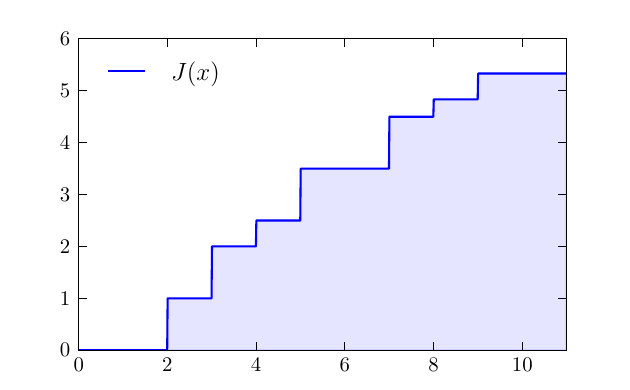}
\caption{The number theoretic function $J(x)$ for $x<11$. }
\label{fig_Jofx}
\end{figure}

Let us first explain the relation \eqref{mob1}   
between $\pi(x)$ and $J(x)$ involving the
M\"obius $\mu$ function.      By definition,  $J(x) =0 $ for $x<2$.   
It jumps by $1/n$ at each $x=p^n$ where $p$ is a prime.     
The expression \eqref{mob1} is always a finite sum since for $n$ large enough
$x^{1/n} <2$.    
Consider for instance the range $x\leq 10$.   $J(x)$ in this range is plotted
in Figure \ref{fig_Jofx}.   

Since $10^{1/4} < 2$,  
the formula \eqref{mob1}  gives 
\beq
\pi(x) =  J(x)  - \inv{2} J(x^{1/2}) - \inv{3}  J(x^{1/3}) .
\eeq
One easily sees that the two subtractions remove from $J(x)$ the 
jumps by $1/2$ at 
$x= 2^2, 3^2$ and the jump by $1/3$ at $x=2^3$ leaving only 
the jumps by one at the primes $2,3,5,7$.   

Let us now derive \eqref{Jzeros}.   Comparing definitions,  one has 
\beq
dJ(x) = \inv{\log x}  d \psi (x) .
\eeq
Integrating this
\beq
J(x) =  \int_0^x  \dfrac{dt}{\log t}   \dfrac{d \psi (t)}{dt}  = \int_0^x   
\dfrac{dt}{\log t}    
\(  1 - \sum_\rho  t^{\rho -1}  - \inv{t (t^2 -1)}  \).
\eeq
Making the change of variables $y = t^\rho$,  
\beq
\int_0^x  \dfrac{dt}{\log t}  t^{\rho -1}  =  
\int_0^{x^\rho}  \dfrac{dy}{\log y}  =  \Li (x^\rho) .
\eeq
Finally using 
\beq
\int_0^x  \dfrac{dt}{\log t}  \inv{(t^2 -1) t}  +  
\int_x^\infty \dfrac{dt}{\log t}  \inv{(t^2 -1) t}  =  \log 2 
\eeq
one obtains the form \eqref{Jzeros}.

\section{An electrostatic analogy}   
\label{Electric}

A complex function is difficult to visualize
since it is a hypersurface in a $4$-dimensional space. 
In this section we construct an electric field and electric potential 
and use them to visualize the RH through a single real scalar
field over the $2$-dimensional $(x,y)$-plane, where $z=x+ i y$.

\subsection{The electric field}  

\begin{figure}
\centering
\includegraphics{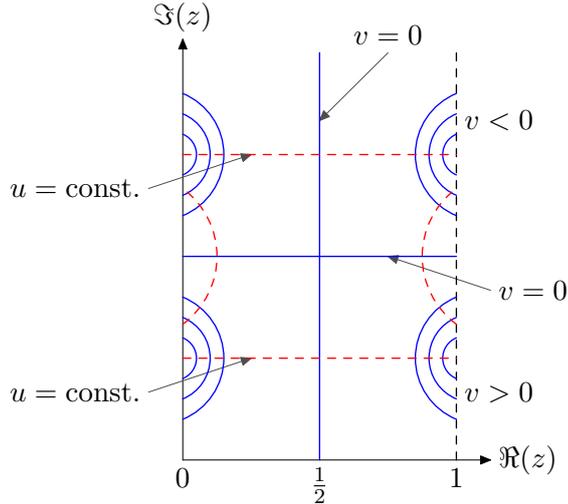}
\caption{$u$ contours are dashed (red) lines and $v$ contours are 
solid (blue) lines.
A Riemann zero on  occurs where a $u=0$ contour spans the entires strip and
passes through the zero on the critical line which is a $v=0$ contour.}
\label{uvcontours} 
\end{figure}

Let us remove the $z=1$ pole in $\chi (z)$ while maintaining its 
symmetry under
$z \to 1-z$ by defining the function
\beq
\label{xidef}  
\xi (z)  \equiv  \tfrac{1}{2} z (z-1) \chi (z) =  \tfrac{1}{2} 
z (z-1) \pi^{-z/2}  \Gamma (z/2) \zeta (z) 
\eeq
which satisfies 
\beq
\label{xisym}
\xi (z) = \xi (1-z).
\eeq    
Let us define the real and imaginary parts of $\xi(z)$ as 
\beq
\label{realimag}
\xi (z) =  u(x,y) + i \, v(x,y)
\eeq
The Cauchy-Riemann equations  
\beq
\d_x u = \d_y v, \qquad \d_y u = - \d_x v
\eeq
are satisfied everywhere since $\xi$ is an entire function.
Consequently,   both $u$ and $v$ are harmonic functions,   i.e.  
solutions of the Laplace equation
$\gradvec^2 u  =  (\d_x^2 + \d_y^2 ) u = 0$ and  $\gradvec^2 v = 0$,  
although they are not completely independent.    
Let us define $u$  or $v$  contours as the curves in the $x$-$y$ 
plane corresponding to $u$ or $v$ equal to
a constant,  respectively. 
The critical line is a $v=0$ contour since $\xi$ is real along it.
As a consequence of the Cauchy-Riemann equations we have
\beq
\gradvec u \cdot \gradvec v = 0.
\eeq
Thus  where the $u$ and $v$ contours intersect,  they are  
necessarily perpendicular,  and this is one aspect
of their dependency.     
A Riemann zero occurs wherever the  $u=0$ and $v=0$ contours intersect,
as illustrated in Figure~\ref{uvcontours}.

From  the symmetry \eqref{xisym}
and $\xi (z)^* = \xi (z^*)$  it follows that 
\beq
\label{symmetries}
u(x,y) = u (1-x, y), \qquad v(x,y) = - v (1-x,y) .
\eeq   
This implies that the $v$ contours do not cross the critical line 
except for $v=0$.    
All  the $u$ contours on the other hand are allowed to cross it by 
the above symmetry.    
Away from the $v=0$ points on the line $\Re (z) = 1$,   since the $u$ 
and $v$ contours are perpendicular,  
the $u$ contours generally cross the critical line and span the whole 
strip due to the symmetry \eqref{symmetries}.     
The $u$ contours that do not cross the critical line must be in 
the vicinity of the $v=0$ contours,
again by the perpendicularity of their intersections.
Figure \ref{uvcontours}   depicts the behavior of the $u$ and $v$ 
contours in regions of the critical strip with no zeros off of the line.

Introduce the vector field 
\beq
\label{Evec}
\vec{E} =   E_x  \, \hat{x}  +  E_y  \, \hat{y}  \equiv   
u (x,y) \,  \hat{x}   -   v(x,y) \,   \hat{y}
\eeq
where $\hat{x}$ and $\hat{y}$ are unit vectors in the $x$ 
and $y$ directions.
This field has zero divergence and curl as a consequence of the 
Cauchy-Riemann equations,
\beq
\label{curl}
\gradvec \cdot \Evec =0 ,\qquad \gradvec \times \Evec = 0, 
\eeq
which are defined everywhere since $\xi$ is entire.   
Thus it satisfies the conditions of a static electric field with 
no charged sources.
We are only interested in the electric field on the critical strip.     
$\Evec$ is not 
a physically realized electric field here,  in that  we do not need to 
specify what kind of charge distribution 
would give rise to such a field.     
All of our  subsequent arguments will be based only on the mathematical 
identities expressed in 
equation \eqref{curl},  and our reference to electrostatics is 
simply a useful analogy.
Since the divergence of $\Evec$ equals zero everywhere,  
the hypothetical electric charge distribution that gives
rise to $\Evec$ should be thought of as existing at infinity.    
Alternatively,  since $u$ and $v$ are harmonic
functions,   one can view them as being determined by their values on 
the boundary of the critical strip.     

\begin{figure}
\centering
\begin{minipage}{0.49\textwidth}
\includegraphics[width=.9\linewidth]{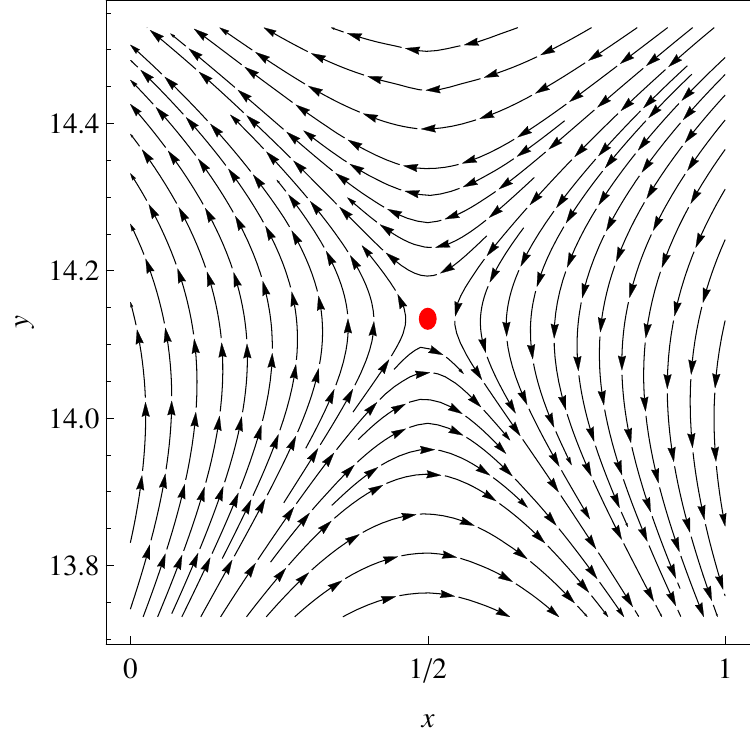}
\end{minipage}
\begin{minipage}{0.49\textwidth}
\vspace{-1em}
\includegraphics{figs/fig9.mps}
\end{minipage}
\caption{
\emph{Left: } 
field lines of $\Evec(x,y)$, equation \eqref{Evec},
in the vicinity of the first Riemann zero. \emph{ Right:} 
illustration of the field $\vec{E}(x,y)$  in the vicinity of two consecutive 
Riemann zeros $\rho_1$ and $\rho_2$  on the critical line.
}
\label{vectorplot} 
\end{figure}

As we now argue,  the main properties of the above $\Evec$  field on 
the critical strip are determined by its behavior near the 
Riemann zeros on the critical line combined with the behavior 
near $\Re(z) =1$.     In particular,  electric field lines do not cross. 
Any Riemann zero on the critical line arises from a $u=0$ contour
that crosses the full width of the strip and thus intersects the 
vertical $v=0$ contour.  
On the $u=0$ contour,  $E_x=0$,  whereas on the $v=0$ contour of
the critical line itself,  $E_y =0$.     Furthermore,  $E_y$ changes 
direction as one crosses the critical line.  
Finally,  taking  into account that $\Evec$ has zero curl,   
one can easily see that there are only two ways that 
all these conditions can be satisfied near the Riemann zero.    
One is shown in Figure~\ref{vectorplot} (left),   
the second has the direction of all arrows reversed.
In short,  Riemann zeros on the critical strip are  manifestly 
consistent with the necessary properties of $\Evec$.

We now turn to the global properties of $\Evec$ along the entire 
critical strip.     
The electric field must alternate in sign from one zero to the next,  
otherwise the curl of $\Evec$ would not be zero 
in a  region between two consecutive zeros.   
Thus there is a form of quasi-periodicity along the 
critical line,  in the sense that zeros alternate between being even 
and odd,  like the integers,  and also  analogous to
the zeros of $\sin x$ at $x= \pi n$  where $e^{i \pi n} = (-1)^n$.
Also,  along  the nearly horizontal $v=0$ contours that cross 
the critical line, 
$\Evec$ is in the $x$ direction.    
This leads to the pattern in Figure~\ref{vectorplot} (right).   
One aspect of  the rendition of this pattern is that it implicitly 
assumes that the $v=0$ and $u=0$ points along the line
$\Re (z) =1$ alternate,  namely,  between two consecutive $v=0$ points 
along this line,  there is only one
$u=0$ point,  which is consistent with the knowledge that there are 
no zeros of $\xi$ along the line $\Re (z) =1$.   
This fact will be clearer when we reformulate our argument in terms 
of the potential $\Phi$ below.

\subsection{The electric potential \texorpdfstring{$\Phi$}{}}

A mathematically integrated   version of the above arguments,  which 
has the advantage  of making manifest the dependency  of $u$ and $v$, 
can be formulated in terms of the  electric potential $\Phi$ which is a 
single real function,  defined to satisfy $\Evec = - \gradvec \Phi$.
Although it contains the same information as the above argument,  
it is more economical.
 
By virtue of $\gradvec \cdot \Evec =0$,  $\Phi$ is also a solution of 
Laplace's equation 
$\d_z \d_\zbar \Phi =0$ where we denote $\zbar =  z^*$.    
The general solution is  
that $\Phi$ is the sum of a function of $z$ and another function of $\zbar$.    Since $\Phi$ must be real, 
\beq
\label{pot}
\Evec =   - \gradvec \Phi, \qquad 
\Phi (x,y) =  \dfrac{1}{2}  \(   \varphi (z)   + \vphibar (\zbar) \) 
\eeq
where $\vphibar (\zbar) =  \(\varphi (z)\)^*$.    
Clearly $\Phi$ is not analytic,  whereas 
$\varphi$ is;   it is useful to work with $\Phi$ since we only have to 
deal with one real function.  
Comparing the definitions of $\Evec$ and $\xi$ in terms of $u$ and $v$,  
one finds 
\beq
u = -\dfrac{1}{2}\(\d_z \vphi + \d_\zbar \vphibar \),
\qquad
v= - \dfrac{i}{2} \( \d_\zbar \vphibar - \d_z \vphi\).
\eeq
This implies 
\beq
\label{Ed}   
\xi (z) = - \dfrac{\d \varphi (z)}{\d z  }
\eeq
This equation can be integrated because $\xi $ is entire.   
Riemann's original paper gave the following integral representation
\beq
\label{xiInt}
\xi (z)  = 4 \int_1^\infty dt  \,  G(t)  
\cosh \[(z-\smallhalf) \log (t) /2 \]  
\eeq
where
\beq
\label{Gtdef}
G(t)  =  t^{-1/4}  \d_t \(  t^{3/2} \d_t g \),  \qquad
g(t) =  \smallhalf \( \vartheta_3 (0, e^{-\pi t}) -1\)  =  
\sum_{n=1}^\infty e^{-n^2 \pi t}.
\eeq
Here,  $\vartheta_3$ is one of the four elliptic theta functions.
Note that the $z \to 1-z$ symmetry is manifest in this expression.   
Using this,  then up to 
an  irrelevant  additive constant  
\beq
\label{phi}
\varphi (z) = - 8 \int_1^\infty   
\dfrac{dt}{\log t}   \, G(t) \sinh \[ \smallhalf (z-\smallhalf) \log t \].
\eeq

Let us now consider the $\Phi = \mbox{const.}$ contours in the 
critical strip.   Using the integral 
representation \eqref{phi},  one finds the symmetry 
$\Phi (x,y) = -\Phi(1-x, y)$.   One sees then  that the $\Phi \neq 0$ contours 
do not cross the critical line,   whereas the $\Phi =0$ contours can and do.    Since $\vphi$ is imaginary along 
the critical line,   the latter is also a $\Phi=0$ contour.   

\begin{figure}
\centering
\begin{minipage}{.49\textwidth}
\includegraphics[width=.8\linewidth]{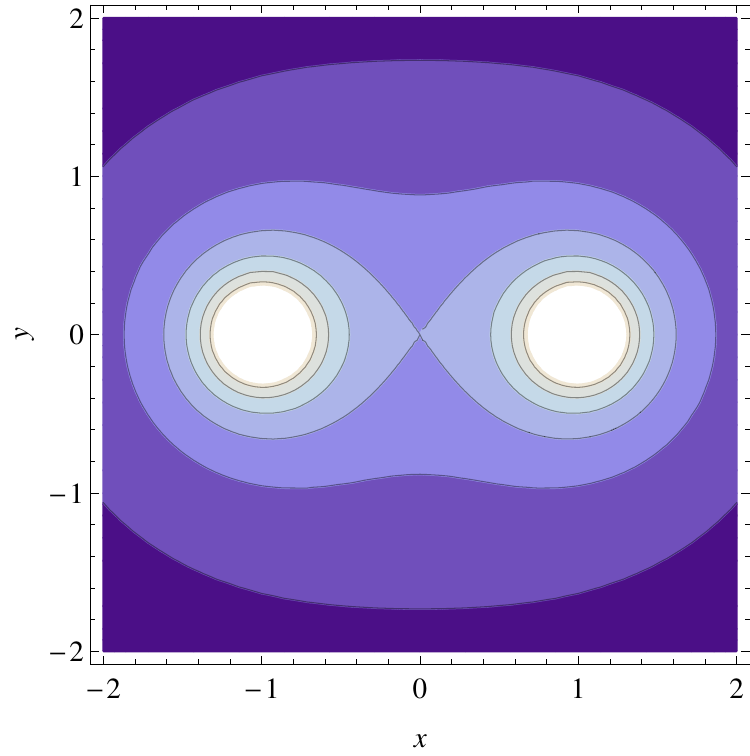}
\end{minipage}
\begin{minipage}{.49\textwidth}
\includegraphics[width=0.8\linewidth]{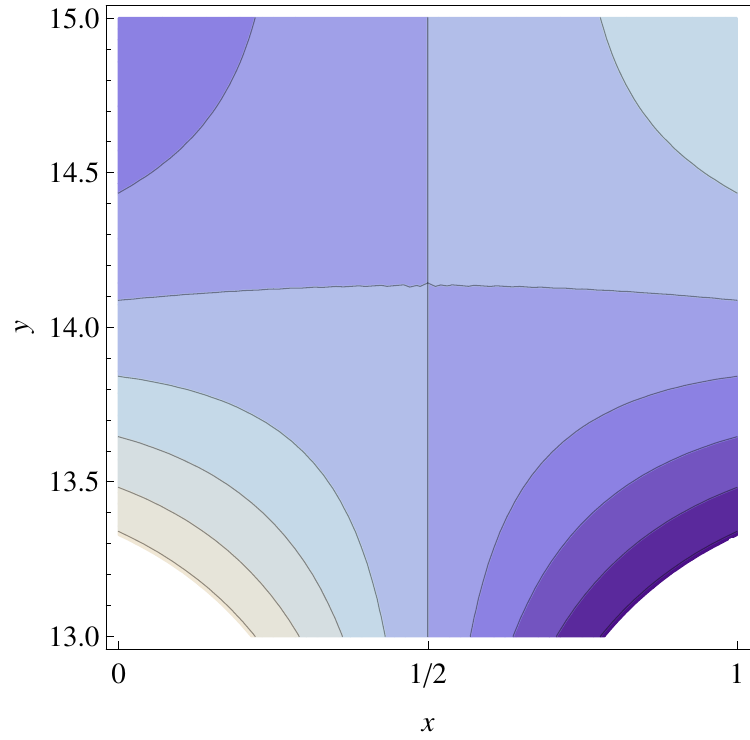} 
\end{minipage}
\caption{\emph{ Left:} the electric potential of two equal point charges.
\emph{Right:} contour plot of the potential $\Phi$  in the vicinity of the 
first Riemann zero  at $\rho=1/2 + (14.1347\dotsc )\, i$.  
The horizontal (vertical) direction is the $x$ ($y$) 
direction, where $z=x + i y$.
The critical line and nearly horizontal line are $\Phi=0$ contours 
and they intersect at the zero.}
\label{fig_PointCharge} 
\end{figure}

All  Riemann zeros  $\rho$ necessarily occur at isolated points,  
which is a property of entire functions.
This is clear from the factorization formula 
$\xi (z) = \xi (0)  \prod_\rho   ( 1- z/\rho )$,  conjectured by Riemann, 
and later proved  by Hadamard.    Where are these zeros located in 
terms of $\Phi$?   
At $\rho$,  $\gradvec \Phi =0$.     
Thus, such isolated zeros occur when two $\Phi$ contours 
\emph{intersect},   which 
can only occur if the two contours correspond to the same value 
of $\Phi$ since $\Phi$ is single-valued.    
A useful analogy is the electric potential for equal point charges.   
The electric field vanishes halfway between them,  and this is the 
unique point where the equi-potential contours vanish.   
The argument is simple:     $\gradvec \Phi$ is perpendicular to 
the $\Phi$ contours,  however as one approaches
$\rho$ along one contour,  one sees that it is not in the same 
direction as inferred from the approach from the other contour.  
The only way this could be consistent is if $\gradvec \Phi = 0$ at $\rho$.  
For purposes of illustration,  we show the electric potential 
contours for two equal point charges in
Figure~\ref{fig_PointCharge} (left).   Here,  the electric field is 
only zero halfway between the charges, 
and indeed this is where two $\Phi$ contours intersect.

With these properties of $\Phi$,  we can now begin to  understand the 
location of the known  Riemann zeros.   
Since the $\Phi=0$ contours intersect  the critical line,  
which itself is also a $\Phi=0$ contour,  a zero 
exists at each such intersection,  and we know there are an 
infinite number of them.   
The contour plot in Figure~\ref{fig_PointCharge} (right) for the actual 
function $\Phi$ constructed above verify these statements.   
We emphasize that there is nothing special about the value $\Phi =0$, 
since $\Phi$ can be shifted by an arbitrary constant without 
changing $\Evec$;  we defined it
such that the critical line corresponds to $\Phi =0$.

\begin{figure}
\begin{center}
\begin{minipage}{.49\textwidth}
\includegraphics{figs/fig12.mps}
\end{minipage}
\begin{minipage}{.49\textwidth}
\includegraphics[width=1\linewidth]{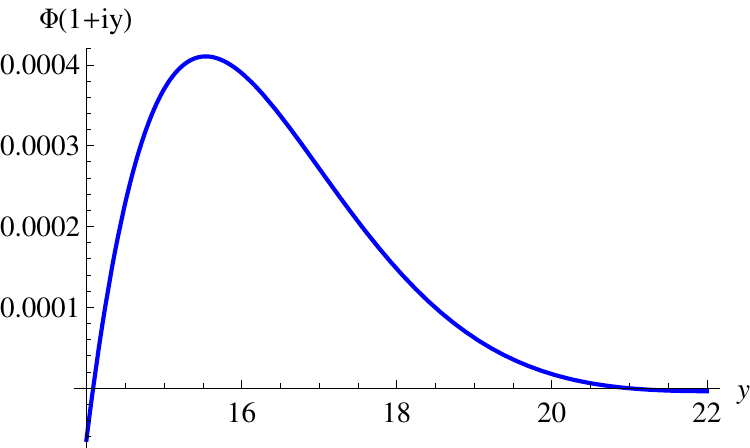}
\end{minipage}
\caption{\emph{Left:} a sketch of the contour plot of the potential 
$\Phi$  in the vicinity of  a \emph{hypothetical} Riemann zero off of 
the critical line.
Such a zero occurs where the contours intersect.  
$\rho_n$ and $\rho_{n+1}$ are consecutive zeros on the 
line.  
\emph{ Right:}
the electric potential between zeros on the boundary of the critical strip 
$\Re (z) =1 $. 
}
\label{fig_PotentialOff} 
\end{center}
\end{figure}

A  hypothetical  Riemann zero off of the critical line would then necessarily  correspond to an intersection of two 
$\Phi \neq 0$ contours.   
For simplicity,  let us assume that only two such contours intersect,  
since our arguments can be
easily extended to more of such intersections.     
Such a situation is depicted in Figure~\ref{fig_PotentialOff} (left). 

This figure   implies that on the line $\Re( z) =1$,  
specifically  $z=1+iy$,    $\Phi$  takes on  the same non-zero value at 
four different values  of $y$ between consecutive  zeros,  
i.e. roots  of the equation $ f(y) = 0$, where 
\beq
\label{fofy}
f(y) \equiv \Phi(1, y) = 
\Re (\vphi (1+iy)).
\eeq
Thus,  the  real function
$f(y)$ would have to have $3$ extrema between two consecutive zeros.
Figure~\ref{fig_PotentialOff} (right) suggests
that this does not occur.   In order to attempt to prove  it,  
let us define a ``regular alternating"  real function  $h(y)$  of a 
real variable 
$y$ as a function that alternates between positive and negative 
values  in the most regular manner possible:  between
two consecutive zeros $h(y)$ has only one maximum, or minimum.     
For example,  the  $\sin (y)$ function is obviously  regular  
alternating.    By the above argument,   if $f(y)$  is regular alternating,
then two $\Phi \neq 0$ contours cannot intersect and there are no Riemann 
zeros off the critical line.  
In Figure~\ref{fig_PotentialOff} (right)  we plot $f(y)$ for low values 
of $y$ in the vicinity of the
first two zeros,  and as expected,  it is regular alternating 
in this region.   

To summarize,   based on the symmetry \eqref{xisym},  and the 
existence of the known infinity 
of Riemann zeros along the critical line,   we have argued 
that $\Evec$ and $\Phi$ satisfy a regular repeating
pattern all along the critical strip,  and the RH would follow from 
such a repeating pattern.       In order to go further,   one obviously  
needs to investigate the  detailed properties of the function
$\xi$,  in particular its large $y$ asymptotic behavior,  and attempt 
to establish this  repetitive behavior, 
more specifically,  that $f(y)$ defined above  is a regular 
alternating function.

\subsection{Analysis}  

In this subsection,  we attempt to 
establish that $f(y)$ of the last section  is a regular alternating 
function,  
however our results will not be conclusive.
If $f(y)$ is a regular alternating function, then so is $\d_y f(y)$:
\beq
\label{df}
\d_y f(y)   =  \Im \[ \xi (1+iy ) \] .
\eeq
Thus,  one only needs to show that $f'(y)$ is regular alternating.   
Using the summation formula for $g(t)$,  one can show 
\beq
\label{xiN}
\begin{split}
\xi (z ) &= \lim_{N\to \infty}  \xi^{(N)} (z) =  
\lim_{N \to \infty} \sum_{n=1}^N  \xi_n (z) , \\
\xi_n (z) &= n^2 \pi \[  4e^{- \pi n^2 }  - 
z\,  E_{ \tfrac{z-1}{2}} (\pi n^2  )  +  (z-1) 
E_{-\tfrac{z}{2}} (\pi n^2  ) \] .
\end{split}
\eeq
where $E_\nu (r)$  is an incomplete $\Gamma$ function
\beq
\label{Er}
E_\nu (r)  = \int_1^\infty dt \,  e^{-rt} t^{-\nu}  =  
r^{\nu-1} \, \Gamma(1-\nu, r) .
\eeq
It is sometimes referred to as the generalized exponential-integral 
function.    
In obtaining the above equation we have used the identity
\beq
r E_\nu (r)  =  e^{-r}  - \nu E_{\nu +1} (r) .
\eeq

The nature of this approximation is that  the roots $\rho$   
of $\xi^{(N)} (\rho )=0$ provide a very good approximation to the 
smaller Riemann zeros for large enough $N$. However 
small values of $N$ are actually sufficient to a good degree of 
accuracy for small $y$. For instance,   
the first root for $\xi^{(3)}$  coincides with the first Riemann 
zero to 15 digits,  and it's sixth root is correct to 8 digits. 
Furthermore $\xi_{n+1}$ is  smaller than $\xi_n$ because of 
the $e^{-n^2 \pi t}$ suppression in the
integrand for  $E_\nu (\pi n^2)$.
   
For $\nu$ large,   one has the series
\beq
\label{Eseries} 
\begin{split}
E_\nu (r) &=  \sqrt{\frac{\pi}{2}}\,   r^{\nu -1}  \, \csc [(1-\nu)\pi] \, 
e^{-(\nu -1/2) \log \nu  + \nu  - 1/12 \nu + O(1/\nu^3)}  \\
& +  \dfrac{e^{-r}}{\nu}  \left\{  1 - \dfrac{(r-1)}{\nu}  + 
\dfrac{(r^2 - 3r +1)}{\nu^2} + O\((r/\nu)^3\)  \right\}  .
\end{split}
\eeq
Using this,  the leading term for large $y$ is
\beq
\label{Imxi}
\Im \[\xi_n (1+iy)\] \approx - \dfrac{y^2\, e^{-\pi y/4} }{\sqrt{2} n }   
\sin \[  \dfrac{y}{2}  \log\(  \dfrac{y} {2 \pi n^2 e} \) \] .
\eeq

To a reasonably good approximation,   
for large $y$,  $\Im \[ \xi (1+iy) \] \approx \Im \[ \xi_1 (1+iy) \]$,  
and \eqref{Imxi}  indeed is a regularly alternating function because 
the argument of the $\sin$ function is monotonic.     
However we cannot completely rule out that including 
the other terms in $\xi^{(N)}$  for $N>1$ could spoil this behavior.

\section{Transcendental equations for zeros of the 
\texorpdfstring{$\zeta$}{Zeta}-function}
\label{sec:zeta_function}

The main new result presented in the next few sections are transcendental 
equations
satisfied by individual zeros of some $L$-functions.   
For simplicity we first consider the Riemann $\zeta$-function,  
which is the simplest
Dirichlet $L$-function. 
Moreover, we first consider the asymptotic equation 
\eqref{FinalTranscendence}, first proposed in \cite{RHLeclair}, since it 
involves more familiar functions.  
This  asymptotic equation  follows  trivially from 
the exact equation \eqref{exact_eq2}, presented later.

\subsection{Asymptotic equation satisfied by the 
\texorpdfstring{$n$}{n}-th zero on 
the critical line}
\label{sec:zeta_asymptotic}

As above,  let us define the function 
\beq
\label{chidef}
\chi(z) \equiv \pi^{-z/2} \, \Gamma\(z/2\) \zeta(z).
\eeq
which satisfies the functional equation 
\beq
\label{chisym}
\chi\(z\) = \chi\(1-z\).
\eeq
Now consider Stirling's approximation
\beq
\Gamma (z)  = \sqrt{2 \pi}  z^{z - 1/2}e^{-z}\(1+O\(z^{-1}\)\)
\eeq
where $z=x+iy$,  which is valid for large $y$. Under this
condition we also have
\beq
z^z = \exp\( i\(y\log y + \dfrac{\pi x}{2}\) + x \log y - \dfrac{\pi y}{2} 
+ x + O\(y^{-1}\) \).
\eeq
Therefore, using the polar representation
\beq
\zeta = |\zeta| e^{i\arg\zeta}
\eeq
and
the above expansions,  we can write
$$\chi(z) = A \, e^{i\theta}$$
 where
\begin{align}
A(x,y) &= \sqrt{2\pi } \, \pi^{-x/2} \(  \dfrac{y}{2} \)^{(x-1)/2} 
e^{- \pi y /4} |\zeta (x + i y)|\(1+O\(z^{-1}\)\) ,  
\label{A_assymp} \\
\theta(x, y) &= \dfrac{y}{2} \log \( \dfrac{y}{2 \pi e} \)  + 
\dfrac{\pi}{4}(x-1) + \arg \zeta(x + i y) + O\(y^{-1}\) .
\label{theta_assymp}
\end{align}
The above approximation is very accurate. For $y$ as low as $100$, 
it evaluates $\chi\(\tfrac{1}{2}  + i y\)$ correctly to one part in $10^6$.
Above we are assuming $y > 0$. The results for $y<0$ follows trivially
from the relation 
$\(\chi(z)\)^* = \chi(z^*)$.

We will need  the result that the argument on the principle branch, $\Arg f(z)$, of an analytic 
function $f(z)$ has a well 
defined limit at a zero $\rho$ where $f(\rho)=0$.    
Let $\curve$  be a curve in the
$z$-plane such that $z\(\curve\)$ approaches the zero $\rho$ in a 
smooth manner, namely,   $z\(\curve\)$ has a well-defined 
tangent at $\rho$.  Without loss of generality,
let $\rho =0$.   If the zero is of order $k$, then near zero
\beq
f(z) =  a_k z^k +   a_{k+1} z^{k+1} +  \dotsm.
\eeq
Then $\Arg ( f(z)/z^k) $ converges to $\Arg a_k$ along the curve $\curve$.   
Since $z(\curve)$ has a tangent at $0$,  
$\Arg z(\curve)$  converges to a limit $t$  
as $\curve $ approaches $\rho$, 
so that $\Arg f(z) \to \Arg (a_k) + k t$  as $\curve\to \rho$.

Now let $\rho = x+iy$ be a Riemann zero. Then $\Arg \zeta(\rho)$ can be 
well-defined by the limit
\beq \label{deltadef} 
\Arg\zeta\(\rho\) \equiv  \lim_{\delta \to 0^+} \Arg \zeta\(x+\delta+iy\).
\eeq
For reasons that are explained below,  
it is important  that $0 < \delta \ll 1$. This limit in general is not zero. 
For instance, for the first Riemann zero at $\rho = \tfrac{1}{2} +iy_1$,  
where $y_1 = 14.1347\dotsc$, 
\beq
\Arg \zeta\(\tfrac{1}{2} + i y_1\) \approx 0.157873919880941213041945.
\eeq

On the critical line $z=\tfrac{1}{2}+iy$, if $y$ does not
correspond to the imaginary part of a zero, the well-known function 
\beq
\label{SofyDef}
S(y) = \frac{1}{\pi}\arg\zeta\(\tfrac{1}{2}+iy\)
\eeq
is defined by continuous variation along the straight 
lines starting from $2$, then up to $2+iy$ and finally to
$\tfrac{1}{2}+iy$, where $\arg\zeta(2)=0$ (see \eqref{S_variation}).
The function $S(y)$ is discussed in greater detail below in section \ref{Soft}. 
On a zero, the standard way to define this term is through the limit
$S(\rho) = \tfrac{1}{2} 
\lim_{\epsilon\to0} \( S\(\rho+i\epsilon\)+S\(\rho-i\epsilon\) \)$.
We have checked numerically that for several zeros on the line, our definition
\eqref{deltadef} gives the same result  as the standard  piecewise integration definition,  \emph{ so  long as $S(y)$ is on the principal branch,  which is almost always true
for low enough,  but relatively large $y$.}   (See section~\ref{Soft} for more discussion on this important point.)

From \eqref{chidef} it follows that $\zeta(z)$ and $\chi(z)$ have
the same zeros on the critical strip, so it is enough
to consider the zeros of $\chi(z)$. 
Let us now consider  approaching  a zero
$\rho = x + iy$ through the $\delta \to 0^+$ limit in $\arg \zeta$.  
Consider first  the simple zeros along the critical line.    Later we will argue that
all such zeros are in fact simple.    
As we now show,  these zeros are in one-to-one correspondence with the
zeros of the cosine, 
\beq
\label{coszero}
\lim_{\delta \to 0^+}   \cos  \theta = 0.
\eeq
The argument goes as follows.\footnote{A more streamlined argument was presented in \cite{LecMussScattering}.}   
On the critical
line $z = \tfrac{1}{2} + i y$, the functional equation 
\eqref{chisym} implies $\chi (z) =  A(\cos\theta + i \sin \theta)$ is real, 
thus for $y$ \emph{not}  the ordinate
of a zero, $\sin\theta = 0$ and $\cos\theta = \pm 1$.
Thus $\cos \theta $ is a discontinuous function. 
Now let $\yzero$ be the ordinate of a \emph{simple zero}. Then close
to such a zero we define 
\beq
c(y) \equiv  \frac{ \chi \big(\tfrac{1}{2} + i y \big) }{\big|\chi \big(\tfrac{1}{2} + i y\big)\big|}  
=  \frac{ y-\yzero}{|y - \yzero|}.
\eeq
For $y > \yzero$ then $c(y)=1$, 
and for $y<\yzero$ then  $c(y) = -1$.  
Thus $c(y)$ is discontinuous precisely at a zero.
In the above polar representation, formally 
$c(y) =  \cos \theta (\tfrac{1}{2}, y)$. 
Therefore, by identifying zeros as the solutions to $\cos \theta =0$,
we are simply defining  the value of the function  $c(y)$ at 
the discontinuity as  $c(\yzero)=0$.
As explained above,  the argument $\theta$ of $\chi(z)$ 
is well defined
on a zero so this leads to equations satisfied by the zeros. 

\begin{figure}
\begin{center}
\begin{minipage}{.49\textwidth}
\includegraphics[width=1\linewidth]{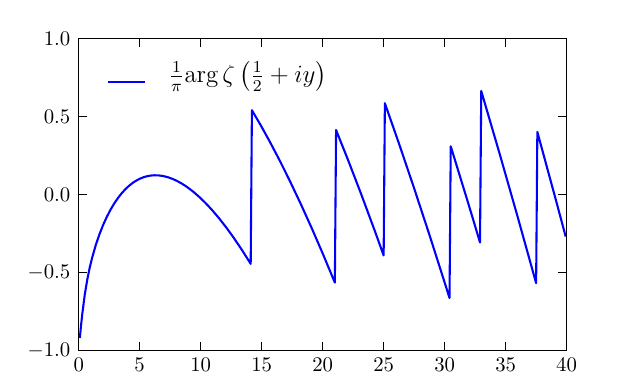}
\end{minipage}%
\begin{minipage}{.49\textwidth}
\includegraphics[width=1\linewidth]{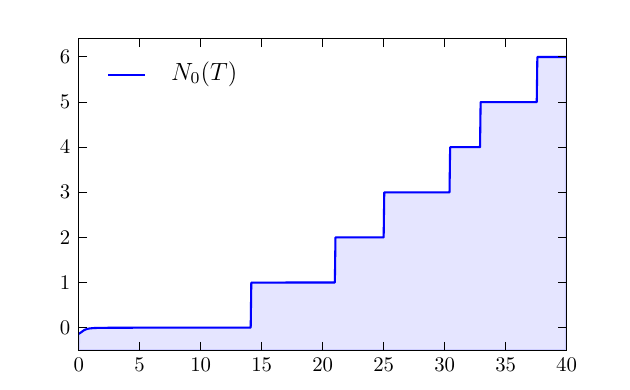}
\end{minipage}
\caption{\emph{ Left:} $\tfrac{1}{\pi}\arg \zeta\(\tfrac{1}{2} + i y\)$  
versus $y$, showing its rapid oscillation. 
The jumps occur on a Riemann zero.
\emph{Right:} $N_0(T)$ versus $T$ in \eqref{counting2}, which is 
indistinguishable from a manual counting of zeros.}
\label{fig:arg_counting}
\end{center}
\end{figure}

The small shift by $\delta$ in \eqref{FinalTranscendence} is essential
since it smooths out $S(y)$,
which is known to jump discontinuously at each zero.
As is  well known, $S(y)$ is a piecewise continuous function, but rapidly 
oscillates around zero with discontinuous jumps, 
as shown in Figure~\ref{fig:arg_counting} (left).
However, when this term  is added to the smooth part 
of $N_0(T)$ (see equations \eqref{counting2} and \eqref{riemann_counting}), 
one obtains an accurate 
staircase function, which 
jumps by one at each zero on the line; see 
Figure~\ref{fig:arg_counting} (right).  
The function $S(y)$ is further discussed in section \ref{Soft}.  
Note that $N_0(T)$ and $N(T)$ are necessarily monotonically increasing 
functions.

The reason $\delta$ needs to be positive
in \eqref{exact_eq2} is the following.     Near a zero $\rho_n$, 
\beq
\zeta (z)  \approx  \(z-\rho_n\) \zeta' \(\rho_n\)
= \(\delta +i \(y-y_n\)\) \zeta'\(\rho_n\).
\eeq
This gives 
\beq
\arg \zeta (z) \approx \arctan\((y-y_n)/\delta\)+ \arg \zeta'(\rho_n).
\eeq
Thus,  with $\delta > 0$,  as one passes through a zero from below,  
$S(y)$ \emph{increases by one}, as it should based on its role in the counting
formula $N(T)$. On the other hand, if $\delta < 0$ then $S(y)$ would decrease
by one instead.

We can now obtain  a precise equation for the location of the zeros on
the critical line.  
The equation \eqref{coszero}, implies 
$\lim_{\delta\to0^+}\theta\(\tfrac{1}{2}+\delta,y\) = \(n+\tfrac{1}{2}\)\pi$, 
for $n=0,\pm 1, \pm 2, \dotsc$, hence
\beq
\label{almost_final_zeta}
n = \dfrac{y}{2\pi}\log\(\dfrac{y}{2\pi e}\) -\dfrac{5}{8}
+ \lim_{\delta\to0^{+}}\dfrac{1}{\pi}\arg \zeta\(\tfrac{1}{2}+\delta+i y\).
\eeq
A closer inspection shows that the RHS of
the above equation  has a minimum in the interval $(-2, -1)$, thus $n$ 
is bounded from below, i.e. $n \ge -1$. 
Establishing the \emph{convention} that zeros 
are labeled by positive integers, $\rho_n = \tfrac{1}{2}+i y_n$ where
$n=1,2,\dotsc$, we must replace $n \to n - 2$ in \eqref{almost_final_zeta}. 
Therefore, the imaginary
parts of these zeros satisfy the transcendental equation
\beq
\label{FinalTranscendence} 
\dfrac{y_n}{2 \pi}  \log \( \dfrac{y_n }{2 \pi e} \)   
+ \lim_{\delta \to 0^{+}}  \dfrac{1}{\pi} 
\arg  \zeta \( \tfrac{1}{2}+ \delta + i y_n  \) = n - \dfrac{11}{8}
\qquad (n=1,2,\dotsc).
\eeq
In summary, we have shown that, asymptotically for now, there are an  infinite 
number of  zeros on the critical line whose ordinates can be determined by 
solving \eqref{FinalTranscendence}. This equation was first proposed in \cite{RHLeclair}. This equation determines
the zeros on the upper half of the critical line. The zeros on the lower
half are symmetrically distributed; if 
$\rho_n = \tfrac{1}{2}+iy_n$ is a zero, so is 
$\rho_n^* = \tfrac{1}{2}-iy_n$.

The LHS of \eqref{FinalTranscendence} is a monotonically  
increasing function of $y$, and the leading term is a smooth function. 
This is clear since the same terms appear in the staircase function 
$N(T)$ described below.   
Possible discontinuities can only come from 
$\tfrac{1}{\pi}\arg\zeta\(\tfrac{1}{2} +iy\)$, and in fact, it has a jump 
discontinuity by one whenever $y$ corresponds to a zero on the line.   
However, if
$\lim_{\delta \to 0^+}\arg\zeta\(\tfrac{1}{2}+\delta+iy\)$ is well
defined for every $y$, then the left hand side of equation 
\eqref{FinalTranscendence} is well defined for any $y$,  and due to 
its monotonicity, there  must be a 
unique solution for every $n$.
Under this assumption, the number of solutions of equation 
\eqref{FinalTranscendence}, up to height $T$, is given by
\beq \label{counting2}
N_0(T) = \dfrac{T}{2 \pi} \log \( \frac{T}{2 \pi e} \) + \frac{7}{8}  + 
\inv{\pi} \arg \zeta \( \tfrac{1}{2} + i T \) + O\(T^{-1}\).
\eeq
This is so because the zeros are already numbered in 
\eqref{FinalTranscendence}, but the left hand side jumps by one at each 
zero, with values $-\tfrac{1}{2}$ to the left and $+\tfrac{1}{2}$ to the 
right of the zero. 
Thus we can replace $n \to N_0 + \tfrac{1}{2}$ and $y_n \to T$, such that the 
jumps correspond to integer values. In this way $T$ will not correspond 
to the ordinate of a zero and $\delta$ can be eliminated.

Using Cauchy's argument principle (see Appendix~\ref{sec:number_zeros}) 
one can derive the Riemann-von Mangoldt 
formula, which gives the number of zeros in the region
$\{ 0 < x < 1, \, 0 < y < T\}$ inside the \emph{critical strip}.
This formula is standard \cite{Edwards, Titchmarsh}: 
\beq \label{riemann_counting}
N(T) = \dfrac{T}{2 \pi} \log \( \frac{T}{2 \pi e} \) + \frac{7}{8}  + 
S(T) + O\(T^{-1}\).
\eeq
The leading $T\log T $ term was already in Riemann's original paper.  
Note that it has the same form as the counting formula  on 
the \emph{critical line} that we have just found \eqref{counting2}.   
Thus, under the assumptions we have described,  we conclude that 
$N_0(T) = N(T)$ asymptotically.    
This means that our particular solution \eqref{particular_sol}, leading
to equation \eqref{FinalTranscendence}, already saturates the counting 
formula on the whole strip and there are no additional zeros from $A=0$ 
in \eqref{sumchi} nor from the  more general equation 
$\theta + \theta' = (2n+1)\pi$ described below. 
This strongly suggests that
\eqref{FinalTranscendence} describes all non-trivial zeros of $\zeta(z)$, 
which must then lie on the critical line.

\subsection{Exact equation satisfied by the 
\texorpdfstring{$n$}{n}-th zero on the critical line}
\label{sec:zeta_exact}

Let us  now repeat the previous analysis  but without
considering an asymptotic expansion. The exact versions of
\eqref{A_assymp} and \eqref{theta_assymp} are
\begin{align}
\label{A_exact}
A(x,y) &= \pi^{-x/2} |\Gamma\(\tfrac{1}{2}(x+iy)\)| |\zeta(x+iy)|, \\
\label{theta_exact}
\theta(x,y) &= 
\arg \Gamma\(\tfrac{1}{2}(x+iy)\) -\dfrac{y}{2}\log\pi + \arg \zeta(x+iy),
\end{align}
Then, as before, zeros are described by 
$\lim_{\delta\to0^+}\cos\theta = 0$, equivalent to
$\lim_{\delta\to0^+}\theta\(\tfrac{1}{2}+\delta,y\) = \(n+\tfrac{1}{2}\)\pi$,
and upon replacing $n \to n-2$ the imaginary parts of these zeros must 
satisfy the exact equation
\beq\label{exact_eq}
\arg\Gamma\(\tfrac{1}{4}+\tfrac{i}{2}y_n\) - y_n \log\sqrt{\pi} 
+ \lim_{\delta \to 0^+} \arg\zeta\(\tfrac{1}{2}+\delta + iy_n\) = 
\(n-\tfrac{3}{2}\)\pi.
\eeq
The Riemann-Siegel $\vartheta$ function is defined by
\beq
\label{riemann_siegel}
\vartheta(y)\equiv \arg\Gamma\(\tfrac{1}{4}+\tfrac{i}{2}y\) - y\log \sqrt{\pi},
\eeq
where the argument is defined such that this function is continuous and
$\vartheta(0)=0$. This can be done through the relation
$\arg \Gamma = \Im \log \Gamma$, and numerically one can use the 
implementation of the ``logGamma'' function. This is equivalent to the analytic
multivalued $\log\(\Gamma\)$ function, but it simplifies its complicated
branch cut structure.
Therefore, there are an infinite  number of zeros in 
the form
$\rho_n=\tfrac{1}{2}+iy_n$, where $n=1,2,\dotsc$, whose imaginary
parts \emph{exactly} satisfy the following equation:
\beq\label{exact_eq2}
\vartheta(y_n) + 
\lim_{\delta\to 0^{+}}\arg\zeta\(\tfrac{1}{2}+\delta +iy_n\) = 
\(n-\tfrac{3}{2}\)\pi \qquad (n=1,2,\dotsc).
\eeq
Expanding the $\Gamma$-function in \eqref{riemann_siegel} through 
Stirling's formula 
\beq
\label{AsymRS}
\vartheta (y) =  \frac{y}{2}  \log \( \frac{y}{2 \pi e} \)  - 
\frac{\pi}{8}  + O(1/y) 
\eeq
one recovers the asymptotic equation 
\eqref{FinalTranscendence}. 

Again, as discussed after \eqref{FinalTranscendence}, the first
term in \eqref{exact_eq2} is smooth and the whole left hand side is
a monotonically increasing function. If 
$\lim_{\delta\to 0^{+}}\zeta\(\tfrac{1}{2}+\delta+iy\)$ is well defined
for every $y$, then equation \eqref{exact_eq2} must have a unique solution
for every $n$.  Under this condition it is valid to replace 
$y_n \to T$ and $n \to N_0 + \tfrac{1}{2}$,  and then
the number of solutions of \eqref{exact_eq2} is given by
\beq\label{counting2_exact}
N_0(T) = \dfrac{1}{\pi}\vartheta(T) + 1 +
\dfrac{1}{\pi}\arg\zeta\(\tfrac{1}{2}+iT\).
\eeq

The exact Backlund counting formula (see Appendix~\ref{sec:number_zeros}), 
which gives the number of zeros on
the critical strip with $0<\Im(\rho)<T$, is given by \cite{Edwards}
\beq\label{backlund}
N(T) = \dfrac{1}{\pi}\vartheta(T) +  1 + S\(T\).
\eeq
Therefore, comparing \eqref{counting2_exact} with the exact counting 
formula on the 
\emph{entire critical strip} \eqref{backlund}, we have $N_0(T) = N(T)$ exactly. 
This indicates that our particular solution, 
leading to equation \eqref{exact_eq2}, captures all the zeros on the strip \emph{assuming there is a solution for every $n$},
indicating that they should all be on the critical line.

In summary, if \eqref{exact_eq2}  has a unique solution for each $n$,
then this saturates the counting formula for the entire 
critical strip and this would establish the validity of the RH.

\subsection{A more general equation} 

The above equation \eqref{exact_eq2}  was first obtained by us with  
a different, and less rigorous,  argument \cite{RHLeclair, FL1}.
It  is  a particular solution of a more general formula
which we now present.  

We will  need the following.   
From \eqref{chidef} we have $\(\chi(z)\)^{*} = \chi\(z^*\)$, thus 
$A(x,-y)=A(x,y)$ and $\theta(x,-y)=-\theta(x,y)$. 
Denoting 
\beq
\chi\(1-z\) = A' \, e^{-i\theta'}
\eeq
we then have
\beq
A'(x,y)=A(1-x,y), \qquad \theta'(x,y)=\theta(1-x,y).
\eeq
From \eqref{chisym} we also have $|\chi(z)| = |\chi(1-z)|$, 
therefore 
\beq\label{AAp}
A(x,y) = A'(x,y)
\eeq
for \emph{any} $z$ on the critical strip.

From \eqref{chisym} we see 
that if $\rho$ is a zero so is $1-\rho$. Then we clearly have
\beq
\label{sumchi}
\lim_{z \to \rho} \left\{ \chi(z) + \chi(1-z) \right\} = \lim_{z\to \rho} 
A(x, y) \, B(x, y) = 0, 
\eeq
where we have defined 
\beq\label{B_def}
B(x,y) =    e^{i \theta(x,y) }  +  e^{-i \theta'(x,y)}.
\eeq
The second equality in \eqref{sumchi} follows from \eqref{AAp}. 
For now,  we do not specify the precise curve $\curve$ through which we 
approach the zero.

The above equation \eqref{sumchi} is identically satisfied on a 
zero $\rho$ since $\lim_{z\to\rho}A \sim |\zeta (\rho)| =0$,
independently of $B$.
However, this by itself  does not provide any more detailed information 
on the zeros.
There is much more information in the phases $\theta$ and $\theta'$.
Consider instead taking the limits in $A$ and $B$ separately,
\beq
\label{ABlimit}
\lim_{z' \to \rho} \, \lim_{z\to \rho} \, A(x',y') \, B(x, y) = 0, 
\eeq
where $z' = x' + i y'$.  
Taking $z \to \rho $ first,  a potential zero occurs when 
\beq
\label{Blimit}
\lim_{z\to \rho}  B(x  , y)  =  \lim_{z \to \rho}   
\( e^{i \theta} + e^{-i\theta'} \) =0.
\eeq
We propose that Riemann zeros satisfy \eqref{Blimit}.  
The equation $B=0$ provides more information on the location of zeros 
than $A=0$ since the phases $\theta$ and $\theta'$
can be well defined at a zero through an appropriate limit.  
We emphasize that we have not yet assumed the RH,  and the above 
analysis is valid on the entire 
complex plane,  except at $z=1$ due to the simple pole of $\chi$.  
We will provide ample evidence  that the equation 
\eqref{Blimit}  is evidently correct even for the example of 
the Davenport-Heilbronn function,  which has zeros off the critical line,  and 
 the RH fails.
Clearly a more rigorous derivation would
be desirable, the delicacy being the limits involved,  but let us proceed.   

The linear combination in \eqref{sumchi}  was chosen 
to be manifestly symmetric under $z\to 1-z$. Had we taken a different 
linear combination in \eqref{sumchi}, such as $\chi(\rho) + b \, \chi(1-\rho)$, 
then $B = e^{i\theta} + b \, e^{-i\theta'}$ for some constant $b$.
Setting the real and imaginary parts of $B$ to zero gives the two equations
$\cos \theta + b\cos \theta' =0$ and $\sin\theta - b  \sin \theta'=0$.    
Summing the squares of these equations one obtains  
$\cos (\theta + \theta')=-(b+1/b)/2$. However, since $b+1/b >1$,  there 
are no solutions except for $b=1$.

The general solution of \eqref{Blimit} is given by 
\beq
\label{Bzero} 
\theta + \theta' = (2n+1)\pi.
\eeq
Note that $\chi(z) = \chi(1-z)$ implies $\theta + \theta' = 2\pi n$ 
for $z\ne \rho$. This 
together with \eqref{Bzero} is analogous to the previous discussion
where $\cos \theta = 1$ for $z\ne \rho$ and $\cos\theta = 0$ for
$z = \rho$.     All zeros satisfying \eqref{Bzero}   are  simple since there 
are in correspondence with zeros of the cosine or sine function.   

The zeros on the \emph{critical line} correspond to the particular solution
\beq \label{particular_sol}
\theta = \theta' = (n+\smallhalf) \pi ,
\eeq
which is equivalent to \eqref{coszero} and \eqref{exact_eq2}.

\begin{figure}
\begin{center}
\includegraphics[width=.5\linewidth]{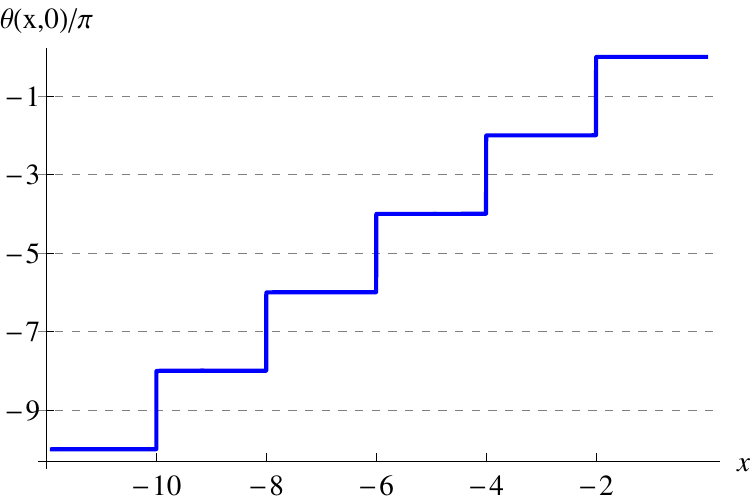}
\caption{Plot of equation \eqref{eq_trivial}.
Note the jumps occuring at $2n+1$ corresponding to the trivial zeros
of the $\zeta$-function.}
\label{trivial_zeros} 
\end{center}
\end{figure}

In fact, the trivial zeros along the negative $x$ axis also 
satisfy \eqref{Bzero},  which 
again strongly supports its validity.   One can show this as follows. 
For $x < 0$, $\theta'(x,0) = 0$. Since we are on the real line,
we approach the zero through the path $z = \rho - i \epsilon$, with
$0<\epsilon \ll 1$. This path smooths out the $\arg\zeta$ term in the
same way as the $\delta\to 0^+$ for zeros on the critical line. Then 
from \eqref{Bzero} we are left with
\beq\label{eq_trivial}
\dfrac{1}{\pi}\theta(x,0) = \dfrac{1}{\pi}\Im\[ \log\Gamma\(x/2\) \] + 
\dfrac{1}{\pi}\lim_{\epsilon\to0^+}{\rm Arg} \, \zeta(x-i\epsilon) + 2 k 
= 2n+1.
\eeq
Note that we write $\arg\zeta(x-i\epsilon) = {\rm Arg} \, \zeta(x-i\epsilon)
+ 2\pi k $, where we have the principal value 
${\rm Arg}\, \zeta(x) \in \{0,\pm \pi\}$ for $x < 0$. The changes in branch
are accounted for by $k$, and  depends  on $x$. If we take these  changes correctly
into account we obtain Figure~\ref{trivial_zeros}, showing
the trivial zeros as solutions to \eqref{eq_trivial}. The first
term in \eqref{eq_trivial}, i.e. $\tfrac{1}{\pi}\Im\[\log\Gamma(x/2)\]$,
is already a staircase function with jumps by $1$ at every negative even $x$.
The other two terms, $\tfrac{1}{\pi}{\rm Arg} \, 
\zeta(x-i\epsilon) + 2 k$, just shift the function by a constant, 
such that these  jumps
coincide with odd integers $(2n+1)$ in such a way that \eqref{eq_trivial}
is satisfied exactly at a zero.
Thus, remarkably, the equation \eqref{Bzero}  characterizes all known zeros 
of $\zeta$. 

\subsection{On possible zeros off of the line} 

\begin{figure}
\centering
\begin{minipage}{0.49\textwidth}
\includegraphics[width=1\linewidth]{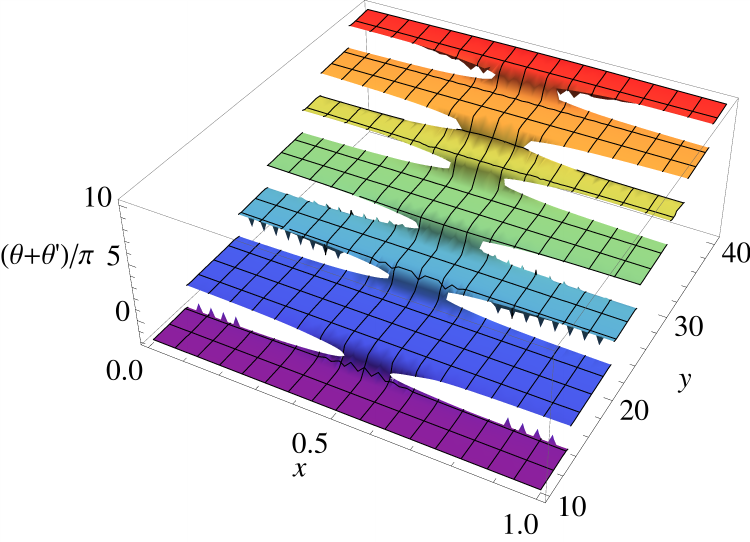}
\end{minipage}
\begin{minipage}{0.49\textwidth}
\includegraphics[width=1\linewidth]{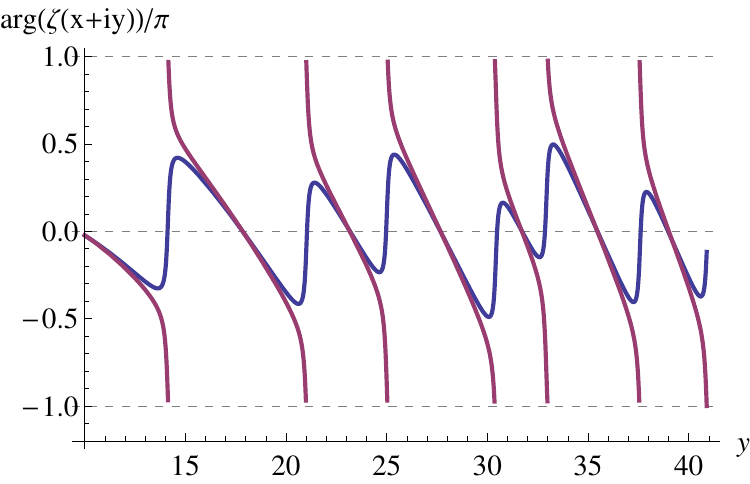}
\end{minipage}
\caption{
\emph{ Left:} $3D$ plot of $\tfrac{1}{\pi}\[\theta(x,y) + \theta' (x,y)\]$.
The function does not vary with $x$ on each plateau, but
it is only well defined between jumps very close to the critical line.
\emph{ Right:}
$\tfrac{1}{\pi}\arg \zeta (x+iy)$  for $x=1/2 + 0.1$ (blue)  
and $x=1/2-0.1$ (purple)  as a function of $y$.}
\label{fig_ArgPlusMinus}
\end{figure}

Suppose one looks for solutions to \eqref{Bzero} off the line.   
In Figure~\ref{fig_ArgPlusMinus} (left)  we plot the RHS
of \eqref{Bzero}  divided by $\pi$ for a region on the critical strip.
One clearly sees 
that precisely where a solution requires
that this equals an odd integer, the function is not well-defined.   
On the other hand, for $x=1/2$ with the $\delta$-prescription  it is 
well-defined and has a unique solution.

The fact that the RHS of \eqref{Bzero} is not well defined for $x>1/2$  
is due to the very different properties of $\arg \zeta (x+iy)$  
for $x>1/2$ verses $x<1/2$.      
In Figure~\ref{fig_ArgPlusMinus} (right)  we plot  $\arg \zeta (x+ iy)$ for 
$x=1/2 \pm 0.1$ as a function of $y$.    
One sees that for $x<1/2$ there are
severe changes of branch where the function is not defined,  
whereas for $x>1/2$ 
it is smooth.    
Since the RHS of \eqref{Bzero}  involves both $\theta$ and $\theta'$, 
the $\theta'$ term is ill-defined for $x>1/2$ and thus neither is  the RHS.
Only on the critical line where $\theta = \theta'$ and $x= 1/2 + \delta$ with 
$\delta \to 0^+$ is the RHS well-defined.     

There  is another very interesting 
aspect of Figure \ref{fig_ArgPlusMinus} (left).
On each plateau the function is mainly constant
with respect to $x$.
This formally follows from
\beq
\d_x \[  \theta (x,y) + \theta' (x,y)\] = 0
\eeq
if one assumes $\theta$ is differentiable at $x$.   
This  leads to the following suggestion.    Suppose that the dependence on
$x$ for $1/2<x\leq 1$  is weak enough that $\theta + \theta'$ is very well 
approximated
by the curve at $x=1$.     
Recall that it is known that there are no zeros along the line
$x=1$.     If the curve for $1/2 < x < 1$ is a smooth and very small 
deformation of the
one at $x=1$,  then there are no solutions to \eqref{Bzero} off of the line,
and if the latter captures all zeros,  then there are no zeros off the line.  
As in section \ref{Electric},   the RH would then be  related to the 
non-existence of zeros
at $x=1$,  which is equivalent to the prime number theorem.
The main problem with this argument is that at a zero off the line,  
probably the above derivative is not well defined.

One sees that the particular solution \eqref{particular_sol}  of the more 
general
$B=0$ is a consequence of the direction in which the zero on the line 
is approached. 
Let $\zeta (z) = u(x,y) + i v(x,y)$ as in section~\ref{Electric}.   
The $u,v=\mbox{const.}$ contours
are sketched in Figure~\ref{uvcontours}.    In the transcendental equation 
\eqref{exact_eq}   the $\delta\to 0^+$ limit approaches the zero on 
the critical line along $u=0$ contours that are  nearly in the 
$x$ direction.   
For potential zeros off of the line where $u,v=0$ contours intersect 
perpendicularly, 
one does not expect that the directions of these contours at the zeros is 
always the same,  in contrast to those of the zeros on the critical line. 
Thus,  for zeros off of the line,  we expect $B=0$ will be satisfied, i.e. 
\eqref{Bzero},  but not  the particular solution $\cos \theta =0$.
 
In section~\ref{sec:Davenport}  we will study an example of an $L$-series that
does not satisfy the RH,  the Davenport-Heilbronn function.   
We will show that the  zeros off the line indeed satisfy \eqref{Bzero}.

\subsection{Further remarks}
\label{sec:zeta_remarks}

It is possible to introduce a new function
$\zeta(z) \to \widetilde{\zeta}(z) = f(z) \zeta(z)$ that also 
satisfies the functional equation \eqref{chisym}, 
i.e. $\tilde{\chi}(z)=\tilde{\chi}(1-z)$, but has zeros 
off of the critical line due to the zeros of $f(z)$. In such a case 
the corresponding functional equation will hold if and only 
if $f(z) = f(1-z)$ for any $z$, and this is a trivial condition 
on $f(z)$, which could have been canceled in the first place. 
Moreover, if $f(z)$ and $\zeta(z)$ have different zeros, the analog of 
equation \eqref{sumchi} has a  factor
$f(z)$, i.e. $\tilde{\chi}(\rho+\delta) + \tilde{\chi}(1-\rho-\delta) = 
f(\rho+\delta)\[\chi(\rho+\delta)+\chi(1-\rho-\delta)\]=0$,  
implying \eqref{sumchi} again  where $\chi(z)$ is the original \eqref{chidef}. 
Therefore, the previous analysis eliminates $f(z)$
automatically and only finds the zeros of $\chi(z)$.
The analysis is  non-trivial precisely  because $\zeta(z)$ satisfies 
the functional equation but $\zeta(z) \ne \zeta(1-z)$. Furthermore, 
it is a well known theorem that the only function which satisfies the 
functional equation \eqref{chisym} and has the same characteristics 
of $\zeta(z)$, is $\zeta(z)$ itself. 
In other words, if $\widetilde{\zeta}(z)$ is required to have  the 
same properties of $\zeta(z)$, then $\widetilde{\zeta}(z) = C \, \zeta(z)$, 
where $C$ is a constant \cite[pg. 31]{Titchmarsh}.

Although equations \eqref{exact_eq2} and \eqref{backlund} have an obvious 
resemblance, it is impossible to derive the former from the later, 
since the later is just a counting formula valid on the entire strip, 
and it is assumed that $T$ is \emph{not} the ordinate of a zero.
Moreover, this would require the assumption of the validity of the RH,
contrary to our approach, where we derived equations \eqref{exact_eq2}
and \eqref{FinalTranscendence} on the critical line, without assuming
the RH.
Despite our best efforts, we were not able to find 
equations \eqref{FinalTranscendence} and \eqref{exact_eq2} 
in the literature. 
Furthermore,  the counting formulas \eqref{counting2} and \eqref{backlund}
have  never been proven  to be valid on the critical line \cite{Edwards}.

\section{The argument of  the Riemann \texorpdfstring{$\zeta$}{Zeta}-function}
\label{Soft}

Let us recall the definition used in section \ref{sec:zeta_function} in conjunction with the Fran\c ca-LeClair equation \eqref{exact_eq2},  namely
\beq
\label{Sydef}
S(y)  =   \lim_{\delta \to 0^+}   
\inv \pi  \arg \zeta \(\tfrac{1}{2} +\delta +   iy \)  
= \lim_{\delta \to 0^+}   
\inv \pi  \Im \[ \log \zeta\(\tfrac{1}{2} +\delta +   iy\)\].
\eeq
Previously, we argued that $S(y_n)$ is well defined at a 
zero $\rho = \tfrac{1}{2} + iy_n$ in the non-zero $\delta \to 0^+$ limit.
A proper understanding of this function is essential in any theory of 
the Riemann zeros because of its role in the counting function $N(T)$,  
and in our equation \eqref{exact_eq2} satisfied by individual zeros.   
It is the fluctuations in $S(y)$ that ``knows'' about the actual zeros.  
As stated above,  if the equation  \eqref{exact_eq2}  has a unique 
solution for every $n$,  then the RH would follow since 
then $N_0 (T) = N(T)$.
In this section we describe some important  properties of $S(y)$  in connection with the equation \eqref{exact_eq2}.    
However this will still not be sufficient to prove  that there is indeed a unique solution to this equation.     
All of these properties are essentially already  known (see for instance \cite{Karatsuba}).
 
The conventional way to define $S(y)$ is by piecewise 
integration of $\zeta' /\zeta$
from $z=2$ to $2+iy$,  then to $1/2 + iy$.
Namely
$\arg \zeta ( \smallhalf + i y) =  \arg \zeta(2+iy)  +  \Delta$,  
where $\Delta = \arg \zeta (\smallhalf + i y) - \arg \zeta ( 2+ i y)$.
The integration to arbitrarily high $y$  along $\Re (z) =2$ is bounded 
and gives something
relatively small on the principal branch.
This can be seen from the Euler product.
For $x>1$, 
\beq
\label{arg2}
\arg \zeta (x+ iy) =  \Im \log \zeta (x+ iy)  =  \inv{2i}  \sum_p \log\(
\frac{  1- p^{-x+iy} }{1-p^{-x-iy} }\) \approx -  
\sum_p \inv{p^x} \sin ( y \log p ) .
\eeq
For $x=2$,
\beq
\label{arg2b}
| \arg \zeta (2+ iy)|  <  \sum_{p}  \inv{p^2}  =  0.452235\dotsc.
\eeq
The above sum obviously converges since it is less than 
$\zeta (2) = \pi^2 /6 = 1.645$.
 
It is known that $S(y)$ is unbounded,   which implies that the logs in \eqref{arg2},  which are Arg's,  can accumulate.  
The roughest estimate is the known fact that $S(y) = O(\log y)$,  which actually was not written explicitly  in Riemann's paper,   but only later by
von Mangoldt. 
Such a logarithmic growth  apparently should come from the short integration 
that gives $\Delta$ of the last paragraph.  
Assuming the RH, the current 
best bound is given by \cite{Milinovich} (see also \cite{Goldston})
$$|S(y)|\le\(\dfrac{1}{4}+o(1)\)\dfrac{\log y}{\log\log y}$$ for 
$y\to \infty$. 
It is important to bear in mind that these are upper bounds, 
and that $S(y)$ may actually be much smaller.

\medskip

The first three properties of $S(y)$ listed below  
are well-known \cite{Edwards,Titchmarsh,Karatsuba}:
 
\begin{enumerate} 
 
\item At each zero $\rho = x +i y$ in the critical strip, 
$S(y)$  jumps by  the multiplicity $m$ of the zero.   
This simply follows from the role of $S(T)$ in the counting
formula $N(T)$  in \eqref{backlund}.
For instance,  simple zeros on the critical line have $m=1$, 
whereas double zeros on the line have $m=2$.
Since zeros off the line at a given height $y$ always occurs in
pairs, i.e. $\rho$ and $1-\rho^*$, if one of such zeros has
multiplicity $m$, then $S(y)$ has to jump by $2m$ at this height $y$.
It is believed that all the Riemann zeros are simple, although this is
largely a completely open problem. 
 
\item Between zeros,  since $N(y)$ is constant, 
\beq
\label{Sprime}
S' (y) = \d_y S(y) = -\inv{\pi}  \vartheta' (y) < 0,
\eeq
where the last inequality follows because
$\vartheta (y)$ is a monotonically increasing function.   

\item  
\label{Sav}
The average $\langle S \rangle$  of $S(y)$ is zero \cite{Edwards},
\beq
\label{Saveeq}
\langle S \rangle =  
\lim_{Y\to \infty}   \inv{Y}   \int_0^Y   dy \,  S(y)  = 0.
\eeq

\item
\label{loglog}
 A celebrated theorem of Selberg \cite{Selberg} states that $S(y)$ over a large interval $0 < y < T$ satisfies a normal distribution with zero mean and variance
\beq
\label{Selberg} 
\overline{S(y)^2}   \equiv \inv{T} \int_0^T  S (y)^2 \, dy = \inv{2 \pi^2}  \log \log T + O\( \sqrt{\log \log T} \).
\eeq
The numerical work  performed for these lectures did not go higher than $y=10^{10}$.   Over this range the variance of $S(y)$ is approximately $0.16$,
which implies it is nearly always on the principal branch.   For this reason it is typically OK  to compute ${\rm Arg} \zeta (\half + i y)/\pi$ 
in solving \eqref{exact_eq2}
since it almost always equals $S(y)$,   especially near the zeros.    One must bear this in mind when trying to solve \eqref{exact_eq2} and treat this issue with some care,  however for the ranges of $t$ numerically considered here this did not present any difficulty.    
For much much higher $y$,  one needs to keep track of branches to determine the actual $\arg$.

\item  
\label{DeltaS}
Let 
\beq
\label{dS}
\Delta S_n \equiv    S(y_n) - S(y_{n+1})  =   
\inv{\pi} \(  \vartheta (y_{n+1}) -  \vartheta(y_n) \)
\eeq
where the equality  follows from \eqref{Sprime}.  
Then, if the RH is true,  $\Delta S_n$ has to compensate the 
jumps by $1$ at each simple zero, and since $\langle S \rangle = 0$,  one expects   
\beq
\label{DSave} 
\langle \Delta S_n  \rangle = 1.
\eeq

There is one more property we will need,  which is a precise statement of
the fact that the real part of $\zeta (\tfrac{1}{2} + iy )$  
is \emph{almost always}  positive.   
Let $\yplus_n$ and $\yminus_n$,  where $n=1,2,\dotsc$,
denote the points where either the imaginary or real  part,  
respectively, of $\zeta (\tfrac12+iy)$ are zero, but not both.    
These points are easy to find since they do not depend on the 
fluctuating $S(y)$.  We have
\beq
\label{ReIm}  
\zeta (\smallhalf - i y)  = \zeta\(\smallhalf+i y\)  \,  G(y), 
\qquad G(y) =  e^{2 i \vartheta (y)}
\eeq
where $\vartheta (y)$ is the smooth Riemann-Siegel function 
\eqref{riemann_siegel}. 
Since the real and imaginary
parts are not both zero,   
at $\yplus_n$ then  $G=1$,  whereas at $\yminus_n$ then $G=-1$.
Thus 
\begin{align}
\Im \[ \zeta \(\smallhalf +i \yplus_n \) \] &=  0
&\mbox{for} \ \vartheta(\yplus_n) &= (n-1)\pi, \label{Yplus}\\
\Re \[ \zeta \(\smallhalf +i \yminus_n \) \] &= 0 
&\mbox{for} \ \vartheta(\yminus_n) &= \(n-\tfrac{1}{2}\)\pi. \label{Yminus}
\end{align}
Our convention is $n=1$ for the first point where this occurs for $y>0$.
Using the approximation \eqref{AsymRS}, equations 
\eqref{Yplus} and \eqref{Yminus} can be written in the form
$\tfrac{y}{2\pi}\log\(\tfrac{y}{2\pi e}\) = A_n$, which through the
transformation $y \to 2\pi A_n x^{-1}$ can be solved in terms of the
Lambert $W$-function (see section~\ref{sec:lambert}). The result is
\beq
y_n^{(+)} \approx  \dfrac{2\pi\(n-7/8\)}{W\[e^{-1}(n-7/8)\]},  \qquad
y_n^{(-)} \approx \dfrac{2\pi\(n-3/8\)}{W\[e^{-1}(n-3/8)\]}, 
\label{Ypm}
\eeq
where above $n=1,2,\dotsc$ and $W$ denotes the principal branch
$W_0$. The $\yplus_n$ are actually the Gram points.
From \eqref{Ypm} we can see that these points are ordered in a regular manner,
\beq 
\yplus_1 < \yminus_1 < \yplus_2 < \yminus_2 < \yplus_3 <  \yminus_3 < \dotsm
\eeq
as illustrated in Figure~\ref{RealoverImag} (left).

\begin{figure}
\centering
\begin{minipage}{0.49\textwidth}
\vspace{1.6em}
\includegraphics[width=1\linewidth]{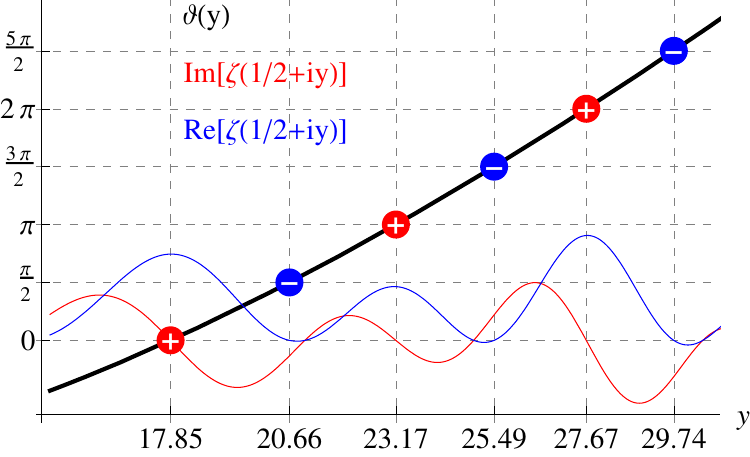}
\end{minipage}
\begin{minipage}{0.49\textwidth}
\includegraphics[width=1\linewidth]{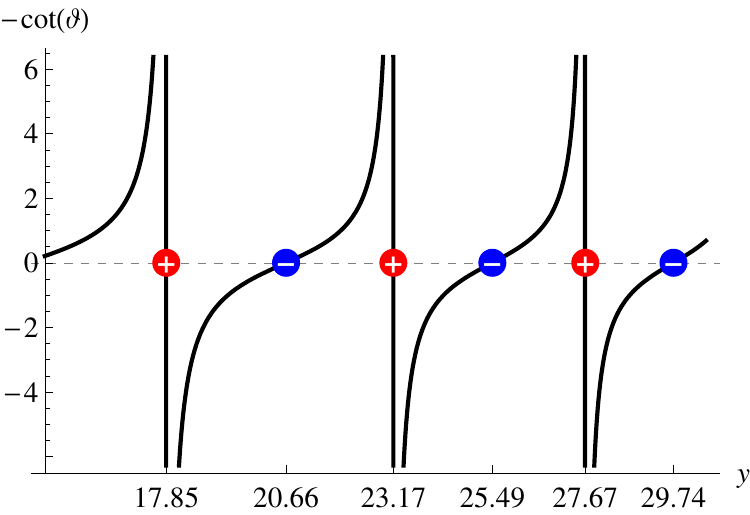}
\end{minipage}
\caption{\emph{ Left:} the $y_n^{(+)}$ (red balls) are the Gram points, where
$\Im \[ \zeta\(\tfrac{1}{2}+iy\) \] = 0$, 
and $y_n^{(-)}$ (blue balls) are the points where 
$\Re \[ \zeta\(\tfrac{1}{2}+iy\) \] = 0$. These points are determined approximately from
formulas \eqref{Ypm}. \emph{ Right: } A plot of \eqref{ratio} indicating
the points \eqref{Ypm}.
}
\label{RealoverImag} 
\end{figure}

The ratio 
\beq
\label{ratio}
\dfrac{\Re\[ \zeta \(\tfrac{1}{2}+iy\) \]}{\Im\[ \zeta \(\tfrac{1}{2}+iy\) \]} 
=  - \cot \vartheta(y) 
\eeq
has a regular repeating pattern,  as can be seen in 
Figure~\ref{RealoverImag} (right),  thus the signs
of the real and imaginary parts are related in a specific manner.    
From this figure one 
sees that when the imaginary part 
is negative the real part is positive.

\item
\label{titchave}
The statement  we need about the real part being mainly
positive concerns the average value of the real part at the 
points $\yplus_n$.
The average of the real part of $\zeta \(\tfrac{1}{2} + i y\)$ at the 
Gram points $\yplus_n$ is \cite{Titchmarsh}
\beq \label{TitchAverage}
\lim_{N\to \infty}   \inv{N}   \sum_{n=1}^N  
\Re \[ \zeta \(\tfrac{1}{2} + i \yplus_n \) \]  = 2.
\eeq
\end{enumerate}

Now, let us consider the behavior  of $S(y)$ starting from the first zero.  
At the first zero,  in the jump by $1$,  $S(y)$ passes through zero and 
remains on the principle branch (see Figure~\ref{fig:arg_counting}).
The branch cut in the
$z$-plane is along the negative $x$-axis, thus on the principle branch 
$-1 <  S(y) \le 1$.
At the points $\yplus_n$, where the imaginary part is zero, 
the vast majority of them have 
$\Re\big[\zeta\big(\tfrac{1}{2}+i\yplus_n\big)\big] > 0$ 
according to item \ref{titchave} above, and thus 
for the most part $S(\yplus_n) =0$. 
Thus $S(y) =0$ at infinitely  many points $\yplus_n$  between zeros,
consistent with $\langle S \rangle = 0$.  
    
\begin{figure}
\centering
\begin{minipage}{0.49\textwidth}
\includegraphics[width=0.9\linewidth]{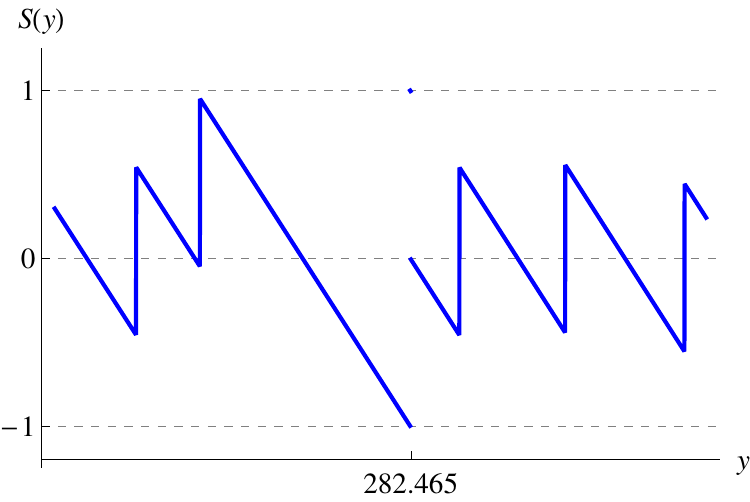}
\end{minipage}
\begin{minipage}{0.49\textwidth}
\includegraphics[width=0.9\linewidth]{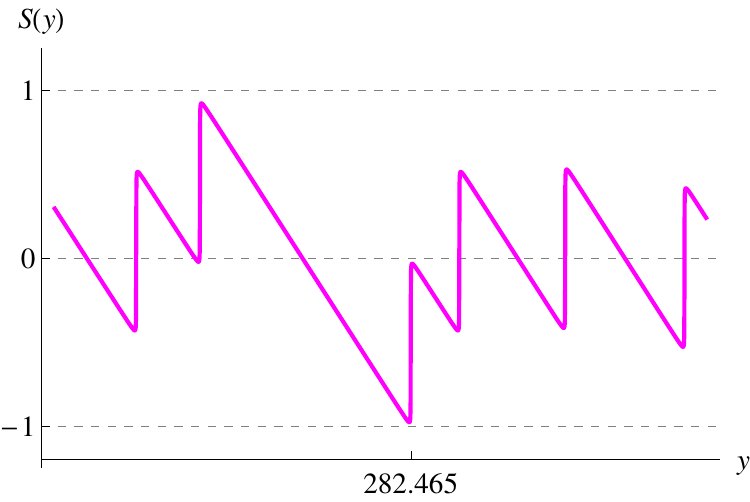}
\end{minipage}
\caption{\emph{ Left:} $S(y)$ in the vicinity of the first point 
where $|S(y)|>1$, at $\yplus_{127} = 282.455$. There is a zero
at $y_{127} = 282.465$. \emph{ Right:} by adding a $\delta$ in 
$S(y) = \tfrac{1}{\pi}\arg\zeta\(\tfrac{1}{2}+\delta + iy\)$ we can smooth
out the curve such that it stays in the principal branch.
}
\label{arg127}
\end{figure}

At the relatively rare  points $\yplus_n$ 
where $\Re\big[\zeta\big(\tfrac{1}{2}+i\yplus_n\big)\big] < 0$,
$S(y)$ crosses one of the lines $S(y) = \pm 1$.
Taking into account  the  properties  and \ref{DeltaS} and  \ref{titchave},  
one expects  that  $S(y)$ primarily stays in the principle branch, 
i.e. it can pass to another branch,  it has to otherwise it cannot grow as $\log y$,  but it must  return to the
principal branch in order for $\langle S \rangle$ = 0. 
An example where  this
occurs  is close to the point $\yplus_{127}=282.455$. The function starts
to change branch, and as soon as it crosses the branch cut, there is a 
Riemann zero at $y_{127} = 282.465$ so $S(y)$ jumps by $1$ coming back
to the principal branch again.
This behavior is shown Figure~\ref{arg127}  and one sees that $S(y)$ 
just barely touches $-1$.
 
\begin{figure}
\centering
\begin{minipage}{.49\textwidth}
\includegraphics[width=.9\linewidth]{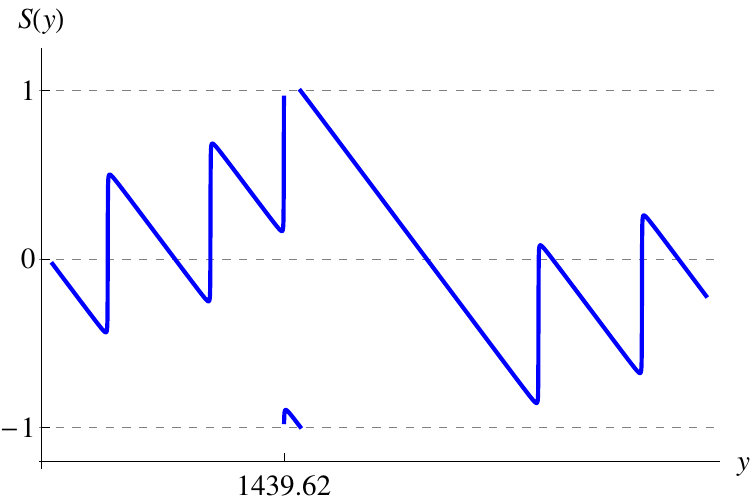}
\end{minipage}
\begin{minipage}{.49\textwidth}
\includegraphics[width=.9\linewidth]{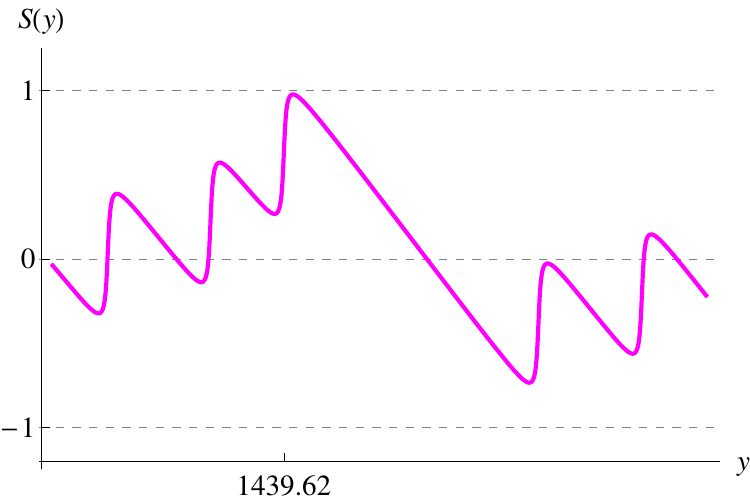}
\end{minipage}
\caption{\emph{ Left:} $S(y)$  on the principle branch in the
vicinity of $\yplus_{1017} = 1439.778$ where $S(y) >1$.
\emph{ Right:} we have included a non-zero $\delta$ to smooth out the function,
and it stays in the principal branch.
}
\label{arg1017}
\end{figure}

In Figure \ref{arg1017} (left) we  plot $S(y)$ on the principle branch  
in the vicinity of another point where $S(y)$ passes to another branch.
This time it passes to another branch while jumping at the zero
$y_{1018} = 1439.623$. The real part is negative for 
$\yplus_{1017} = 1439.778$. Since $S'(y) < 0$ it comes back to the principal
branch pretty quickly. 
The interpretation of this figure is that $S(y)$ has changed branch:   
the dangling part of the curve at the bottom should be shifted by $2$ to 
make $S(y)$ continuous.
By including a $\delta$ we can smooth out the curve to make it continuous
and to stay in the principal branch, as shown in Figure~\ref{arg1017} (right),
and this is a better rendition of the actual behavior.

\begin{figure}
\centering
\begin{minipage}{.49\textwidth}
\includegraphics[width=.9\linewidth]{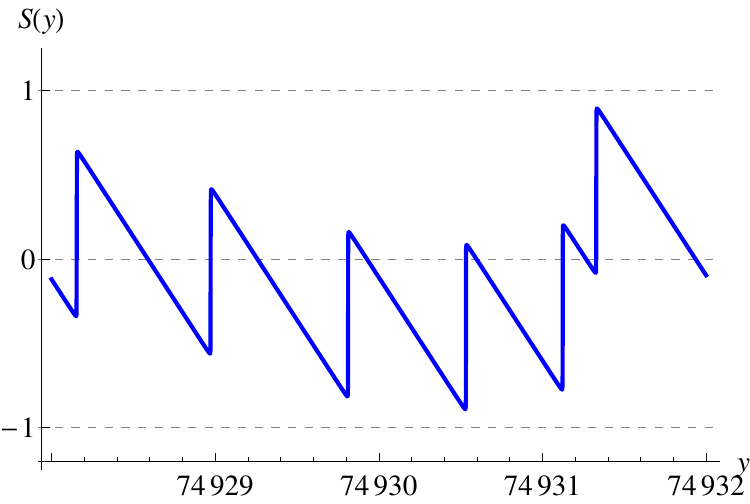}
\end{minipage}
\begin{minipage}{.49\textwidth}
\includegraphics[width=0.9\textwidth]{figs/fig7.mps}
\end{minipage}
\caption{\emph{ Left:} behavior of $S(y)$ around the $10^5$-th zero.
\emph{ Right:}  hypothetical, potentially  more  dramatic behavior of $S(y)$ although we have not observed this numerically in the range of $y$ 
we have explored.}
\label{10to5St}
\end{figure}

Note that at the rare points where $| S(y) | \ge 1$,  
it strays off the principle branch
but quickly returns to it.   
In Figure~\ref{10to5St} (left)  we  plot $S(y)$ around the  $10^5$-th zero,
and one sees that it is still on the principle branch.   

There is a well-known counter example to the RH based on the 
Davenport-Heilbronn function.    
It has a functional equation like $\zeta$, but is
known to have zeros off of the critical line.   
We will study this function in  section~\ref{sec:Davenport},
and explain how the properties of $S(y)$ described in this section 
are violated.
In short,  at a zero off of the line there is a change of branch in 
such a manner that the analog of $S(y)$ is ill defined,  
and there is thus no solution to the transcendental equation
at these points,  so the argument that $N_0 (T) = N(T)$ fails.

\section{Zeros of Dirichlet \texorpdfstring{$L$}{L}-functions}
\label{sec:dirichlet}

\subsection{Some properties of Dirichlet $L$-functions }  

We now consider the generalization of the previous results for the
$\zeta$-function to Dirichlet $L$-functions.    
The main arguments are the same as for $\zeta$,  
thus we do not repeat all of the statements in section~\ref{sec:zeta_function}. 

Much less is known about the zeros of $L$-functions in 
comparison with the $\zeta$-function,  however let us mention a few works.   
Selberg \cite{Selberg1} obtained the analog of Riemann-von Mangoldt 
counting formula \eqref{riemann_counting} for Dirichlet $L$-functions.
Based on this result, Fujii \cite{Fujii} gave an estimate for the number 
of zeros on the critical strip with the ordinate between $[T+H, T]$. 
The distribution of low lying zeros of $L$-functions near and at the 
critical line was examined in \cite{Iwaniec}, assuming the \emph{generalized
Riemann hypothesis} (GRH).
The statistics of the zeros,  i.e. the analog of the 
Montgomery-Odlyzko conjecture, were studied in \cite{Conrey3,Hughes}.
It is also known that more than half of the non-trivial zeros 
of Dirichlet $L$-functions are on the critical line \cite{Conrey4}.
For a more detailed introduction to $L$-functions see \cite{Bombieri2}.

Let us  first  introduce the basic ingredients and definitions  
regarding this class of functions, which are all well known \cite{Apostol}.
Dirichlet $L$-series are defined as 
\beq
\label{Ldef}
L(z, \chi) = \sum_{n=1}^\infty \frac{\chi (n)}{n^z} \qquad (\Re (z) > 1)
\eeq
where the arithmetic function $\chi(n)$ is a Dirichlet character.  
They enjoy an Euler product formula
\beq
L(z,\chi) =  \prod_{p} \inv{1 - \chi (p) \, p^{-z}} \qquad \(Re (z)>1\). 
\eeq
They  can all be analytically continued to the entire 
complex plane, except for a simple pole at $z=1$, and are then 
referred to as  Dirichlet $L$-functions. 

There are an infinite number of distinct Dirichlet 
characters which are primarily characterized by their modulus $k$,
which determines their periodicity. They can
be defined axiomatically,  which leads to specific properties,
some of which we now describe.
Consider a Dirichlet character $\chi$ mod $k$,  and let 
the symbol $(n,k)$ denote the greatest  common divisor of the
two integers $n$ and $k$. Then $\chi$  has the 
following properties:
\begin{enumerate}
\item $\chi(n+k) = \chi (n)$.
\item $\chi(1) = 1$ and  $\chi(0) = 0$.
\item $\chi(n  m) = \chi (n) \chi (m)$.
\item $\chi(n) = 0$ if  $(n, k) > 1$ 
and $\chi(n) \ne 0$ if $(n, k) = 1$.
\item \label{root1} If $(n, k) = 1$  then  $\chi(n)^{\varphi(k)} = 1$, where
$\varphi(k)$  is the Euler totient arithmetic function. This implies  
that $\chi(n)$ are roots of unity.
\item If $\chi$ is a Dirichlet  character so is the 
complex conjugate $\chi^*$.
\end{enumerate}
For a given modulus $k$ there are $\varphi(k)$ distinct Dirichlet 
characters,  which essentially follows from property \ref{root1} above.
They can thus be labeled as $\chi_{k,j}$ where $j= 1, 2,\dotsc, \varphi(k)$ 
denotes an arbitrary ordering. 
If $k=1$ we have the \emph{trivial} character where
$\chi(n)=1$ for every $n$, and \eqref{Ldef} reduces to the Riemann 
$\zeta$-function.
The \emph{principal } character,  usually denoted by $\chi_1$,   
is defined as 
$\chi_1(n) = 1$ if $(n,k) = 1$ and zero otherwise. In the above notation 
the principal character is always $\chi_{k,1}$.

Characters can be classified as \emph{primitive} or \emph{non-primitive}.
Consider the Gauss sum
\beq\label{tau}
G(\chi) = \sum_{m=1}^{k}\chi(m)e^{2\pi i m / k}.
\eeq
If the character $\chi$ mod $k$ is primitive, 
then 
\beq
|G(\chi)|^2 = k.
\eeq
This is no longer valid for a non-primitive character.
Consider a non-primitive character $\bar{\chi}$ 
mod $\bar{k}$. Then it can be expressed in terms of
a primitive character of smaller modulus
as $\bar{\chi}(n) = \bar{\chi_1}(n) \chi(n)$, where $\bar{\chi_1}$ is
the principal character mod $\bar{k}$ and $\chi$ is a primitive
character mod $k < \bar{k}$, where $k$ is a divisor of $\bar{k}$. More
precisely, $k$ must be the \emph{conductor} of $\bar{\chi}$ 
(see \cite{Apostol} for further details).
In this case the two $L$-functions are related as
$L(z, \bar{\chi}) = L(z,\chi) \Pi_{p | \bar{k}}\(1 - \chi(p)/p^z\)$. 
Thus $L(z,\bar{\chi})$ has the same zeros as $L(z,\chi)$.
The principal character is only primitive when $k=1$, which yields the
$\zeta$-function. The simplest example of non-primitive characters
are all the principal ones for $k \ge 2$, whose zeros are the same as 
the $\zeta$-function. Let us consider another example with $\bar{k}=6$,
where $\varphi(6)=2$,  
namely $\bar{\chi}_{6,2}$, whose components are\footnote{Our 
enumeration convention for the $j$-index of $\chi_{k,j}$ 
is taken from Mathematica.} 
\beq
\begin{tabular}{@{}c|cccccc@{}}
$n$             & $1$ & $2$ & $3$ & $4$ & $5$ & $6$ \\
\midrule[0.3pt]
$\bar{\chi}_{6,2}(n)$ & $1$ & $0$ & $0$ & $0$ & $-1$ & $0$
\end{tabular}
\eeq
In this case, the only divisors are $2$ and $3$. Since $\chi_1$ mod $2$
is non-primitive, it is excluded. We are left with $k=3$ which is the
conductor of $\bar{\chi}_{6,2}$. Then we have  
two options; $\chi_{3,1}$ which is the non-primitive 
principal character mod $3$, thus excluded, 
and $\chi_{3,2}$ which is primitive. Its components are
\beq
\begin{tabular}{@{}c|ccc@{}}
$n$             & $1$ & $2$ & $3$ \\
\midrule[0.3pt]
$\chi_{3,2}(n)$ & $1$ & $-1$ & $0$ 
\end{tabular}
\eeq
Note that $|G(\chi_{6,2})|^2 = 3 \ne 6$ and $|G(\chi_{3,2})|^2 = 3$.
In fact one can check that 
$\bar{\chi}_{6,2}(n) = \bar{\chi}_{6,1}(n) \chi_{3,2}(n)$, 
where $\bar{\chi}_{6,1}$ is the principal character 
mod $\bar{k}=6$. Thus the zeros of $L(z,\bar{\chi}_{6,2})$ are the same as 
those of $L(z,\chi_{3,2})$. Therefore, it suffices to consider 
primitive characters, and we will henceforth do so.  

We will need the functional equation satisfied by $L(z,\chi)$.   
Let $\chi$ be a \emph{primitive} character. 
Define its \emph{order} $a$ such that 
\beq\label{order}
a \equiv \begin{cases}
1 &\mbox{if $\chi(-1)=-1$ (odd),} \\
0 &\mbox{if $\chi(-1)=1$ (even).}
\end{cases} 
\eeq
Let us define the  meromorphic  function 
\beq
\label{Lambda}
\Lambda(z,\chi) \equiv \( \dfrac{k}{\pi} \)^{\tfrac{z+a}{2}}   
\Gamma\(\dfrac{z+a}{2}\) L(z, \chi).   
\eeq
Then $\Lambda$  satisfies the well known functional equation \cite{Apostol}
\beq
\label{FELambda}
\Lambda (z, \chi) =  \dfrac{ i^{-a}  \, G(\chi)}{\sqrt{k}} \,
\Lambda (1-z, \chi^*).
\eeq
The above equation is only valid for primitive characters.

\subsection{Exact equation for the \texorpdfstring{$n$}{n}-th zero}
\label{sec:dirichlet_equation}

For a primitive character, since $|G(\chi)| = \sqrt{k}$,  
the factor on the right hand side of \eqref{FELambda} is a phase. 
It is  thus possible to obtain a more symmetric form  
through a new function defined as
\beq
\label{xi}
\xi(z, \chi) \equiv \dfrac{i^{a/2} \, k^{1/4} }{\sqrt{G\(\chi\)}} \, \, 
\Lambda (z, \chi). 
\eeq
It  then satisfies 
\beq
\label{FExi}
\xi(z, \chi)  = \xi^{*}(1-z , \chi ) \equiv   \(\xi(1-z^*, \chi)\)^*.
\eeq
Above, the function $\xi^*$  of $z$ is defined as the complex conjugation of
all coefficients that  define $\xi$,  
namely  $\chi$ and the $i^{a/2}$ factor, evaluated at a non-conjugated $z$.
   
Note that $(\Lambda(z, \chi))^* = \Lambda(z^*,\chi^*)$. Using the well
known result 
\beq G\(\chi^*\) = \chi(-1) \( G(\chi) \)^*\eeq 
we conclude that \beq\label{xiconju}
\(\xi(z,\chi)\)^* = \xi\(z^*,\chi^*\).
\eeq
This implies that if the character is real, then if $\rho$ is a 
zero of $\xi$  so is $\rho^*$, and one needs only 
consider $\rho$ with  
positive imaginary part. On the other hand if $\chi \neq \chi^*$, 
then the zeros with negative imaginary part are different than $\rho^*$.
For the trivial character where $k=1$ and $a=0$, implying 
$\chi(n) = 1$ for any $n$,  then $L(z,\chi)$  reduces to the
Riemann $\zeta$-function and  \eqref{FExi} yields the well
known functional equation \eqref{chisym}.

Let  $z=x+i y$. Then 
the function \eqref{xi} can be written as
\beq
\xi(z,\chi) = A e^{i \theta}\eeq
where
\begin{align}
A(x,y,\chi) & =
\(\dfrac{k}{\pi}\)^{\tfrac{x+a}{2}}
\left|\Gamma\(\dfrac{x+a+iy}{2}\)\right| 
\left|L(x+iy,\chi)\right|, \label{Achi} \\
\theta(x,y,\chi) &= \arg \Gamma\(\dfrac{x+a+iy}{2}\) - 
\dfrac{y}{2}\log \(\dfrac{\pi}{k}\) 
- \dfrac{1}{2}\arg G(\chi) + 
\arg L(x+iy,\chi) + \dfrac{\pi a}{4}. \label{thetachi}
\end{align}
From \eqref{xiconju} we also conclude that
$A(x,y,\chi) = A(x,-y,\chi^*)$ and $\theta(x,y,\chi) = -\theta(x,-y,\chi^*)$.
Denoting 
\beq
\xi^*(1-z, \chi) = A' e^{-i \theta'}
\eeq
we have 
\beq
A'(x,y,\chi) = A(1-x,y, \chi), \qquad 
\theta'(x,y, \chi) = \theta(1-x,y,\chi).
\eeq
Taking the modulus of \eqref{FExi} we also have that
$A(x,y,\chi)=A'(x,y,\chi)$ for any $z$.

On the critical strip, the functions $L(z,\chi)$ and $\xi(z,\chi)$ 
have the same zeros.
Thus on a zero
we clearly have
\beq\label{onzero}
\lim_{z \to \rho } 
\left\{ \xi (z, \chi)  +  \xi^* (1-z,\chi ) \right\} =0.
\eeq
Let us define
\beq
\label{Bdef}
B(x,y,\chi) \equiv e^{i\theta (x,y,\chi)} + e^{- i\theta' (x,y,\chi)} .
\eeq
Since $A=A'$ everywhere, and taking separate limits in \eqref{onzero}
we therefore have
\beq\label{theta_zero}
\lim_{z'\to \rho}  \lim_{z\to \rho} A(x',y', \chi) \, B(x,y,\chi) = 0.
\eeq
Considering the $z\to \rho$ limit, a potential zero occurs when 
\beq\label{dir_bzero}
\lim_{z\to\rho} B(x,y,\chi) = 0.
\eeq
The general solution of this equation is thus given by
\beq\label{gensol_dir}
\theta + \theta' = (2n+1)\pi.
\eeq
Until now, the path to approach the zero $z\to \rho$ was not specified.
Now we put ourselves on the critical line $x=1/2$, 
and the path will be choosen as $z = \rho + \delta$ with $0 < \delta \ll 1$.
Then $\theta = \theta'$ and \eqref{gensol_dir} yields
\beq\label{critsol_dir}
\lim_{\delta\to0^+}
\theta\(\tfrac{1}{2}+\delta,y\)=\(n+\tfrac{1}{2}\)\pi.
\eeq

Let us define the function
\beq\label{RSgen}
\vartheta_{k,a} (y) \equiv  
\Im\[\log\Gamma\( \dfrac{1}{4} + \dfrac{a}{2} +i \,  \dfrac{y}{2} \)\]
- \dfrac{y}{2} \log \( \dfrac{\pi}{k} \).
\eeq
When $k=1$ and $a=0$, the function \eqref{RSgen} is just the usual 
Riemann-Siegel $\vartheta$ function \eqref{riemann_siegel}.
Thus \eqref{critsol_dir} gives the equation
\beq\label{almost_final}
\vartheta_{k,a}(y_n) + 
\lim_{\delta\to 0^{+}}\arg L\(\tfrac{1}{2}+\delta+iy_n,\chi\)
- \dfrac{1}{2}\arg G\(\chi\) + \dfrac{\pi a }{4} = 
\(n+\dfrac{1}{2}\)\pi .
\eeq
Analyzing the left hand side of \eqref{almost_final} we can see that it
has a minimum, thus we shift $n \to n - (n_0 + 1)$ for a given $n_0$, 
to label the zeros according to the convention that the first positive 
zero is labelled by $n=1$. 
Thus the upper half of the critical line will have the zeros
labelled by $n=1,2,\dotsc$ corresponding to positive $y_n$, while the 
lower half will have the negative values $y_n$ 
labelled by $n=0,-1,\dotsc$. The integer $n_0$ depends 
on $k$, $a$ and $\chi$, and should be chosen according to each specific case.
In the cases we analyze below $n_0=0$, whereas for the trivial 
character $n_0 =1$. In practice, the value of $n_0$ can always be 
determined by  plotting \eqref{almost_final} with $n=1$, 
passing all terms to its left hand side. Then it is trivial to adjust the 
integer $n_0$ such that the graph passes through the point 
$(y_1, 0)$ for the first jump, corresponding to the first positive solution.
Henceforth we will \emph{omit} the integer $n_0$ in the equations, 
since all cases analyzed in the following have $n_0=0$. 
Nevertheless, the reader
should bear in mind that for other cases, it may  be necessary
to replace $n \to n - n_0$ in the following equations.

In summary,  these zeros have the form $\rho_n = \tfrac{1}{2}+i y_n$, 
where for a given  
$n \in \mathbb{Z}$, the imaginary part $y_n$ is the solution of 
the equation
\beq\label{exact}
\vartheta_{k,a}(y_n) + 
\lim_{\delta\to0^{+}} \arg L\(\tfrac{1}{2}+\delta+iy_n,\chi\)
- \dfrac{1}{2}\arg G\(\chi\) = \(n - \dfrac{1}{2} - \dfrac{a}{4} \)\pi.
\eeq

\subsection{Asymptotic equation for the \texorpdfstring{$n$}{n}-th zero}
\label{sec:dirichlet_asymptotic}

From Stirling's formula we have the following 
asymptotic form for $y \to \pm \infty$:
\beq
\label{approxRS}
\vartheta_{k,a} (y) = \mathrm{sgn}(y)
\left( \dfrac{|y|}{2} \log\(  \dfrac{k |y|}{2 \pi e} \)  +  
\dfrac{2a-1}{8} \pi
+O(1/y) \right).
\eeq
The first order approximation of \eqref{exact}, i.e.
neglecting terms of $O(1/y)$, is given by
\beq\label{asymptotic}
\sig_n  \dfrac{|y_n|}{2\pi}\log\(\dfrac{k \, |y_n|}{2\pi e }\)
+\dfrac{1}{\pi}\lim_{\delta\to 0^{+}}
\arg L\(\tfrac{1}{2}+ \delta + i \sig_n |y_n|, \chi\)
- \dfrac{1}{2\pi}\arg G\(\chi\) = \alpha_n,
\eeq
where
\beq
\alpha_n = n + \dfrac{\sig_n -4 - 2a(1+\sig_n)}{8}.
\eeq
Above $\sig_n = 1$ if $n>0$ and $\sig_n = -1$ if $n\le 0$. For $n>0$
we have $y_n = |y_n|$ and for $n\le 0$ $y_n = - |y_n|$.

\subsection{Counting formulas}
\label{sec:dirichlet_counting}

Let us define $N^+_0 (T,\chi)$ as  the number of zeros on the critical 
line with $0< \Im( \rho) < T$ and  $N^-_0 (T,\chi)$ as the number of zeros 
with $-T < \Im( \rho)  < 0$. As explained before,
$N^+_0(T,\chi) \neq N^-_0 (T,\chi)$
if the characters are complex numbers, since  the zeros are not 
symmetrically distributed between the upper and lower half of 
the critical line. 

\begin{figure}
\centering
\includegraphics[width=0.6\linewidth]{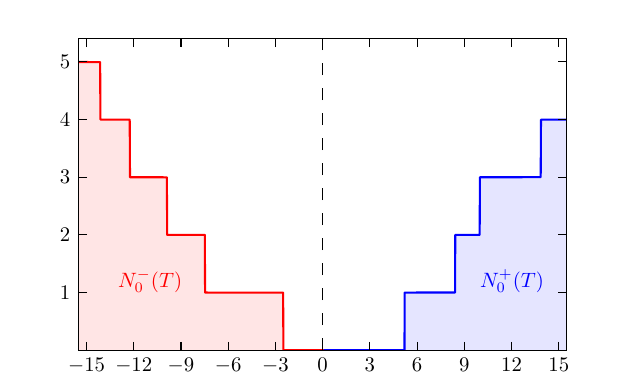}
\caption{Exact counting formulae \eqref{count_exact1}  
and \eqref{count_exact2}.
Note that they are not symmetric with respect to the origin, since the
$L$-zeros for complex $\chi$ are not complex conjugates.
We used $\chi=\chi_{7,2}$ \eqref{char72}.}
\label{fig:counting}
\end{figure}

The counting formula $N^+_0(T,\chi)$  is obtained 
from \eqref{exact} by replacing $y_n \to T$ and $n \to N^+_0 + 1/2 $, 
therefore
\beq\label{count_exact1}
N^+_0(T,\chi) = \dfrac{1}{\pi}\vartheta_{k,a}(T) + 
\dfrac{1}{\pi}\arg L\(\tfrac{1}{2}+ i T, \chi\)
- \dfrac{1}{2\pi}\arg G\(\chi\) + \dfrac{a}{4}.
\eeq
The passage from 
\eqref{exact} to \eqref{count_exact1} is justified under the assumptions
already discussed in connection with 
\eqref{counting2} and \eqref{counting2_exact}, i.e. 
assuming that \eqref{exact} has a unique solution for every $n$.
Analogously, the counting formula on the lower half line is
given by
\beq\label{count_exact2}
N_0^-(T,\chi) = \dfrac{1}{\pi}\vartheta_{k,a}(T) - 
\dfrac{1}{\pi}\arg L \(\tfrac{1}{2} - i T, \chi\)
+\dfrac{1}{2\pi}\arg G(\chi) - \dfrac{a}{4}.
\eeq
Note that in \eqref{count_exact1} and \eqref{count_exact2} $T$ is positive.
Both cases  are plotted in Figure~\ref{fig:counting} for the character 
$\chi_{7,2}$ shown in \eqref{char72}.
One can notice that they are precisely staircase functions, jumping by one
at each zero. Note also that the functions
are not symmetric about  the origin, since for a complex $\chi$ the
zeros on upper and lower half lines are not simply complex conjugates.

From \eqref{approxRS} we also have the first order approximation
for $T \to \infty$,
\beq\label{count_approx}
N^+_0(T,\chi) = \dfrac{T}{2\pi}\log\(\dfrac{k \, T}{2\pi e }\) + 
\dfrac{1}{\pi}\arg L\(\tfrac{1}{2}+i T, \chi\)
-\dfrac{1}{2\pi}\arg G\(\chi\) - \dfrac{1}{8} + \dfrac{a}{2} .
\eeq
Analogously, for the lower half line we  have
\beq\label{count_approx2}
N^-_0(T,\chi) = \dfrac{T}{2\pi}\log\(\dfrac{k \, T}{2\pi e }\) - 
\dfrac{1}{\pi}\arg L\(\tfrac{1}{2}-i T, \chi\)
+\dfrac{1}{2\pi}\arg G\(\chi\) - \dfrac{1}{8}.
\eeq
As  in \eqref{exact}, again we are omitting $n_0$ since in the cases 
below $n_0=0$, but for other cases one may need to include $\pm n_0$ on
the right hand side of $N_0^{\pm}$, respectively.

It is known that the number of zeros on the \emph{entire critical strip} 
up to height $T$, i.e. in the region
$\{0 < x < 1,\, 0 < y < T\}$, is given by \cite{Montgomery}
\beq
\label{count_strip}
N^+(T,\chi) = \dfrac{1}{\pi} \vartheta_{k,a}\(T\)
+\dfrac{1}{\pi}\arg L\(\tfrac{1}{2}+iT,\chi\) -
\dfrac{1}{\pi}\arg L\(\tfrac{1}{2},\chi\).
\eeq
This formula follows from a straightforward generalization of the
method shown in Appendix~\ref{sec:number_zeros} for the $\zeta$-function.
From Stirling's approximation and using 
\beq
2a-1=-\chi(-1),
\eeq
for $T\to\infty $ we obtain the asymptotic 
approximation \cite{Selberg1,Montgomery}
\beq
\label{selberg_counting}
N^+(T,\chi) = \dfrac{T}{2\pi}\log\(\dfrac{k \, T}{2\pi e} \)
+ \dfrac{1}{\pi}\arg L\(\tfrac{1}{2} + i T,\chi\) 
- \dfrac{1}{\pi}\arg L\(\tfrac{1}{2},\chi\) 
- \dfrac{\chi(-1)}{8} + O(1/T).
\eeq
Both formulas \eqref{count_strip} and \eqref{selberg_counting} are
exactly the same as \eqref{count_exact1} and \eqref{count_approx}, 
respectively. This can be seen as follows. 
From \eqref{FExi} we conclude that $\xi$ is real on the critical line.
Thus 
\beq
\arg\xi\(\tfrac{1}{2}\)=0=
-\dfrac{1}{2}\arg G\(\chi\)+\arg L\(\tfrac{1}{2},\chi\)+\dfrac{\pi a}{4}.
\eeq
Then, replacing $\arg G$ in \eqref{exact} we obtain
\beq\label{exact_v2}
\vartheta_{k,a}\(y_n\) + 
\lim_{\delta\to0^+}\arg L\(\tfrac{1}{2}+\delta+iy_n, \chi\)
-\arg L\(\tfrac{1}{2},\chi\) = \(n -\tfrac{1}{2}\)\pi.
\eeq
Replacing $y_n\to T$ and $n\to N_0^+ + 1/2$ in 
\eqref{exact_v2} we have precisely the expression
\eqref{count_strip}, and also \eqref{selberg_counting} for $T\to \infty$.
Therefore, we conclude that $N_0^+(T,\chi) = N^+(T,\chi)$ exactly.
From \eqref{xiconju} we see that negative zeros for character
$\chi$ correspond to positive zeros for character $\chi^*$.
Then for $-T < \Im(\rho) < 0$ the counting on the strip
also coincides  with the counting on the line, since
$N_0^-(T,\chi) = N_0^+(T,\chi^*)$ and $N^-(T,\chi) = N^+(T,\chi^*)$.
Therefore, the number of zeros on the whole \emph{critical strip}
is the same as the number of zeros on the \emph{critical line} obtained
as solutions of \eqref{exact}. This is valid under the assumption that
\eqref{exact} has a unique solution for every $n$.

\section{Zeros of \texorpdfstring{$L$}{L}-functions based on modular forms}
\label{sec:modular}

Let us  generalize the previous results to $L$-functions based on level 
one modular forms.   We first recall some basic definitions and properties.  
The \emph{modular group} can be represented by the set
of $2 \times 2$ integer matrices
\beq \label{SL2Z}
\SL =  \left\{ 
A = \begin{pmatrix}
a & b \\
c & d 
\end{pmatrix}  \, \big\vert \, \,  a,b,c,d \in \mathbb{Z}, \ \det A = 1
\right\},
\eeq
provided each matrix $A$ is identified with $-A$, i.e. $\pm A$ are regarded
as the same transformation. Thus for $\tau$ in the upper half 
complex plane, it transforms as 
$$\tau  \mapsto A \tau = \frac{a \tau + b}{c\tau + d}$$ 
under the action of the modular group.
A \emph{modular form} $f$ of weight $k$ is a function that is analytic in the 
upper half complex plane which satisfies the functional relation 
\cite{Apostolmodular} 
\beq \label{fform}
f\(\dfrac{a\tau  + b}{c\tau + d} \) =  \(c\tau +d\)^k   f(\tau).
\eeq
If the above equation is satisfied for all of $\SL$,  
then $f$ is referred to as being of \emph{level one}.
It is possible to define higher level modular forms  which satisfy 
the above equation for a subgroup of $\SL$.   Since our results are easily
generalized to the higher level case,  henceforth we will  only consider 
level one forms.     

For the $\SL$ element 
$\big( \begin{smallmatrix} 1 & 1 \\ 0 & 1 \end{smallmatrix} \big)$, 
the above implies  the periodicity $f(\tau ) = f(\tau+1)$, 
thus it has a Fourier series
\beq \label{Fourier}
f(\tau ) = \sum_{n=0}^\infty  a_f (n) \, q^n  ,  
\qquad q \equiv  e^{2 \pi i \tau}.
\eeq
If $a_f(0) = 0$ then $f$ is called a \emph{cusp form}. 

From the Fourier coefficients, one can define the Dirichlet series
\beq \label{Lmod}
L_f \(z\) = \sum_{n=1}^\infty \dfrac{a_f \(n\)}{n^z}.
\eeq
The functional equation relates  $L_f\(z\)$ to 
$L_f\(k-z\)$,  so that the critical line is $\Re (z) = k/2$, 
where $k \ge 4$ is an even integer.   One can always shift the critical line
to $1/2$ by  replacing $a_f (n) \to a_f (n)/n^{(k-1)/2}$,  
however we will not do this here.   
Let us define 
\beq \label{ximod}
\Lambda_f(z) \equiv \(2\pi\)^{-z} \, \Gamma \( z \) \, L_f (z).
\eeq
Then the functional equation is given by \cite{Apostolmodular}
\beq \label{LFE}
\Lambda_f(z) = (-1)^{k/2} \Lambda_f(k-z).
\eeq

There are only two cases to consider since $k/2$ can be
an even or an odd  integer.  As in \eqref{xi} we can absorb
the extra minus sign  factor for the odd case. 
Thus we define $\xi_f(z) \equiv \Lambda_f(z)$ for
$k/2$ even, and  we have $\xi_f(z) = \xi_f(k-z)$,  
and $\xi_f(z) \equiv e^{-i \pi/2}\Lambda_f(z)$ for $k/2$ odd, 
implying $\xi_f(z) = \xi_f^*(k-z)$. 
Representing $\xi_f(z) = |\xi_f| \, e^{i \theta }$ where 
$z=x+i y$, we follow exactly the same steps as in the previous sections.
From the solution \eqref{gensol_dir} we conclude
that there are infinite zeros
on the critical line $\Re(\rho) = k/2$ determined by
$\lim_{\delta\to0^+}\theta\(\tfrac{k}{2}+\delta,y, \chi\)=
\(n-\tfrac{1}{2}\)\pi$.
Therefore, these zeros have
the form $\rho_n = \tfrac{k}{2} + iy_n$, where $y_n$ is the solution
of the equation
\beq \label{exactmod}
\vartheta_k (y_n) + \lim_{\delta \to 0^+} 
\arg L_f \( \tfrac{k}{2} + \delta + i y_n \) =  
\( n - \dfrac{1 + (-1)^{k/2}}{4} \) \pi 
\eeq
for $n=1,2,\dotsc$, and we have defined
\beq \label{RSL}
\vartheta_k (y)  \equiv \Im \[ \log \Gamma\( \tfrac{k}{2} + i y \)\] 
- y \log 2\pi.
\eeq
This implies that the number of solutions of \eqref{exactmod} with 
$0< y < T$ is given by
\beq \label{NTmod}
N_0\(T\) =  \inv{\pi} \vartheta_k (T) + 
\inv{\pi} \arg L_f \( \tfrac{ k}{2} + iT \) - \dfrac{1 - (-1)^{k/2}}{4}.
\eeq
In the limit of large $y_n$, neglecting terms of $O(1/y)$,
the equation \eqref{exactmod} becomes 
\beq \label{asymod}
y_n \log \(  \frac{y_n}{2 \pi e} \)  + 
\lim_{\delta \to 0^+}  \arg L_f \( \tfrac{k}{2} +\delta + i y_n\) 
= \( n - \dfrac{k+ (-1)^{k/2}}{4} \) \pi.
\eeq

\section{The Lambert \texorpdfstring{$W$}{W}-function}

The following section,  and even more, section \ref{sec:saddle}, 
will involve the Lambert-$W$ function,  thus we review its most
important properties in this section.   
The basic facts about the $W$-function that we  present, 
including some history, are essentially based on \cite{Corless}, 
where the reader can also find more details.

In $1758$ Lambert solved the equation 
$x=q+x^m$,
expressing $x$, and also powers $x^\alpha$, as power series in $q$. 
A few years later Euler considered a more 
symmetric version of this equation through the substitution
$x \to x^{-\beta}$, 
$q \to (\alpha - \beta) v$ and
$m \to \tfrac{\alpha}{\beta}$. Then taking the limit $\beta \to \alpha$
one obtains $\log x = v x^{\alpha}$. This equation can be written
in the form $\log x  = v x$ through $x^\alpha \to x$ and $\alpha v \to v$.
Thus, exponentiating the last equation and introducing 
the variables $z = -v$ and $W(z) = -v x$,
we obtain the equation
\beq\label{W_def_eq}
W(z) \, e^{W(z)} = z.
\eeq
$W(z)$ is called the \emph{Lambert function}, and
\eqref{W_def_eq} is its defining equation. 
Although its power series approximation was consider for the first 
time by Lambert and Euler, this function only started to be effectively
studied during the past $20$ years. The $W$-function should be considered
as a new elementary function, since it cannot be expressed in terms
of the other known elementary functions.

Johan Heinrich Lambert was born in Mulhouse (a French city
 bordering with Switzerland and Germany) in $1728$, and died in 
Berlin in $1777$. Lambert was a self educated mathematician
having a broad range of scientific interests. 
He contributed to number theory, 
geometry, statistics, astronomy, cosmology and philosophy, just to name
a few areas. He is responsible for the modern notation of
hyperbolic functions, and was the first one to prove 
the irrationality of $\pi$. When Euler considered the solution
to equation $\log x = v x$, he gave credit to Lambert. The translation
of his paper's title is ``On a series of Lambert and some of its significant
properties''. It seems that Euler learned about Lambert's result
from a conversation between both when Lambert was travelling 
from Z\" urich to Berlin. Euler described his excitement in a letter 
to Goldbach in 1764.

Leaving history aside, let us now consider the $W$-function 
defined through \eqref{W_def_eq}
over a complex field, $z \in \mathbb{C}$.
$W(z)$ is not a single-valued function, so we need to introduce its
branch structure. This is done following the same branch structure
of $\log z$, that we now recall. If $u = \log z$ then for every $z$ 
we have $e^{u} = z$. However, 
note that $u$ is not uniquely defined since $u \to u + 2\pi k i$
for $k \in \mathbb{Z}$ will give the same $z$. Thus it is necessary
to divide the complex $u$-plane in regions which are single-valued related
to the $z$-plane. Each of these regions is called a \emph{branch} and
is labelled by $k$. The way one partitions the $u$-plane is a convention.
The most convenient way is to define the $k$-th branch as the
region on the $u$-plane limited by
\beq\label{log_branch}
(2k - 1)\pi < \Im(u) \le (2k + 1)\pi, \qquad k\in\mathbb{Z}.
\eeq
Each boundary
in \eqref{log_branch} is mapped on the negative real line $(-\infty, 0]$
of the $z$-plane. The line $(-\infty, 0]$ is called the \emph{branch cut}
and $z=0$ is the \emph{branch point}. 
We adopt the counterclockwise direction such that
the branch cut closes on top, i.e. for $z=re^{i\theta}$, 
$-\pi < \theta \le \pi$ then $\log z$, for each branch $k$, is a continuous 
function of $z$.

\begin{figure}
\begin{center}
\begin{minipage}{.5\textwidth}
\includegraphics[width=.95\linewidth]{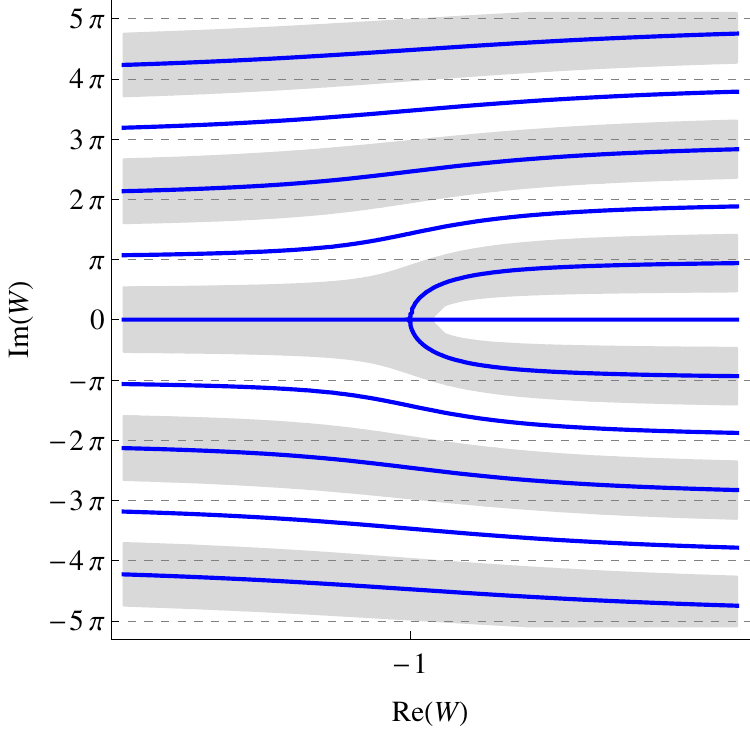}
\end{minipage}%
\begin{minipage}{.5\textwidth}
\includegraphics[width=.95\linewidth]{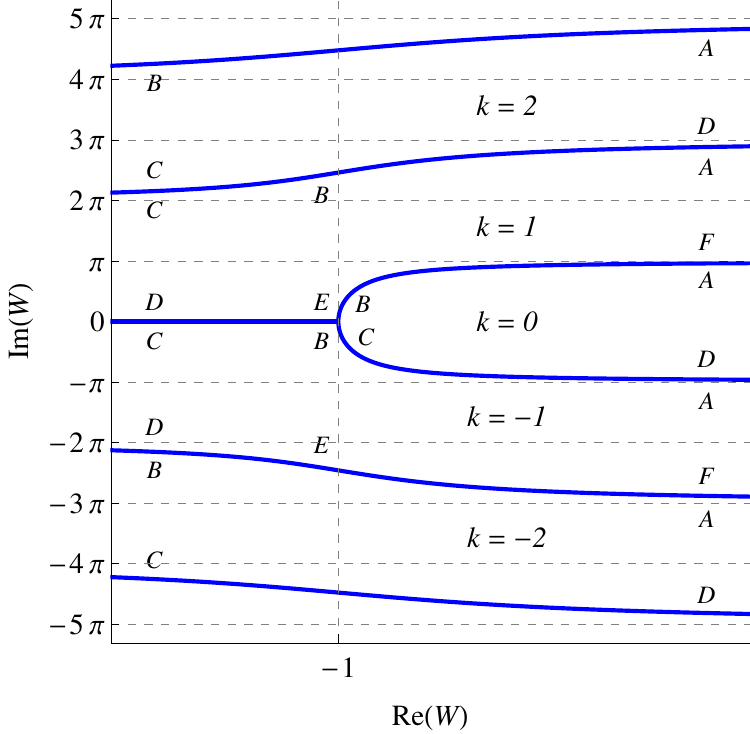}
\end{minipage}
\caption{\emph{ Left:} the solid (blue) lines correspond to the condition
$y=0$ in \eqref{coordW} and the shaded (gray) regions to
$x \le 0$. \emph{ Right:} the boundaries of each branch are given by the
lines which are contained in the shaded regions. 
This is the branch structure of
$W(z)$ viewed from the $W$-plane.
}
\label{fig:W_branches}
\end{center}
\end{figure}

Now let us turn back to $W(z)$. We refer to the $W$-plane and the $z$-plane,
and denote
\beq
W = u + iv, \qquad z = x + iy.
\eeq
From the defining equation \eqref{W_def_eq} we have
\beq\label{coordW}
x = e^u\( u\cos v - v\sin v \), \qquad y = e^{u}\( u\sin v + v\cos v \).
\eeq
We require the branch cut of $W(z)$ to be similar to $\log z$, thus we
define the branch cut on the $z$-plane to be the line $(-\infty, 0]$.
Imposing $y=0$ and $x \le 0$ on relations \eqref{coordW} we obtain
the boundaries and regions shown in Figure~\ref{fig:W_branches} (left). 
Thus the boundaries of each branch are given by the curves which
are contained in the shaded regions, which are shown
in Figure~\ref{fig:W_branches} (right). Each of these lines are mapped
to the branch cut $(-\infty, 0]$ in the $z$-plane. For each branch $k$  
the lower boundary is open, and not included, while the upper boundary
is closed, and included, in its defining region. 
From now on we denote $W(z)$ by $W_k(z)$ when referring to the specific $k$-th
branch. Note that the boundary lines are asymptotic to $W= n \pi i$.

\begin{figure}
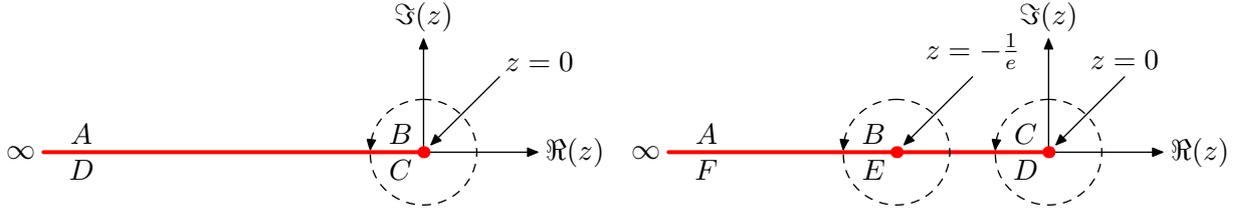

\begin{center}
\begin{minipage}{.5\textwidth}
\includegraphics{figs/fig10.mps}
\end{minipage}%
\begin{minipage}{.5\textwidth}
\includegraphics{figs/fig11.mps}
\end{minipage}
\caption{\emph{ Left:} Mapping from $W_k$ 
for $k=\pm 2, \pm 3, \dotsc$ 
(Figure~\ref{fig:W_branches} (right))
into the $z$-plane. Note how
the points $A, \, B, \, C\, $ and $D$ are mapped.
We close the contour on top, i.e. $z=re^{i\theta}$ 
where $-\pi <\theta\le \pi$, and $W_k(z)$ is a continuous
single-valued function.
\emph{ Right:} The branch cuts for $W_{\pm 1}$.
}
\label{fig:W2_branch}
\end{center}
\end{figure}

The branches $W_k$ for $k=\pm 2, \pm 3, \dotsc$ in 
Figure~\ref{fig:W_branches} (right)
are mapped into the $z$-plane according to Figure~\ref{fig:W2_branch} (left). 
The curve in the $W$-plane separating $W_k$ 
from $W_{k+1}$ for $k=1,2,\dotsc$ is given by
$\left\{ -v\cot v + iv \ | \ 2 k \pi < v < (2k+1)\pi \right\}$,
and the curve separating $W_k$ from $W_{k-1}$ for $k=-1,-2,\dotsc$ is given by
$\left\{ -v\cot v + iv \ | \ (2 k-1)\pi < v < 2k\pi \right\}$.

The branch structure for $W_0$ and $W_{\pm 1}$ is different from
above.
The curves shown in Figure~\ref{fig:W_branches} (right), 
separating $W_0$ from $W_{\pm 1}$ 
are given by
$\left\{ -v\cot v + iv \ | \ -\pi < v < \pi \right\}$.
The point $W=-1$, corresponding to $z=-1/e$, is a double branch point,
linking $W_0$ and $W_{\pm 1}$. 
The curve separating $W_1$ and $W_{-1}$ is $(-\infty, -1]$.
Note that $W_1$ does not include the real line, since the lower boundary
is not included in its region. The branches $W_0$ and $W_{-1}$ are 
special in the sense that they are the only ones which include part of 
the real line. Thus only $W_{0}$ and $W_{-1}$ can yield real values.

The branches $W_1$ and $W_{-1}$ have a double branch cut in
the $z$-plane. One is the
line $(-\infty, 0]$ and the other is $(-\infty, -1/e]$. Both branches
are mapped into the $z$-plane
according to Figure~\ref{fig:W2_branch} (right). Note how the
points $A, \, B,\, \dotsc F$ are related.

\begin{figure}
\begin{center}
\includegraphics[width=.5\textwidth]{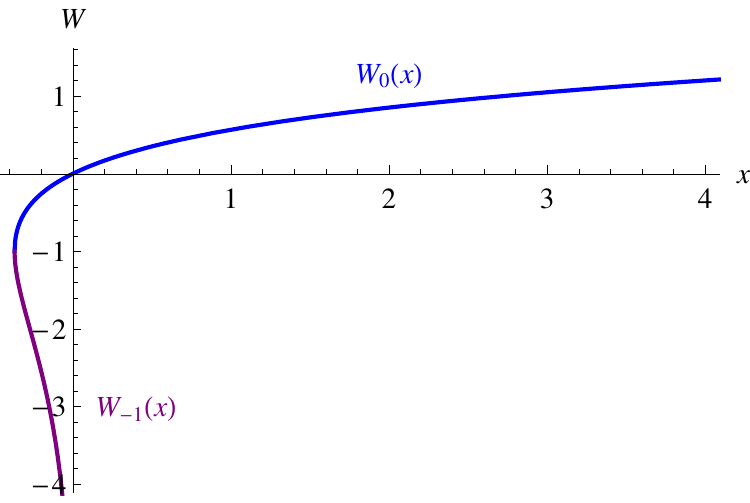}
\caption{The two branches $W_0(x)$ and $W_{-1}(x)$ for real $x$.
These are the only branches giving real values. The domain of
$W_0(x)$ is $x\in[-1/e,\infty)$ and the domain of $W_{-1}(x)$ is 
$x\in[-1/e,0)$.
}
\label{fig:W_real}
\end{center}
\end{figure}

The branch $W_0$ has only one branch cut which is $(-\infty, -1/e]$.
It is exactly like in Figure~\ref{fig:W2_branch} (left) but with
the point $z=0$ replaced by $z=-1/e$.
For real values of $z \in [-1/e,\infty)$ we see from \eqref{W_def_eq} 
that $W_0(-1/e) = -1$, $W_0(0) = 0$ and
$W_0(\infty) = \infty$. For the branch $W_{-1}$ we see that the function
is real for $z \in [-1/e, 0)$, having the values $W_{-1}(-1/e)=-1$ and
$W_{-1}(0) = -\infty$. Thus for real $z$ we have the picture shown
in Figure~\ref{fig:W_real}. Sometimes the principal branch $W_0$ is
denoted simply by $W$ for short, when there is no chance of confusion.

Now it only remains to show that each branch maps bijectively into 
the $z$-plane.
This will be true provided the
Jacobian determinant of the transformation \eqref{coordW} does
not vanish. We therefore have
\beq
J = \dfrac{\partial(x,y)}{\partial(u,v)} = 
e^u \begin{pmatrix}
\cos v + u\cos v - v\sin v  &&  - v\cos v - \sin v -u \sin v \\
\sin v + v\cos v + u\sin v  &&  \cos v + u \cos v -v \sin v
\end{pmatrix}.
\eeq
Then, $\det J = e^{2u}\(v^2 + (1+u)^2\)$. This can only be zero if
$u=-1$ and $v=0$, which correspond exactly to the double branch point.
Thus for every branch $W_k$ we have $\det J > 0$.

Except on the branch cut we have the symmetry
\beq
\(W_k(z)\)^* = W_{-k}\(z^*\).
\eeq
This can easily be seen from Figure~\ref{fig:W_branches} (right), where
the branches are symmetric by $v \to -v$, and analyzing how the points
are mapped into the $z$-plane by complex conjugation on
the $W$-plane.

\section{Approximate zeros in terms of the Lambert 
\texorpdfstring{$W$}{W}-function}
\label{sec:lambert}

\subsection{Explicit formula}

We now show that it is possible to obtain an approximate solution to
the previous transcendental equations with an explicit formula.
Let us  start with the zeros of the $\zeta$-function, described
by equation \eqref{FinalTranscendence}. Consider its leading order 
approximation, or equivalently its average since
$\langle \arg\zeta\(\tfrac{1}{2}+iy\) \rangle = 0$. Then we have
the transcendental equation
\beq \label{ApproxTranscendence}
\dfrac{\tilde{y}_n}{2\pi}\log\(\dfrac{\tilde{y}_n}{2\pi e}\) = 
n - \dfrac{11}{8}.
\eeq
Through the transformation 
\beq
\tilde{y}_n = 2\pi\(n-\dfrac{11}{8}\) \dfrac{1}{x_n}
\eeq
the equation \eqref{ApproxTranscendence} can be written as 
\beq
x_n e^{x_n} = \dfrac{n-\tfrac{11}{8}}{e}.
\eeq
Comparing with 
\eqref{W_def_eq} we thus we obtain 
\beq \label{Lambert}
\tilde{y}_n = 
\dfrac{2\pi\(n-\tfrac{11}{8}\)}{W\[e^{-1}\(n-\tfrac{11}{8}\)\]}
\eeq
for $n=1,2,\dotsc$, where above $W$ denotes the principal branch $W_0$.

Although the inversion from \eqref{ApproxTranscendence} to \eqref{Lambert}
is rather simple, it is very  convenient since  it is indeed an 
explicit formula depending only on $n$,  
and $W$ is included in most numerical packages.
It gives an approximate solution for the ordinates of the Riemann zeros
in closed form. The values computed from \eqref{Lambert} are much closer 
to the Riemann zeros than Gram points,  
and one does not have to deal with violations of 
Gram's law (see below).     

Analogously, for Dirichlet $L$-functions, after neglecting
the $\arg L$ term, the equation \eqref{asymptotic} yields a transcendental
equation which can be solved explicitly as
\beq\label{approx}
\tilde{y}_n = \dfrac{2\pi \sig_n A_{n}\(\chi\)}{
W\[ k \, e^{-1} A_{n}\(\chi\) \]} 
\eeq
for $n=0,\pm1,\pm2,\dotsc$ and where
\beq
A_{n}\(\chi\) = \sig_n \(n + \dfrac{1}{2\pi}\arg G(\chi)\)
+\dfrac{1 - 4\sig_n - 2a\(\sig_n + 1\)}{8}.
\eeq
In the above formula $n=1,2,\dotsc$ correspond to positive $y_n$
solutions,  while $n=0,-1,\dotsc$ correspond to negative $y_n$ solutions.
Contrary to the $\zeta$-function, in general, the zeros are not
conjugate related along the critical line.

In the same way, ignoring the small $\arg L_f$ term in \eqref{asymod},
the approximate solution for the imaginary part of the zeros of $L$-functions 
based on level one modular forms is given by
\beq \label{yLambertmod}
\tilde{y}_n = \dfrac{A_n \pi }{W\[ (2e)^{-1}A_n \]} 
\eeq
where $n=1,2,\dotsc$ and
\beq
A_n = n - \dfrac{k+(-1)^{k/2}}{4}. 
\eeq

\subsection{Further remarks}
\label{sec:lambert_remarks}

Let us focus on the approximation \eqref{Lambert} 
regarding zeros of the $\zeta$-function.
Obviously the same arguments apply to the zeros of the other classes of 
functions, based on formulas \eqref{approx} and  \eqref{yLambertmod}.

The estimates given by \eqref{Lambert} can be calculated to high accuracy for 
arbitrarily large $n$, since $W$ is a standard elementary function.  
Of course, the  $\tilde{y}_n$ are not as accurate as the solutions 
$y_n$ including
the $\arg \zeta$ term, as we will see in section \ref{sec:numerical}.
Nevertheless, it is indeed a good estimate, especially if one considers 
very high zeros,  where  traditional methods have 
not previously estimated  such high values.
For instance, formula \eqref{Lambert} can easily estimate the 
zeros shown in Table~\ref{highn}, and much higher if desirable.
\begin{table}
\small
\def\arraystretch{1.2}
\centering
\begin{tabular}{@{}ll@{}}
\toprule[1pt]
$n$ & $\tilde{y}_n$ \\ 
\midrule[1pt]
$10^{22}+1$ &  
$1.370919909931995308226770\times 10^{21}$ \\
$10^{50}$ &  
$5.741532903784313725642221053588442131126693322343461\times 10^{48}$ \\
$10^{100}$  & 
$2.80690383842894069903195445838256400084548030162846045192360059224930$
\\[-0.5ex] 
&$922349073043060335653109252473234\times 10^{98}$ \\
$10^{200}$ &
$1.38579222214678934084546680546715919012340245153870708183286835248393$
\\[-0.5ex]
&$8909689796343076797639408172610028651791994879400728026863298840958091$
\\[-0.5ex]
&$288304951600695814960962282888090054696215023267048447330585768
\times10^{198}$ \\
\bottomrule[1pt]
\end{tabular}
\caption{Formula \eqref{Lambert} can easily estimate very high Riemann zeros.
The results are expected to be correct up to the decimal point,  i.e. to the 
number of digits in the integer part.   
The numbers are shown with three digits beyond the integer part.}
\label{highn}
\end{table}

The numbers in this table are accurate approximations to the $n$-th zero 
to the number of digits shown,  which is approximately the number of 
digits in the integer part.
For instance,  the approximation to the $10^{100}$ zero is
correct to $100$ digits.   
With Mathematica we easily calculated the first million digits
of the $10^{10^6}$ zero.\footnote{The result is $200$ pages long 
and available at   
\url{http://www.lepp.cornell.edu/~leclair/10106zero.pdf}.}

Using the asymptotic behaviour $W(x) \sim \log x$ for large $x$,
the $n$-th zero is approximately $\tilde{y}_n \approx  2 \pi n/ \log n$,
as already known \cite{Titchmarsh}.
The distance between consecutive zeros is $2\pi/\log n$, which
tends to zero when $n\to \infty$.

The solutions \eqref{Lambert} are reminiscent of the so-called Gram 
points $g_n$,  which are solutions to 
$\vartheta (g_n) = n\pi$ where $\vartheta$ is given by \eqref{riemann_siegel}. 
Gram's law is the tendency for Riemann zeros to lie between 
consecutive Gram points, but it is known to fail for about $1/4$ 
of all 
Gram intervals. Our $\tilde{y}_n$ are intrinsically different from 
Gram points. It is an approximate solution for the ordinate of
the zero itself.
In particular,  the  Gram point $g_0 = 17.8455$ is the closest to 
the first Riemann zero, whereas $\tilde{y}_1 = 14.52$ is 
already much closer to the true zero which is $y_1 = 14.1347\dotsc$.
The traditional method to compute the zeros is based on 
the Riemann-Siegel formula, 
$\zeta\(\tfrac{1}{2}+iy\) = Z(y)\(\cos\vartheta(y) - i \sin\vartheta(y)\)$,
and the empirical observation that the real part of this equation 
is almost always positive, except when Gram's law fails, and $Z(y)$ has 
the opposite
sign of $\sin\vartheta$. Since $Z(y)$ and $\zeta\(\tfrac{1}{2}+iy\)$ have
the same zeros, one looks for the zeros of $Z(y)$ between two
Gram points, as long as Gram's law holds $(-1)^nZ\(g_n\)>0$. 
To verify the RH numerically, the counting formula \eqref{backlund} must 
also be used, to assure that the number of zeros on the critical line 
coincide with the number of zeros on the strip. 
The detailed procedure is throughly 
explained in \cite{Edwards,Titchmarsh}.
Based on this method, amazingly accurate solutions and high zeros on 
the critical line were computed 
\cite{OdlyzkoSchonhage,Odlyzko,Odlyzko2,Gourdon}.
Nevertheless, our proposal is \emph{fundamentally} different.
We claim that \eqref{exact_eq2}, or its asymptotic 
approximation \eqref{FinalTranscendence}, is the equation that determines
the Riemann zeros on the critical line.
Then, one just needs to find its solution for a given
$n$. We will compute the Riemann zeros  in this way in the next 
section, just by solving the equation numerically, starting
from the approximation given by the explicit formula \eqref{Lambert},
without using Gram points nor the Riemann-Siegel $Z$ function.
Let us emphasize that our goal is not to provide a more efficient algorithm
to compute the zeros \cite{OdlyzkoSchonhage},  although the method 
described here may very well be,   but to justify the 
validity of equations \eqref{FinalTranscendence} and \eqref{exact_eq2}.

\section{Numerical analysis: \texorpdfstring{$\zeta$}{Zeta}-function}
\label{sec:numerical}

Instead of solving the exact equation \eqref{exact_eq2} we will 
initially consider its first order approximation, 
which is equation \eqref{FinalTranscendence}. 
As we will see, this approximation already yields surprisingly accurate 
values for the Riemann zeros.    These approximate solutions will 
be used to study the GUE statistics and prime number counting function in
the next two sections.   

\begin{figure}
\begin{center}
\begin{minipage}{.5\textwidth}
\includegraphics[width=1\linewidth]{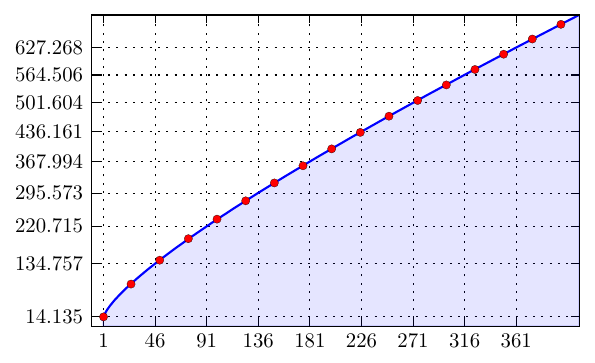}
\end{minipage}%
\begin{minipage}{.5\textwidth}
\includegraphics[width=1\linewidth]{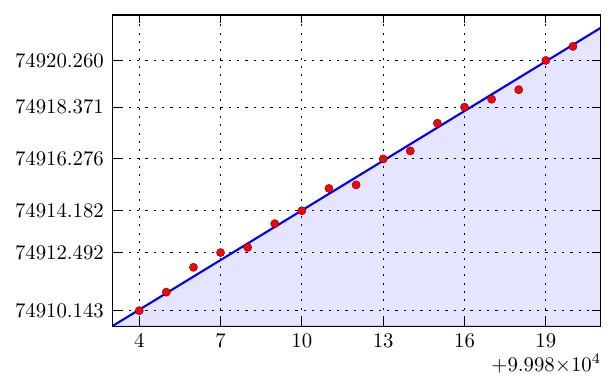}
\end{minipage}
\caption{
Comparison of \eqref{Lambert} (blue line) and 
\eqref{FinalTranscendence} (red dots). We plot $y_n$ against $n$.
\emph{ Left:} here $n \in [1,\dotsc,400]$. 
\emph{ Right:} if we focus on a small range we can see
the solutions of \eqref{FinalTranscendence} oscillating around the line  
\eqref{Lambert} due to the fluctuating $\arg \zeta$ term. Here
$n \in [99984,\dotsc,10^5]$.
}
\label{fig:trans_zeros}
\end{center}
\end{figure}

Let us first consider how the approximate solution given by
\eqref{Lambert} is modified by the presence of the $\arg\zeta$ term
in \eqref{FinalTranscendence}.  
 Based on the discussion in section~\ref{Soft},   up to rather high $y <10^{10}$,  
$S(y)$ is almost always on the principal branch, thus we numerically compute $\arg\zeta$  by taking
its principal value.  This is not necessarily correct given the growth of $S(y)$,   thus we check it is indeed a zero to a high
number of digits by verifying $|L(\rho)| = 0$ to this precision. 
 Numerically this prescription  of taking the $\Arg$ rather than the $\arg$ seems to be better near a zero $y_n$ than for typical $y$.\footnote{We remind the reader that $\arg \zeta = \Arg \zeta  + 2\pi n$,  with $-\pi < \Arg \zeta <\pi $,  where $n$ is an integer.}
As already discussed, the function 
$\arg\zeta\(\tfrac{1}{2} + i y\)$ oscillates around zero, as 
shown in Figure~\ref{fig:arg_counting} (left).
At a zero it can be well-defined by the limit \eqref{deltadef}.  For example, for the first Riemann zero 
$y_1 = 14.1347\dotsc$, 
\beq
\lim_{\delta \to 0^+} \arg \zeta \( \tfrac{1}{2}+\delta +  i y_1 \) = 
0.157873919880941213041945.....
\eeq
The $\arg \zeta$ term plays an essential  role and 
indeed improves the estimate of the $n$-th zero.
This can be seen in Figure~\ref{fig:trans_zeros}, where we compare the
estimate given by \eqref{Lambert} with the numerical solutions of
\eqref{FinalTranscendence}.

Since equation \eqref{FinalTranscendence} typically alternates in sign 
around a zero, we can apply a root finder method in an appropriate interval, 
centered around the approximate solution $\tilde{y}_n$ given by 
formula \eqref{Lambert}.  
Some of the solutions obtained in this way are presented in 
Table~\ref{some_zeros} (left), and are accurate up to the number of 
decimal places shown. 
We used  only Mathematica or some very simple algorithms to perform 
these  numerical computations, taken from standard open source numerical 
libraries. 

\begin{table}
\footnotesize
\def\arraystretch{1.2}
\centering
\begin{minipage}{.49\textwidth}
\begin{tabular}{@{}lrr@{}}
\toprule[1pt]
$n$ & $\tilde{y}_n$ & $y_n$ \\ 
\midrule[1pt] 
$1$       &         $14.52$ &         $14.134725142$ \\
$10$      &         $50.23$ &         $49.773832478$ \\
$10^2$    &        $235.99$ &        $236.524229666$ \\
$10^3$    &       $1419.52$ &       $1419.422480946$ \\
$10^4$    &       $9877.63$ &       $9877.782654006$ \\
$10^5$    &      $74920.89$ &      $74920.827498994$ \\
$10^6$    &     $600269.64$ &     $600269.677012445$ \\
$10^7$    &    $4992381.11$ &    $4992381.014003179$ \\
$10^8$    &   $42653549.77$ &   $42653549.760951554$ \\
$10^9$    &  $371870204.05$ &  $371870203.837028053$ \\
$10^{10}$ & $3293531632.26$ & $3293531632.397136704$ \\
\bottomrule[1pt]
\end{tabular}
\end{minipage}
\begin{minipage}{.49\textwidth}
\begin{tabular}{@{}ll@{}}
\toprule[1pt]
$n$ & $y_n$ \\ 
\midrule[1pt] 
$1$ & $14.13472514173469379045725198356247$ \\
$2$ & $21.02203963877155499262847959389690$ \\
$3$ & $25.01085758014568876321379099256282$ \\
$4$ & $30.42487612585951321031189753058409$ \\
$5$ & $32.93506158773918969066236896407490$ \\
$6$ & $37.58617815882567125721776348070533$ \\
$7$ & $40.91871901214749518739812691463325$ \\
$8$ & $43.32707328091499951949612216540681$ \\
$9$ & $48.00515088116715972794247274942752$ \\
$10$& $49.77383247767230218191678467856372$ \\
$11$& $52.97032147771446064414729660888099$ \\
\bottomrule[1pt]
\end{tabular}
\end{minipage}
\caption{Numerical solutions of equation \eqref{FinalTranscendence}.
\emph{ Left: } solutions accurate up to the $9$-th decimal place 
and agree with \cite{Odlyzko,Oliveira}.
\emph{ Right:} although it  was derived for high $y$,
it provides accurate numbers even for the lower zeros.
}
\label{some_zeros}
\end{table}

Although the formula for $y_n$ was solved  for large $n$, it is 
surprisingly accurate even for the lower zeros, as shown in
Table~\ref{some_zeros} (right).   
It is actually easier to solve numerically for low 
zeros since  $\arg \zeta $ is better behaved.    These numbers are correct 
up to the number of digits shown, and the precision was improved simply 
by decreasing the error tolerance.

\section{GUE statistics}
\label{sec:GUE}

The link between the Riemann zeros and random matrix theory
started with the pair correlation of zeros, proposed by 
Montgomery \cite{Montgomery},  and
the observation of Dyson that it is the same as the $2$-point correlation
function predicted by the gaussian unitary ensemble (GUE) for large 
random  matrices \cite{Dyson}.

\begin{figure}
\begin{center}
\includegraphics{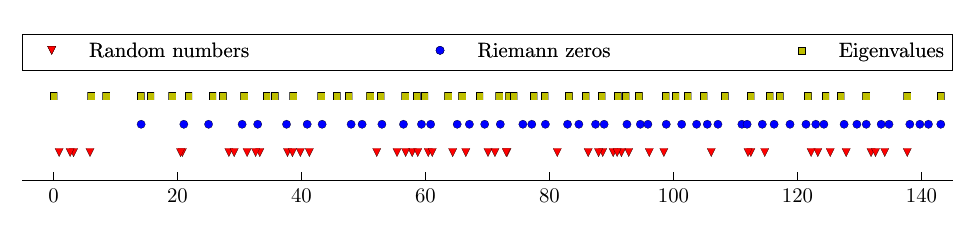}
\vspace{-3em}
\caption{The first $50$ Riemann zeros at $x=1/2$,  
in comparison with 
$50$ random numbers to the left and  
the eigenvalues of a random $50\times 50$ hermitian matrix.}
\label{GUE} 
\end{center}
\end{figure}

The main result of the GUE random matrix theory is that the eigenvalues of
a random hermitian matrix are not completely random,  for instance there 
are correlations in the spacings of eigenvalues,  commonly referred to as 
level-repulsion.     In Figure~\ref{GUE}  we show 3 collections of points: 
the first 50 Riemann zeros, $50$ random real numbers 
between 0 and the ordinate of the $50$-th zero,  and the eigenvalues
of a $50\times 50$ hermitian matrix where each element of the matrix is 
also random. 
One clearly sees the statistical resemblance of the statistics of 
the Riemann zeros and that of the GUE,  in contrast to the $50$ 
random numbers.

The main purpose of this section is to test whether our 
approximation \eqref{FinalTranscendence} to the zeros 
is accurate enough to reveal this statistics.
Whereas formula \eqref{Lambert} 
is a valid estimate of  the zeros, 
it is not sufficiently accurate to reproduce 
the GUE statistics,  since it does not have the 
oscillatory  $\arg \zeta$ term. 
On the other hand, the solutions to equation \eqref{FinalTranscendence} 
are accurate enough,
which again  indicates the importance of the $\arg \zeta$ term.

Montgomery's pair correlation conjecture can be stated as follows:
\beq
\label{montgomery}
\dfrac{1}{N(T)}
\sum_{\substack{
0\le y,y'\le T \\[0.4em]
\alpha < d(y,y') \le \beta}} 
\hspace{-1em} 1
\ \sim \ \int_{\alpha}^{\beta}du
\(1 - \dfrac{\sin^2\(\pi u\)}{\pi^2 u^2}\)
\eeq
where 
$d(y,y') = \tfrac{1}{2\pi}\log\(\tfrac{T}{2\pi}\)\(y-y'\)$,
$ 0 < \alpha<\beta$, $N(T)\sim \tfrac{T}{2\pi}\log\(\tfrac{T}{2\pi}\)$ 
according to \eqref{riemann_counting}, and the statement is valid in 
the limit $T\to \infty$. The
right hand side of \eqref{montgomery} is the 2-point GUE correlation 
function. The average spacing between consecutive zeros is given by 
$\tfrac{T}{N} \sim 2\pi/\log\(\tfrac{T}{2\pi}\)\to 0$  
as $T\to \infty$. This can also be seen from \eqref{Lambert} for very
large $n$, i.e. $\tilde{y}_{n+1}-\tilde{y}_n \to 0$ as $n\to\infty$. 
Thus $d(y,y')$ is a normalized distance. 

While \eqref{montgomery} can be applied if we start from the first
zero on the critical line, it is unable to provide a test if we are centered
around a given high zero on the line. To deal with such a situation,
Odlyzko \cite{Odlyzko2} proposed a stronger version of Montgomery's 
conjecture, by taking into account the large density of zeros 
higher on the line. This is done by replacing
$d(y,y')$ in \eqref{montgomery} by a sum of normalized distances over
consecutive zeros in the form
\beq
\delta_n = \dfrac{1}{2\pi}\log\(\dfrac{y_n}{2\pi}\)\(y_{n+1}-y_n\).
\eeq
Thus \eqref{montgomery} is replaced by
\beq
\label{odlyzko_pair}
\dfrac{1}{\(N-M\)\(\beta-\alpha\)} \hspace{-1em} \sum_{\substack{
M \le m,n \le N \\[0.4em]
\alpha < \sum_{k=1}^{n}\delta_{m+k} \le \beta
}} \hspace{-1.5em} 1
\ \approx \ \dfrac{1}{\beta-\alpha}\int_{\alpha}^{\beta}du
\(1 - \dfrac{\sin^2\(\pi u\)}{\pi^2 u^2}\),
\eeq
where $M$ is the label of a given zero on the line and
$N > M$. In this sum it is assumed that $n > m$ also, and we
included the correct normalization on both sides. The conjecture
\eqref{odlyzko_pair} is already  well supported by extensive 
numerical analysis \cite{Odlyzko2,Gourdon}. 

\begin{figure}
\begin{center}
\begin{minipage}{.49\textwidth}
\includegraphics[width=1\linewidth,trim={15 10 15 10}]{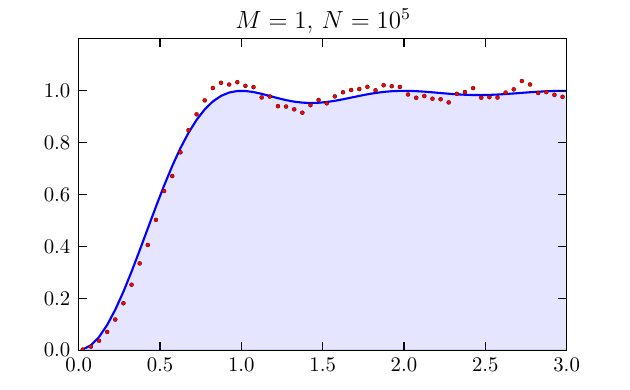}\\[-1em]
\end{minipage}
\begin{minipage}{.49\textwidth}
\includegraphics[width=1\linewidth,trim={15 10 15 10}]{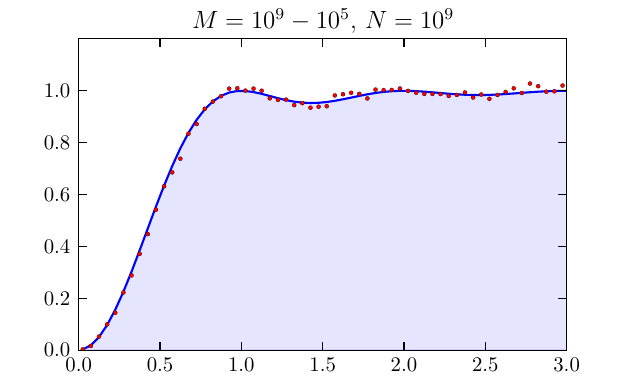}\\[-1em]
\end{minipage}
\end{center}
\caption{
The solid line represents the RHS of \eqref{odlyzko_pair}
and the dots its LHS, computed from equation \eqref{FinalTranscendence}.
The parameters are $\beta = \alpha + 0.05$,
$\alpha=(0,\, 0.05,\,\dotsc,\,3)$ and the $x$-axis is given by
$x=\tfrac{1}{2}\(\alpha+\beta\)$. \emph{ Left:} we use the first $10^5$ zeros.
\emph{ Right:} the same parameters but using zeros in the middle of the
critical line; $M=10^9-10^5$ and $N=10^9$.
}
\label{fig:gue}
\end{figure}

Odlyzko's conjecture \eqref{odlyzko_pair} is a very strong constraint on
the statistics of the zeros. Thus we submit the numerical solutions of 
equation \eqref{FinalTranscendence}, as discussed
in the previous section, to this test. In 
Figure~\ref{fig:gue} (left) we can see the result for $M=1$ and $N=10^{5}$, 
with $\alpha$ ranging from $0\dotsc 3$ in steps of $s=0.05$, and
$\beta=\alpha+s$ for each value of $\alpha$, i.e. 
$\alpha = (0.00, \, 0.05, \, 0.10,  \dotsc,\,3.00)$ and 
$\beta = (0.05,\, 0.10, \, \dotsc, \, 3.05)$. 
We compute the left hand side of \eqref{odlyzko_pair} for each 
pair $(\alpha, \beta)$ and plot
the result against $x = \tfrac{1}{2}\(\alpha + \beta\)$.
In Figure~\ref{fig:gue} (right) we do the same thing but with
$M=10^9-10^5$ and $N=10^9$.

Clearly, the numerical solutions of \eqref{FinalTranscendence} reproduce
the correct statistics. In fact, Figure~\ref{fig:gue} (left) is identical
to the one in \cite{Odlyzko2}.     In Table \ref{high_values} we provide 
the solutions to \eqref{FinalTranscendence}  at the end of the ranges 
considered,
so that the reader may compare with \cite{Odlyzko,Oliveira}.  

\begin{table}
\def\arraystretch{1.2}
\centering
\begin{minipage}{.49\textwidth}
\centering
\begin{tabular}{@{}ll@{}}
\toprule[1pt]
$n$ & $y_n$ \\ 
\midrule[1pt] 
$10^5-5$ & $74917.719415828$ \\
$10^5-4$ & $74918.370580227$ \\
$10^5-3$ & $74918.691433454$ \\
$10^5-2$ & $74919.075161121$ \\
$10^5-1$ & $74920.259793259$ \\
$10^5$   & $74920.827498994$ \\
\bottomrule[1pt]
\end{tabular}
\end{minipage}
\begin{minipage}{.49\textwidth}
\centering
\begin{tabular}{@{}ll@{}}
\toprule[1pt]
$n$ & $y_n$ \\ 
\midrule[1pt] 
$10^9-5$ & $371870202.244870467$ \\
$10^9-4$ & $371870202.673284457$ \\
$10^9-3$ & $371870203.177729799$ \\
$10^9-2$ & $371870203.274345928$ \\
$10^9-1$ & $371870203.802552324$ \\
$10^9$   & $371870203.837028053$ \\
\bottomrule[1pt]
\end{tabular}
\end{minipage}
\caption{Last numerical solutions to \eqref{FinalTranscendence} around
$n=10^5$ and $n=10^9$.}
\label{high_values}
\end{table}

\section{Prime number counting function revisited}  
\label{sec:prime}

\begin{figure}[b]
\begin{center}
\begin{minipage}{.5\textwidth}
\includegraphics[width=1\linewidth,trim={15 10 15 10}]{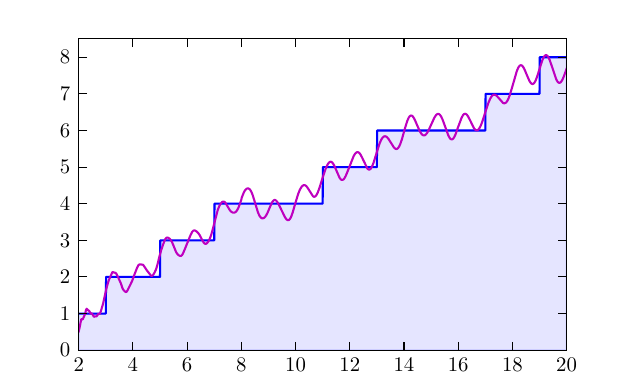}\\[-1em]
\end{minipage}%
\begin{minipage}{.5\textwidth}
\includegraphics[width=1\linewidth,trim={15 10 15 10}]{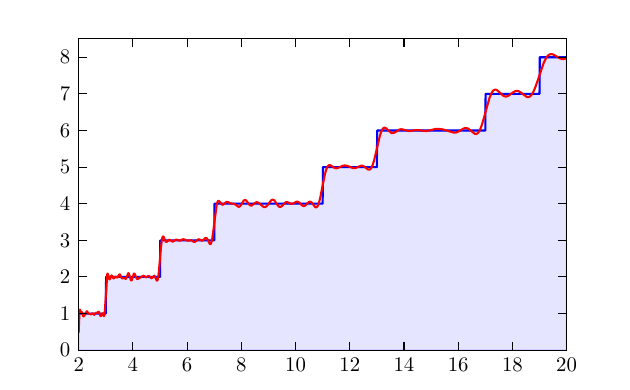}\\[-1em]
\end{minipage}
\caption{The prime number counting function $\pi (x)$ with the 
first $50$ Riemann
zeros. \emph{ Left:} zeros approximated by the Lambert formula \eqref{Lambert}.
\emph{ Right:} zeros obtained from  numerical solutions to the 
equation \eqref{FinalTranscendence}.}
\label{fig:prime}
\end{center}
\end{figure}

In this section we explore whether our approximations to the Riemann zeros
are accurate enough to reconstruct the prime number counting function.     
In Figure~\ref{fig:prime} (left)  we plot $\pi(x)$  from equations \eqref{mob1}
and \eqref{Jzeros}, computed with the first $50$ zeros in the 
approximation $\rho_n = \tfrac{1}{2} + i \tilde{y}_n$ given by 
\eqref{Lambert}. Figure~\ref{fig:prime} (right)  shows the same plot with zeros 
obtained from the numerical solutions of equation \eqref{FinalTranscendence}.
Although with the Lambert  approximation $\tilde{y}_n$ the curve is trying 
to follow the steps in
$\pi (x)$,  once again, one clearly sees the importance 
of the $\arg  \zeta$ term.

\section{Solutions to the exact equation}
\label{sec:numerical_exact}

In the previous sections we have computed numerical solutions of
\eqref{FinalTranscendence} showing that, actually, this first
order approximation to \eqref{exact_eq2} is very good and already 
captures  some  interesting properties of the Riemann zeros, 
like the GUE statistics and ability
to reproduce the prime number counting formula. Nevertheless, by simply 
solving \eqref{exact_eq2} it is possible to obtain values for the zeros
as accurately  as desirable. 
The numerical procedure is performed  as follows:
\begin{enumerate}
\item \label{step1} We apply a root finder method on \eqref{exact_eq2} 
looking for the solution in a region centered around the 
number $\tilde{y}_n$ provided by \eqref{Lambert}, with a not so 
small $\delta$, for instance $\delta \sim 10^{-5}$.
\item \label{step2} We solve \eqref{exact_eq2} again but now centered around 
the solution obtained in step \ref{step1} above, and we decrease $\delta$, for 
instance $\delta \sim 10^{-8}$.
\item We repeat the procedure in step \ref{step2} above, 
decreasing $\delta$ again.
\item Through successive iterations, and decreasing $\delta$ 
each time, it is possible to obtain solutions as accurate as desirable.
In carrying this out,  
it is important to not allow $\delta$ to be exactly zero.   
\end{enumerate}
An actual implementation of the above procedure in Mathematica is shown
in Appendix~\ref{sec:mathematica}.   
The first few zeros computed in this way are shown 
in Table~\ref{lower_precise}.
Through successive iterations it is possible achieve even much higher
accuracy than  these values.

\begin{table}
\def\arraystretch{1.2}
\centering
\begin{tabular}{@{}ll@{}}
\toprule[1pt]
$n$ & $y_n$ \\ 
\midrule[1pt] 
$1$ & $14.1347251417346937904572519835624702707842571156992431756855$ \\
$2$ & $21.0220396387715549926284795938969027773343405249027817546295$ \\
$3$ & $25.0108575801456887632137909925628218186595496725579966724965$ \\
$4$ & $30.4248761258595132103118975305840913201815600237154401809621$ \\
$5$ & $32.9350615877391896906623689640749034888127156035170390092800$ \\
\bottomrule[1pt]
\end{tabular}
\caption{The first few numerical solutions to \eqref{exact_eq2},
accurate to $60$ digits ($58$ decimals).}
\label{lower_precise}
\end{table}

It is known that the first zero where  Gram's law fails is for $n=126$.
Applying the same method, like for any other $n$, the solution of 
\eqref{exact_eq2} starting with the approximation \eqref{Lambert} does
not present any difficulty. We  easily found the following number:
\begin{flalign*}
y_{126} = 279.229250927745189228409880451955359283492637405561293594727
\end{flalign*}
Just to illustrate, and to convince the reader that  the solutions 
of \eqref{exact_eq2} can be made arbitrarily precise, we compute the 
zero $n=1000$ accurate up to $500$ decimal places, also using the same 
simple approach\footnote{Computing this number to $500$ digit accuracy  
took a few minutes on a standard personal laptop 
computer. It only takes a few seconds to obtain 100 digit accuracy.}:
\begin{flalign*}
y_{1000} = 
1419.&42248094599568646598903807991681923210060106416601630469081468460\\
     &86764175930104179113432911792099874809842322605601187413974479526\\
     &50637067250834288983151845447688252593115944239425195484687708163\\
     &94625633238145779152841855934315118793290577642799801273605240944\\
     &61173370418189624947474596756904798398768401428049735900173547413\\
     &19116293486589463954542313208105699019807193917543029984881490193\\
     &19367182312642042727635891148784832999646735616085843651542517182\\
     &417956641495352443292193649483857772253460088
\end{flalign*}
Substituting precise Riemann zeros into \eqref{exact_eq2} one
can check that the equation is identically satisfied.
These results corroborate that \eqref{exact_eq2} is an
exact equation for the Riemann zeros.

\section{Numerical analysis:  Dirichlet \texorpdfstring{$L$}{L}-functions}
\label{sec:numerical_lfunc}

We perform exactly the same numerical procedure as described
in the previous section \ref{sec:numerical_exact}, 
but now with equation \eqref{exact} and \eqref{approx} for 
Dirichlet $L$-functions.  

We will illustrate our formulas with the primitive 
characters $\chi_{7,2}$ and $\chi_{7,3}$,   
since they possess the full  generality of $a=0$ and $a=1$ and 
complex components.  
There are actually $\varphi(7)=6$ distinct characters mod $7$.   

\paragraph*{\bf Example \boldmath{$\chi_{7,2}$}.}
Consider $k=7$ and $j=2$, i.e. we are computing
the Dirichlet character $\chi_{7,2}(n)$. For this case $a=1$.
Then we have the following components:
\beq\label{char72}
\begin{tabular}{@{}c|ccccccc@{}}
$n$             & $1$ & $2$ & $3$ & $4$ & $5$ & $6$ & $7$ \\
\midrule[0.3pt]
$\chi_{7,2}(n)$ &  
$1$ & $e^{ 2\pi i /3}$ & $e^{\pi i / 3}$ & $e^{-2\pi i / 3}$ &
$e^{-\pi i / 3}$ & $-1$ & $0$ 
\end{tabular}
\eeq
The first few zeros, positive and negative, obtained by solving
\eqref{exact} are shown in Table~\ref{zeros_1} 
(see Appendix~\ref{sec:mathematica}). 
The solutions shown are
easily obtained with $50$ decimal places of accuracy, and agree with
the ones in \cite{Oliveira2}, which were computed up to $20$ decimal places.

\begin{table}
\def\arraystretch{1.2}
\centering
\begin{tabular}{@{}rrr@{}}
\toprule[1pt]
$n$ & $\tilde{y}_n$ & $y_n$ \\
\midrule[1pt]
$10$ &  $25.57$ &  $25.68439458577475868571703403827676455384372032540097$ \\
$9$ &   $23.67$ &  $24.15466453997877089700472248737944003578203821931614$ \\
$8$ &   $21.73$ &  $21.65252506979642618329545373529843196334089625358303$ \\
$7$ &   $19.73$ &  $19.65122423323359536954110529158230382437142654926200$ \\
$6$ &   $17.66$ &  $17.16141654370607042290552256158565828745960439000612$ \\
$5$ &   $15.50$ &  $15.74686940763941532761353888536874657958310887967059$ \\
$4$ &   $13.24$ &  $13.85454287448149778875634224346689375234567535103602$ \\
$3$ &   $10.81$ &   $9.97989590209139315060581291354262017420478655402522$ \\
$2$ &    $8.14$ &   $8.41361099147117759845752355454727442365106861800819$ \\
$1$ &    $4.97$ &   $5.19811619946654558608428407430395403442607551643259$ \\
$0$ &   $-3.44$ &  $-2.50937455292911971967838452268365746558148671924805$ \\
$-1$ &  $-7.04$ &  $-7.48493173971596112913314844807905530366284046079242$ \\
$-2$ &  $-9.85$ &  $-9.89354379409772210349418069925221744973779313289503$ \\
$-3$ & $-12.35$ & $-12.25742488648921665489461478678500208978360618268664$ \\
$-4$ & $-14.67$ & $-14.13507775903777080989456447454654848575048882728616$ \\
$-5$ & $-16.86$ & $-17.71409256153115895322699037454043289926793578042465$ \\
$-6$ & $-18.96$ & $-18.88909760017588073794865307957219593848843485334695$ \\
$-7$ & $-20.99$ & $-20.60481911491253262583427068994945289180639925014034$ \\
$-8$ & $-22.95$ & $-22.66635642792466587252079667063882618974425685038326$ \\
$-9$ & $-24.87$ & $-25.28550752850252321309973718800386160807733038068585$ \\
\bottomrule[1pt]
\end{tabular}
\caption{Numerical solutions of \eqref{exact} starting
with the approximation \eqref{approx}, for the
character \eqref{char72}. The solutions are
accurate to $50$ decimal places and verified to
$\left|L\(\tfrac{1}{2} + i y_n\)\right| \sim 10^{-50}$.}
\label{zeros_1}
\end{table}

\paragraph*{\bf Example \boldmath{$\chi_{7,3}$}.} 
Consider $k=7$ and $j=3$, such that $a=0$.
In this case the components of $\chi_{7,3}(n)$ are the following:
\beq\label{char73}
\begin{tabular}{@{}c|ccccccc@{}}
$n$             & $1$ & $2$ & $3$ & $4$ & $5$ & $6$ & $7$ \\
\midrule[0.3pt]
$\chi_{7,3}(n)$ &  
$1$ & $e^{ -2\pi i /3}$ & $e^{2 \pi i / 3}$ & $e^{2\pi i / 3}$ &
$e^{-2 \pi i / 3}$ & $1$ & $0$ 
\end{tabular}
\eeq
The first few solutions of \eqref{exact} are shown in 
Table~\ref{zeros_2} and are accurate up to $50$ decimal places, and
agree with the ones obtained in \cite{Oliveira2}.

\begin{table}
\def\arraystretch{1.2}
\centering
\begin{tabular}{@{}rrr@{}}
\toprule[1pt]
$n$ & $\tilde{y}_n$ & $y_n$ \\
\midrule[1pt]
$10$ &  $25.55$  &  $26.16994490801983565967242517629313321888238615283992$ \\
$9$ &   $23.65$  &  $23.20367246134665537826174805893362248072979160004334$ \\
$8$ &   $21.71$  &  $21.31464724410425595182027902594093075251557654412326$ \\
$7$ &   $19.71$  &  $20.03055898508203028994206564551578139558919887432101$ \\
$6$ &   $17.64$  &  $17.61605319887654241030080166645399190430725521508443$ \\
$5$ &   $15.48$  &  $15.93744820468795955688957399890407546316342953223035$ \\
$4$ &   $13.21$  &  $12.53254782268627400807230480038783642378927939761728$ \\
$3$ &   $10.79$  &  $10.73611998749339311587424153504894305046993275660967$ \\
$2$ &    $8.11$  &   $8.78555471449907536558015746317619235911936921514074$ \\
$1$ &    $4.93$  &   $4.35640162473628422727957479051551913297149929441224$ \\
$0$ &   $-5.45$  &  $-6.20123004275588129466099054628663166500168462793701$ \\
$-1$ &  $-8.53$  &  $-7.92743089809203774838798659746549239024181788857305$ \\
$-2$ & $-11.15$  & $-11.01044486207249042239362741094860371668883190429106$ \\
$-3$ & $-13.55$  & $-13.82986789986136757061236809479729216775842888684529$ \\
$-4$ & $-15.80$  & $-16.01372713415040781987211528577709085306698639304444$ \\
$-5$ & $-17.94$  & $-18.04485754217402476822077016067233558476519398664936$ \\
$-6$ & $-20.00$  & $-19.11388571948958246184820859785760690560580302023623$ \\
$-7$ & $-22.00$  & $-22.75640595577430793123629559665860790727892846161121$ \\
$-8$ & $-23.94$  & $-23.95593843516797851393076448042024914372113079309104$ \\
$-9$ & $-25.83$  & $-25.72310440610835748550521669187512401719774475488087$ \\
\bottomrule[1pt]
\end{tabular}
\caption{Numerical solutions of \eqref{exact} starting
with the approximation \eqref{approx}, for the character \eqref{char73}. 
The solutions are
accurate to $50$ decimal places and verified to
$\left|L\(\tfrac{1}{2} + i y_n\)\right| \sim 10^{-50}$.}
\label{zeros_2}
\end{table}

As stated previously, the solutions to equation \eqref{exact} can be 
calculated to any desired level of accuracy. For instance, continuing with 
the character $\chi_{7,3}$, we can easily compute the following 
number for $n=1000$, accurate to $100$ decimal places, i.e. $104$ digits:
\beq\nonumber
\begin{split}
y_{1000} = 1037.&56371706920654296560046127698168717112749601359549 \\
                &01734503731679747841764715443496546207885576444206
\end{split}
\eeq

We also have been able to solve the equation for high zeros to 
high accuracy, up to the millionth zero, some of which are 
listed in Table~\ref{highzeros},   and were previously unknown.   

\begin{table}
\def\arraystretch{1.2}
\centering
\begin{tabular}{@{}crr@{}}
\toprule[1pt]
$n$ & $\tilde{y}_n$ & $y_n$ \\
\midrule[1pt]
$10^3$ & $1037.61$ & 
  $1037.563717069206542965600461276981687171127496013595490$ \\
$10^4$ & $7787.18$ &      
  $7787.337916840954922060149425635486826208937584171726906$ \\
$10^5$ & $61951.04$ &  
 $61950.779420880674657842482173403370835983852937763461400$ \\
$10^6$ & $512684.78$ & 
$512684.856698029779109684519709321053301710419463624401290$ \\
\bottomrule[1pt]
\end{tabular}
\caption{
Higher zeros for the Dirichlet character \eqref{char73}.  
These solutions to \eqref{exact}  are
accurate to $50$ decimal places.}
\label{highzeros}
\end{table}

\section{Modular \texorpdfstring{$L$}{L}-function based on Ramanujan 
\texorpdfstring{$\tau$}{tau}}
\label{sec:ramanujan}

 In this  section we study an $L$-function based on a level one modular form
 related to the Ramanujan $\tau$ function,  and numerically solve the equations
\eqref{exactmod} and \eqref{yLambertmod}.

\subsection{Definition of the function} 

Here we will consider an example of a modular form of weight $k=12$.
The simplest example is based on the Dedekind $\eta$-function
\beq \label{eta}
\eta (\tau) = q^{1/24} \, \prod_{n=1}^\infty (1-q^n) .
\eeq
Up to a simple factor,  $\eta$ is the inverse of the chiral  partition 
function of the free boson conformal field theory \cite{CFT},  where $\tau$ is
the modular parameter of the torus.   
The modular discriminant 
\beq \label{Delta} 
\Delta (\tau ) = \eta (\tau )^{24}  =  \sum_{n=1}^\infty \, \tau(n) \, q^n
\eeq
is a weight $k=12$ modular form.   
It is closely related to the inverse of the partition function 
of the bosonic string in $26$ dimensions,  where $24$ is the number 
of light-cone degrees of freedom \cite{StringTheo}.     
The Fourier coefficients  $\tau (n)$  correspond to  the 
Ramanujan $\tau$-function, and the first few are 
\beq
\begin{tabular}{@{}c|cccccccc@{}}
$n$             & $1$ & $2$ & $3$ & $4$ & $5$ & $6$ & $7$ & $8$ \\
\midrule[0.3pt]
$\tau(n)$ &  $1$ & $-24$ & $252$ & $-1472$ & $4830$ & $-6048$ & 
$-16744$ & $84480$  
\end{tabular}
\eeq

\begin{figure}
\centering
\includegraphics[width=0.6\linewidth]{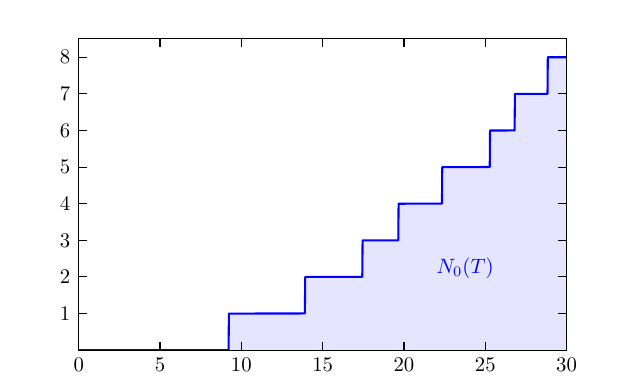}
\caption{Exact counting formula \eqref{countRam} based on the
Ramanujan $\tau$-function.}
\label{fig:ram_counting}
\end{figure}

We then define the Dirichlet series
\beq \label{LRam}
L_{\Delta}(z)  =  \sum_{n=1}^\infty  \,  \frac{\tau (n) }{n^z}.
\eeq
Applying \eqref{exactmod}, the zeros are $\rho_n = 6 + i y_n$,  
where  the $y_n$ satisfy the exact equation
\beq \label{exactRam}
\vartheta_{12} (y) + 
\lim_{\delta \to 0^{+}} \arg L_{\Delta}(6 + \delta + i y_n ) = 
\(n-\tfrac{1}{2}\) \pi.
\eeq
The counting function \eqref{NTmod} and its asymptotic approximation are
\begin{align} 
\label{countRam}
N_0(T) &= \inv{\pi} \vartheta_{12} (T) + 
\inv{\pi} \arg L_{\Delta} (6 + i T)  \\
&\approx
\dfrac{T}{\pi} \log \( \dfrac{T}{2 \pi e} \)  + \inv{\pi}\arg L_\Delta(6+iT) +
\dfrac{11}{4}.
\end{align}
A plot of \eqref{countRam} is shown in 
Figure~\ref{fig:ram_counting}, and we can see that it is a perfect staircase
function.

The approximate solution \eqref{yLambertmod} now has  the form
\beq \label{LambertRam} 
\tilde{y}_n = \dfrac{ \(n- \tfrac{13}{4} \) \pi }{W\[(2e)^{-1}
\(n-\tfrac{13}{4}\)\]}
\eeq
for $n=2,3,\dotsc$.
Note that the above equation is valid for $n > 1$, 
since $W(x)$ is not defined for $x < -1/e$.

\begin{table}[t]
\def\arraystretch{1.2}
\centering
\begin{tabular}{@{}rrr@{}}
\toprule[1pt]
$n$ & $\tilde{y}_n$ & $y_n$ \\
\midrule[1pt]
$1$   &          &   $9.22237939992110252224376719274347813552877062243201$ \\
$2$   &  $12.46$ &  $13.90754986139213440644668132877021949175755235351449$ \\
$3$   &  $16.27$ &  $17.44277697823447331355152513712726271870886652427527$ \\
$4$   &  $19.30$ &  $19.65651314195496100012728175632130280161555091200324$ \\
$5$   &  $21.94$ &  $22.33610363720986727568267445923624619245504695246527$ \\
$6$   &  $24.35$ &  $25.27463654811236535674532419313346311859592673122941$ \\
$7$   &  $26.60$ &  $26.80439115835040303257574923358456474715296800497933$ \\
$8$   &  $28.72$ &  $28.83168262418687544502196191298438972569093668609124$ \\
$9$   &  $30.74$ &  $31.17820949836025906449218889077405585464551198966267$ \\
$10$  &  $32.68$ &  $32.77487538223120744183045567331198999909916163721260$ \\
$100$ & $143.03$ & $143.08355526347845507373979776964664120256210342087127$ \\
$200$ & $235.55$ & $235.74710143999213667703807130733621035921210614210694$ \\
$300$ & $318.61$ & $318.36169446742310747533323741641236307865855919162340$ \\
\bottomrule[1pt]
\end{tabular}
\caption{Non-trivial zeros of the modular $L$-function based
on the Ramanujan $\tau$-function, obtained from \eqref{exactRam} starting
with the approximation \eqref{LambertRam}.  
These solutions are accurate to $50$ decimal places.}
\label{zerosRam}
\end{table}

\subsection{Numerical analysis}

We follow  the same procedure, previously discussed in
section \ref{sec:numerical_exact} 
and  implemented in Appendix \ref{sec:mathematica}, to solve
the equation \eqref{exactRam} starting with the approximation
provided by \eqref{LambertRam}. Some of these solutions are shown
in Table~\ref{zerosRam} and are accurate to $50$ decimal places.

\section{Counterexample of Davenport-Heilbronn}  
 \label{sec:Davenport}
 
The \emph{Davenport-Heilbronn} function has almost all the same 
properties of $\zeta$, such as a functional equation, except that it has
no Euler product formula. It is well known that such a function
has zeros in the region $\Re\(z\) > 1$ and zeros
in the critical strip  $0 \le \Re(z) \le 1$, of which infinitely many of them
lie on \emph{the critical line} $\Re(z) = 1/2$,  but it also \emph{has
zeros off of the critical line}, thus violating the RH.
Whether this is due to the absence of an  Euler product is not understood,  although this is likely to be the reason.
It is very  interesting to apply  the formalism of the previous 
sections to this function
and to understand clearly why the RH fails here.  Such an exercise 
sharpens our understanding of the RH.
For a detailed study of this function and numerical computation
of the zeros see \cite{BombieriGhosh}.

This function is somewhat contrived since it is linear combination of 
Dirichlet $L$-functions,   engineered  to satisfy a functional equation.   
It is  defined by
\beq\label{davenport_def}
\CD(z) \equiv \dfrac{(1-i \kappa)}{2} L\(z,\chi_{5,2}\) + 
\dfrac{(1+i \kappa)}{2} L\(z,\chi_{5,2}^{*}\)
\eeq
with 
\beq
\kappa =\dfrac{\sqrt{ 10-2\sqrt{5}}-2 }{\sqrt{5} - 1},
\eeq
and the Dirichlet
character is given by
\beq\label{char52}
\begin{tabular}{@{}c|ccccccc@{}}
$n$             & $1$ & $2$ & $3$ & $4$ & $5$ \\
\midrule[0.3pt]
$\chi_{5,2}(n)$ &  
$1$ & $i$ & $-i$ & $-1$ & $0$
\end{tabular}
\eeq
and $\chi_{5,2}(-1) = -1$ thus $a=1$. This function satisfies the functional
equation
\beq\label{davenport_func_eq}
\xi(z) = \xi(1-z), \qquad 
\xi(z) \equiv \(\dfrac{\pi}{5}\)^{-z/2}\, \Gamma\(\dfrac{1+z}{2}\) \CD(z).
\eeq
$\CD (z)$ has no Euler product because it is a linear combination 
of functions that do.

Now we repeat the analysis we presented for the Riemann $\zeta$-function.
We have 
\beq
\xi(z) = A e^{i\theta}, \qquad
\xi(1-z) = A'e^{-i\theta'}, 
\eeq 
with $A(x,y)=A'(x,y)$ and $\theta'(x,y)=\theta(1-x,y)$, where
\beq\label{davenport_theta}
\theta(x,y) = 
\Im \[ \log \Gamma\(\dfrac{1+x+i y}{2}\) \] - 
\dfrac{y}{2}\log\(\dfrac{\pi}{5}\) + \arg \CD(x+iy).
\eeq
As for previous $L$-functions,  zeros should be characterized by our equation 
\eqref{Bzero}, namely
\beq \label{Bzero2}
\theta + \theta ' = (2n+1) \pi.
\eeq
The zeros on the critical line  correspond to 
\beq
\theta = \theta'  \quad \Rightarrow \quad \theta = \(n+ \tfrac{1}{2}\) \pi.
\eeq
As before, adopting the convention that
the first positive zero is labeled by $n=1$, we shift $n\to n-1$, then
the equation yielding zeros on the critical line is given by
\beq\label{davenport_transcendental}
\Xi(y_n) + \lim_{\delta\to0^+} \arg \CD\(\tfrac{1}{2}+\delta+iy_n\) = 
\(n-\tfrac{1}{2}\)\pi
\eeq
where
\beq\label{davenport_rs}
\Xi(y) \equiv \Im \[ \log\Gamma\(\dfrac{3}{4}+i\dfrac{y}{2}\) \] - 
\dfrac{y}{2}\log\(\dfrac{\pi}{5}\).
\eeq
Expanding $\Gamma(z)$ through Stirling's formula and
neglecting the $\arg \CD$ term in \eqref{davenport_transcendental} 
it is possible
to obtain an explicit approximate solution given by
\beq\label{davenport_lambert}
\tilde{y}_n = \dfrac{2\pi\(n-\tfrac{5}{8}\)}{W\[ 5e^{-1}(n-\tfrac{5}{8}) \]},
\eeq
where $n=1,2,\dotsc$ and $W$ denotes the principal 
branch $W_0$ of the Lambert function \eqref{W_def_eq}.

Now we can numerically solve \eqref{davenport_transcendental} starting with the 
approximation given by \eqref{davenport_lambert}. 
The first few solutions are shown in Table~\ref{zeros_davenport} (left).

\begin{table}[b]
\footnotesize
\def\arraystretch{1.2}
\begin{minipage}{.39\textwidth}
\begin{tabular}{@{}rrr@{}}
\toprule[1pt]
$n$ & $\tilde{y}_n$ & $y_n$ \\
\midrule[1pt]
$1$   &  $5.32$  &   $5.094159844584467267$ \\
$2$   &  $8.96$  &   $8.939914408100472858$ \\
$3$   &  $11.93$ &   $12.133545425790163309$ \\
$4$   &  $14.60$ &   $14.404003112292645158$ \\
$5$   &  $17.08$ &   $17.130239400567288918$ \\
$6$   &  $19.43$ &   $19.308800174241700381$ \\
$7$   &  $21.68$ &   $22.159707765035018919$ \\
$8$   &  $23.85$ &   $23.345370112090190151$ \\
$9$   &  $25.95$ &   $26.094967346227912542$ \\
$10$  &  $28.00$ &   $27.923798821611878096$ \\
\bottomrule[1pt]
\end{tabular}
\end{minipage}
\begin{minipage}{.6\textwidth}
\begin{tabular}{@{}ccc@{}}
\toprule[1pt]
$\rho$ & $\tfrac{1}{\pi}\theta$ & $\tfrac{1}{\pi}(\theta+\theta')$ \\
\midrule[1pt]
 $0.8085171825 + i\,
 85.6993484854$         &   $44.092$   &   $89$    \\
 $0.6508300806 + i\,
114.1633427308$         &   $64.026$   &   $127$   \\
 $0.5743560504 + i\,
166.4793059132$         &   $103.023$  &   $207$   \\
 $0.7242576946 + i\,
176.7024612429$         &   $111.075$  &   $223$   \\
 $0.8695305796 + i\,
240.4046723514$         &   $163.055$  &   $325$   \\
 $0.8195495921 + i\,
320.8764896688$         &   $232.106$  &   $465$   \\
 $0.7682231236 + i\,
331.0502594079$         &   $241.098$  &   $483$   \\
 $0.6285081083 + i\,
366.6409075762$         &   $273.027$  &   $545$   \\
 $0.8158736778 + i\,
411.7967375490$         &   $314.133$  &   $629$   \\
 $0.7088822242 + i\,
440.4845107397$         &   $341.017$  &   $681$   \\
\bottomrule[1pt]
\end{tabular}
\end{minipage}
\caption{\emph{ Left:} first few zeros of \eqref{davenport_def} on the critical 
line, $\rho_n = \smallhalf+iy_n$, computed from 
\eqref{davenport_transcendental} 
starting with the approximation \eqref{davenport_lambert}.
\emph{ Right: } We can see that equation \eqref{Bzero2} is indeed verified 
for the 
first few zeros off of the critical line. Note that $\tfrac{1}{\pi}\theta$
can be any real number, while the combination 
$\tfrac{1}{\pi}\(\theta+\theta'\)$ always gives an odd integer at a 
non-trivial zero.
}
\label{zeros_davenport}
\end{table}

\begin{figure}
\begin{center}
\begin{minipage}{.5\textwidth}
\includegraphics[width=.8\linewidth]{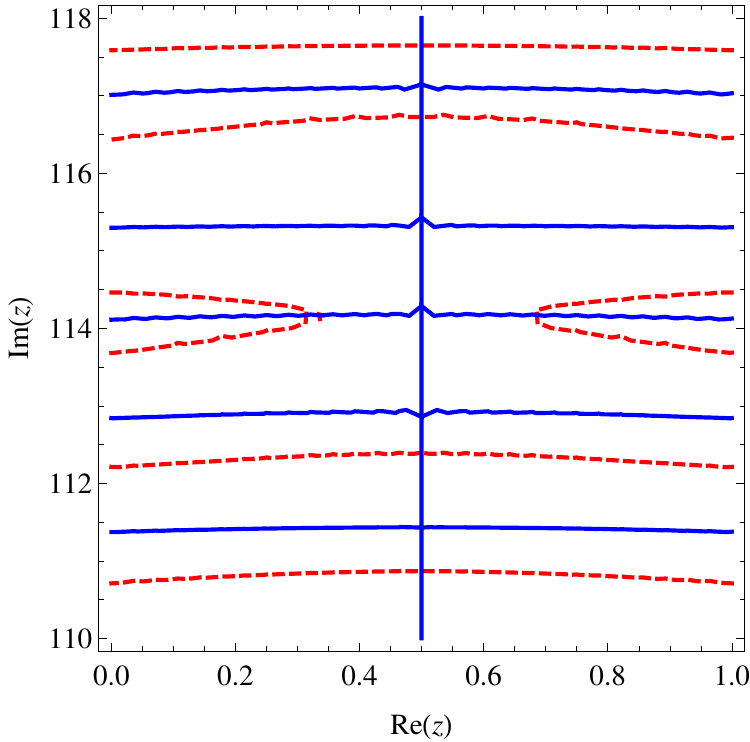}
\end{minipage}%
\begin{minipage}{.5\textwidth}
\includegraphics[width=.8\linewidth]{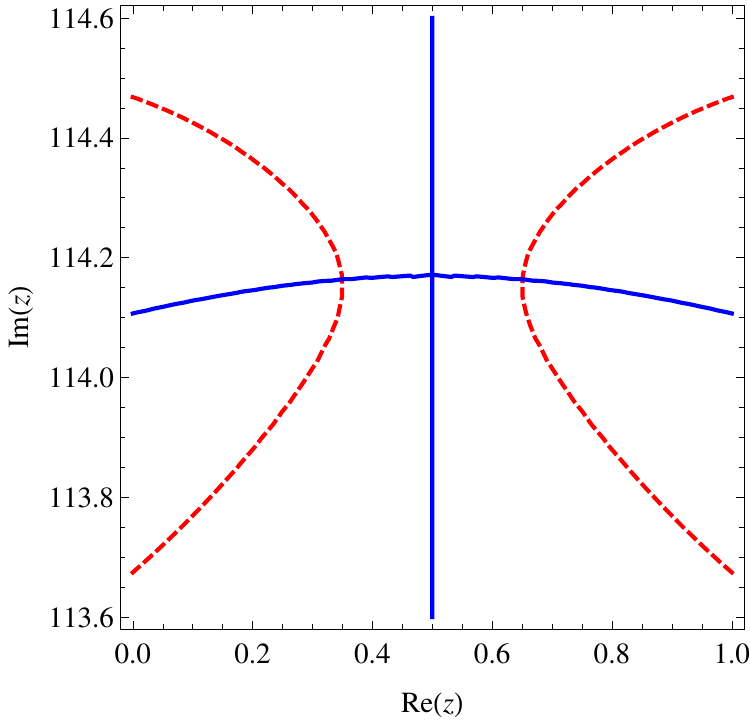}
\end{minipage}
\caption{Contour plot of $\xi=u + iv$, 
equation \eqref{davenport_func_eq}, for $x=0\dotsc1$
and $y=110\dotsc118$. Note the zeros off-line for $y\approx 114.1633$.
The solid (blue) lines correspond to curves $v=0$ and the dashed (red) 
lines to $u=0$.}
\label{dav_contour}
\end{center}
\end{figure}

As for the trivial zeros of $\zeta$ we expect that 
zeros off of the line also satisfy \eqref{Bzero2}.   
We indeed verified this. It is more difficult to find 
these zeros since they are at scattered values of $x$,
but it is in fact feasible.  In Table~\ref{zeros_davenport} (right)
we show some of the lower zeros and the values of 
$\tfrac{1}{\pi}(\theta+\theta')$, which are odd integers.
In Figure~\ref{dav_contour} we show the contour 
lines of \eqref{davenport_func_eq}, i.e. $\xi(z) = u(x,y) + i v(x,y)$,
and we consider the lines $u=0$ and $v=0$. On the critical line
$x=1/2$ we have a $v=0$ contour everywhere, and we approach the
zero on a $u=0$ contour through the $\delta$ limit, as shown
in equation \eqref{davenport_transcendental}. In this case the $\delta$
smooths out the discontinuity. However, note how
the curves $u=0$ are very different in nature for zeros off of the critical
line. 
In this case the $\delta$ limit does not smooth out the function, 
since the path to approach the zero is more involved. Nevertheless, equation
\eqref{Bzero2} is still satisfied for these zeros off of the critical line.

\begin{figure}
\begin{center}
\begin{minipage}{.5\textwidth}
\includegraphics[width=.9\linewidth]{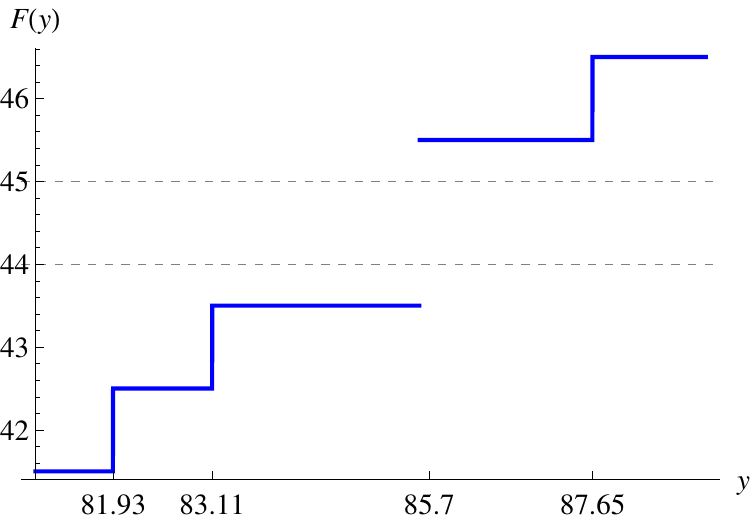}
\end{minipage}%
\begin{minipage}{.5\textwidth}
\includegraphics[width=.9\linewidth]{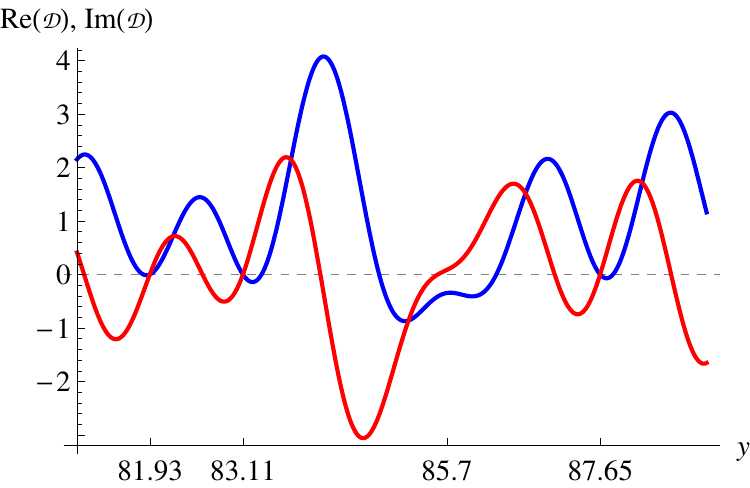}
\end{minipage}
\caption{
\emph{ Left: } 
$F(y) = \tfrac{1}{\pi}\(\Xi(y) + \arg \CD\(\tfrac{1}{2}+\delta+iy\)\) + 
\tfrac{1}{2}$ versus $y$. Note the
discontinuity in the graph for  $y\approx85.6993$, corresponding 
to $n=44$ and $n=45$, where the equation \eqref{davenport_transcendental} 
has no solution. 
\emph{ Right: } the blue line is
$\Re\(\CD\(\tfrac{1}{2}+iy\)\)$ and the red line is 
$\Im\(\CD\(\tfrac{1}{2}+iy\)\)$. Equation 
\eqref{davenport_transcendental} is not defined since 
$\arg\CD\(\tfrac{1}{2}+iy\)$ changes branch. Note that
$\Re(\CD) < 0$ and $\Im(\CD) = 0$ at $y\approx 85.6993$, where there
are two zeros off the line at this height.}
\label{davenport_43}
\end{center}
\end{figure}

If there is a unique solution to 
\eqref{davenport_transcendental}  for every $n$, 
then as for the $\zeta$-function we can determine $N_0 (T)$ which 
counts zeros on the line.
However, the equation \eqref{davenport_transcendental} \emph{is not defined}
for every $n$. For instance, for $n=44$ and $n=45$ this equation has
no solution, as illustrated in Figure~\ref{davenport_43}. The same thing 
happens
again for $n=64$ and $n=65$, for $n=103$ and $n=104$, and so on.

\begin{figure}
\begin{center}
\begin{minipage}{0.49\textwidth}
\includegraphics[width=1\linewidth]{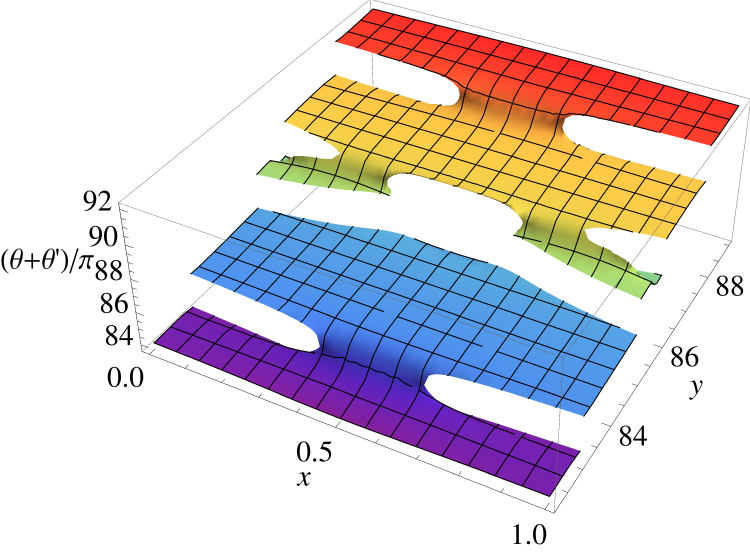}
\end{minipage}
\begin{minipage}{0.49\textwidth}
\includegraphics[width=.9\linewidth]{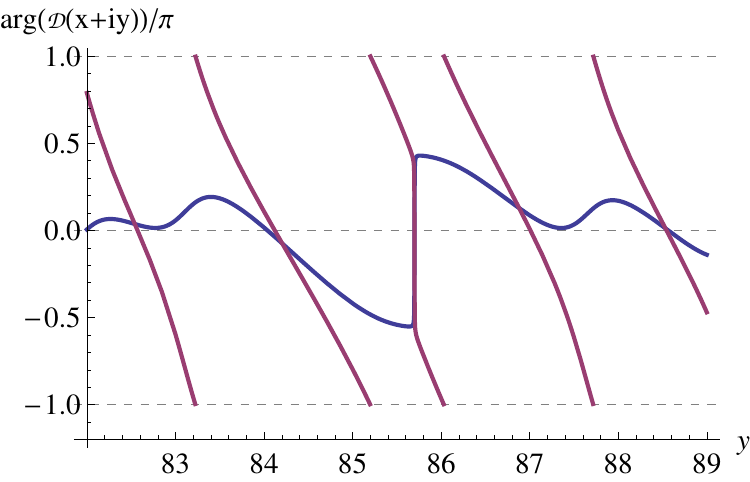}
\end{minipage}
\caption{
\emph{ Left:} The zeros off-line do not allow the function to be smoothed
on the critical line.
\emph{ Right:} the blue line is $\tfrac{1}{\pi}\arg\CD(x+iy)/\pi $ versus $y$, 
in the vicinity of the first zero off the line, 
where $x\approx0.8085$.
The purple line is $\tfrac{1}{\pi}\arg\CD(1-x+iy)$.
}
\label{theta3d}
\end{center}
\end{figure}

If $N(T)$ counts zeros on the entire strip,  then clearly $N_0 (T) \neq N(T)$.
For these values of $y$, corresponding to zeros off the line,  
$\lim_{\delta\to0^+}\arg\CD\(\tfrac{1}{2}+\delta+iy\)$ is not defined, i.e.
the $\delta$ limit \emph{does not smooth out} the function since there is a
severe change of branch. This is why the 
equation \eqref{davenport_transcendental} is not defined in the vicinity
of such $y$'s.
In Figure~\ref{theta3d} (left) we show a plot of 
$\tfrac{1}{\pi}\[\theta(x+\delta,y) + \theta'(x+\delta,y)\]$ 
for $x=0\dotsc1$ and
$y=81\dotsc90$. Note that there are zeros off the line at $y\approx85.6993$.
We included a $\delta\sim 10^{-1}$. Note how the function can be
smoothed on the critical line $x=1/2$ if there are only zeros on the line. 
However, when there are zeros off of the critical line the function 
cannot be made continuous (note  the ``big hole'' around the critical line).

Finally, in Figure~\ref{theta3d} (right) we plot the analog of 
Figure~\ref{fig_ArgPlusMinus} (right).
In the vicinity of the first zero off  the critical line,   
both $\theta$ and $\theta'$ are well-defined
and there is a solution to \eqref{Bzero2}.

For Riemann $\zeta$, we provided arguments, though not a proof,   that there is a unique solution 
to the analogous equation \eqref{exact_eq2}, implying that $N_0(T) =N(T)$.
If $\arg\zeta\(\tfrac{1}{2}+iy\)$ does not change branch, or if it 
changes ``very little'', the $\delta$ limit smooths out the function, making the equation well defined. 
On the other hand, if there is a severe change of branch 
it is impossible to make this function continuous, which is the case
here for the function \eqref{davenport_def} (see Figures~\ref{davenport_43}
and \ref{theta3d}). In the case of the Davenport-Heilbronn function,
Karatsuba \cite{KaratsubaDavenport} showed that at 
least $T (\log T)^{1/2-\epsilon}$ zeros
lie on the critical line.

\begin{figure}
\begin{center}
\begin{minipage}{.49\textwidth}
\includegraphics[width=.9\linewidth]{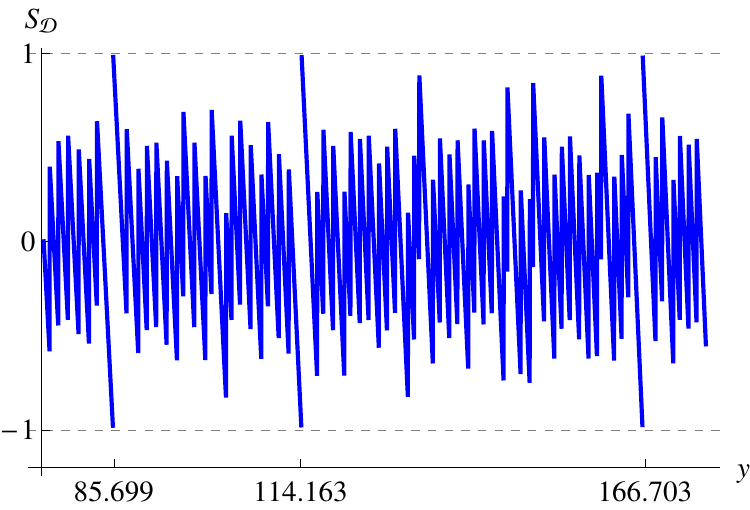}
\end{minipage}%
\begin{minipage}{.49\textwidth}
\includegraphics[width=.9\linewidth]{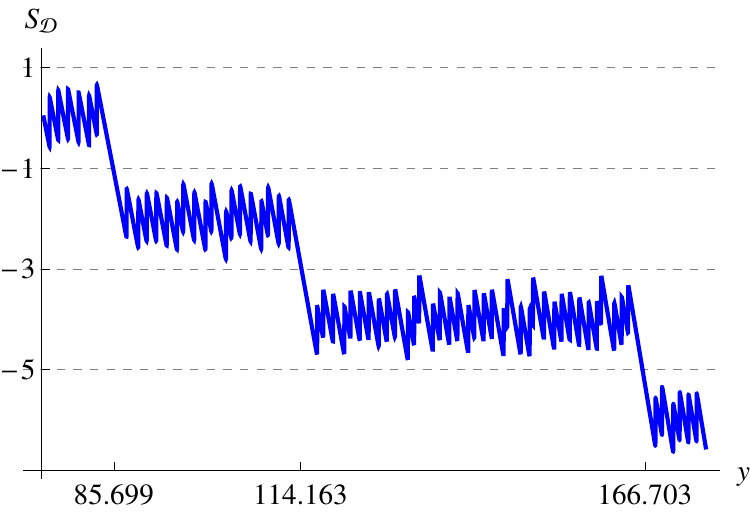}
\end{minipage}
\caption{\emph{ Left:} $S_{\CD}$ in 
the principal branch. 
\emph{ Right: } $S_{\CD}$ explicitly showing the changes in branch.
 $\langle S_{\CD} \rangle \ne 0$ in this  range.}
\label{sy_davenport}
\end{center}
\end{figure}

The important difference with $\zeta$ is that 
\beq
S_\CD (y) =  \lim_{\delta \to 0^+}  \inv{\pi}  
\arg \CD \(\smallhalf + \delta + y\)
\eeq
has very different properties in comparison to 
$S(y)=\lim_{\delta\to0^+}\tfrac{1}{\pi}\arg\zeta\(\tfrac{1}{2}+\delta+iy\)$
or the argument of the other $L$-functions we have considered.      
The properties conjectured for $S(y)$ in section \ref{Soft} 
should not hold. The repeated changes in branch suggests that $S_\CD (y)$ has
the behavior sketched in Figure~\ref{sy_davenport},  however we have not proven this.  

\section{Saddle point approximation}  
\label{sec:saddle} 

The original aim of the work in this section was to obtain an 
approximate analytic 
expression for $\arg \zeta (\smallhalf + iy )$  because of its importance 
to the whole theory of the zeros.
For  the smooth part, $\arg \Gamma$,  one has such an expression due to 
the Stirling's approximation,   
and we will present something similar for $\zeta$.

The saddle-point method, also known as steepest descent method, approximates
integrals of the type
\beq\label{saddle_integral}
I(\lambda) = \int_{\CC} e^{\lambda f(u)} du
\eeq
for large $\lambda$. The saddle-points are  solutions to
$f'(u) = 0$, provided $f''(u) \ne 0$,  where $f'$ and $f''$ are the 
first and second derivatives of $f$. 
For simplicity let us consider only 
one saddle point denoted by $s$. The contour $\CC$ must be deformed
in such a way as to pass through $u=s$ with $\Re(u)$ having the steepest
descent (hence the name).   If $f(u)$ is an \emph{analytic} function, 
using  $f(u) \approx  f(s) + \tfrac{1} {2}  f''(s) (u-s)^2 $ one is 
left with a gaussian integral and obtains 
\beq\label{saddle_approx}
I(\lambda) \approx \sqrt{\dfrac{2\pi}{-f''\(s\)}} \ e^{\lambda f(s)}.
\eeq

Let us consider a very illustrative example,  the result of which we will 
need later. The $\Gamma$-function
has the integral representation
\beq\label{gamma_int_def}
\Gamma(z) = \int_{0}^{\infty} u^{z-1} e^{-u} du, \qquad \Re(z) > 0.
\eeq
Writing the integrand as $e^{f(u)}$ where $f(u) = -u + (z-1)\log u$, then
$f'(u) = 0$ implies $s = z-1$. Then \eqref{saddle_approx} yields
\beq\label{stirling2}
\Gamma(z) \approx \sqrt{2\pi} (z-1)^{z-1/2}\,e^{-(z-1)}.
\eeq
For $n$  a positive integer,  
$\Gamma(n+1) = n!   \approx \sqrt{2\pi} \, n^{n+1/2}\, e^{-n}$, 
which is the very useful Stirling's approximation.   
Since the above approximation is an analytic function of $z$,  
it is valid in the whole complex plane through analytic continuation.    
It is very useful for instance to determine the asymptotic expansion of
the Riemann-Siegel $\vartheta$ function,
\beq
\label{RS}  
\vartheta (y)  = \arg \Gamma\( \tfrac{1}{4} + \tfrac{i y}{2} \)  - 
y \log  \sqrt{\pi}  
\approx \frac{y}{2} \log \( \frac{y}{2 \pi e} \)  - \frac{\pi}{8} 
\eeq
where as usual $\arg \Gamma = \Im \log \Gamma$.   

Now let us consider the $\zeta$-function, which
has the well known integral representation
\beq\label{zeta_integral}
\zeta(z) = \dfrac{1}{\Gamma(z)}\int_{0}^{\infty}
\dfrac{u^{z-1}}{e^u-1} \, du, \qquad \Re(z) > 1.
\eeq
If $z \gg 1$ is real, then the saddle point of the integrand in 
\eqref{zeta_integral} is far from the origin, and one can approximate
$e^u-1 \approx e^u$, then this integral is approximately the
same as \eqref{gamma_int_def}, showing that for real $z$ 
then $\zeta(z)\to 1$ as $z\to \infty$.
The above integral is badly  behaved at the origin $u=0$,
thus let us introduce a small parameter $\mu$, having in mind that we
can always take the limit $\mu \to 0^+$.   
The above integral arises in quantum statistical 
physics,  and $\mu$ is minus the chemical potential.   Thus we introduce
\beq\label{mu_integral}
\zeta(z) = \lim_{\mu \to 0^+}  \, \dfrac{1}{\Gamma(z)}\int_{0}^{\infty}
\dfrac{u^{z-1}}{e^{u+\mu}-1} \, du.
\eeq

We want to estimate this integral through the saddle-point method.
For this aim we can write the integrand in the form $e^{f(u)}$ where
\beq\label{ft}
f(u) = (z-1)\log u - \log\(e^{u+\mu} - 1\). 
\eeq
The  condition $f'(u) = 0$ yields the transcendental equation
\beq\label{trans}
\dfrac{z-1}{u} = \dfrac{e^{u+\mu}}{e^{u+\mu}-1},
\eeq
which can be solved explicitly, determining the saddle points 
in the following form:
\beq\label{saddle_points}
s_k(z) = z-1 +  w_k (z),
\eeq
with  $$w_k (z) =  W_k\[ (1-z)e^{1-z-\mu} \] ,$$  
where $W_k$  is the $k$-th branch of the multi-valued 
Lambert $W$ function \eqref{W_def_eq}.
Note that setting $w_k =0$,  one recovers the saddle point 
for the $\Gamma$-function.  

One can easily show that 
\beq \label{fp} 
f(s_k) = - s_k  - \mu + (z-2) \log s_k + \log (z-1)
\eeq
and 
\beq \label{fpp}
f''(s_k) =  \dfrac{(1+w_k)}{(1-z) \(1+\tfrac{w_k}{z-1}\)^2}.
\eeq

For real $z$ \eqref{saddle_points} yields only one saddle-point given 
by the principal branch $k=0$. For complex $z$ we have an infinite
number of saddle-points. However, not all of them contribute to the 
integral \eqref{mu_integral},
since according to the integration path we must choose the ones
which satisfy $\Re(s_k) > 0$.  Using equations \eqref{fp}, \eqref{fpp} 
and dividing by the Stirling approximation to $\Gamma (z)$,  
equation \eqref{stirling2}, one then obtains 
\beq\label{zeta_saddle}
\zeta(z) \approx \sum_{k}{}^{'} \zeta_k (z)
\eeq
where
\beq\label{zetak}
\zeta_k(z) = 
\exp\left\{ 
(z-1)\log\(1+\dfrac{w_k}{(z-1)}\) - w_k -\mu - \dfrac{1}{2}\log\(w_k+1\)
\right\}.  
\eeq
In the above formula \eqref{zeta_saddle} the $\mu \to 0^+$ is implicit.
The sum over integers $k$ must be taken
according to the condition $\Re(s_k) > 0$.

In obtaining \eqref{zeta_saddle} it is important to note that 
\eqref{ft} is \emph{not} an analytic function, since $\log u$ has 
branches
for $u\in \mathbb{C}$. Thus we are taking into account the contribution 
of each different branch $k$ for the saddle-points. 
This goes a little beyond the assumptions of the saddle-point
method \eqref{saddle_approx}.

\begin{figure}
\begin{center}
\begin{minipage}{.5\textwidth}
\includegraphics[width=.85\linewidth]{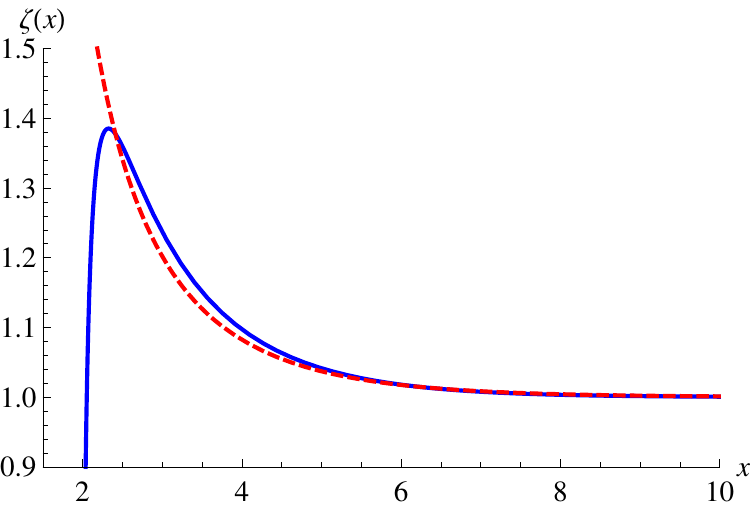}
\end{minipage}%
\begin{minipage}{.5\textwidth}
\includegraphics[width=.85\linewidth]{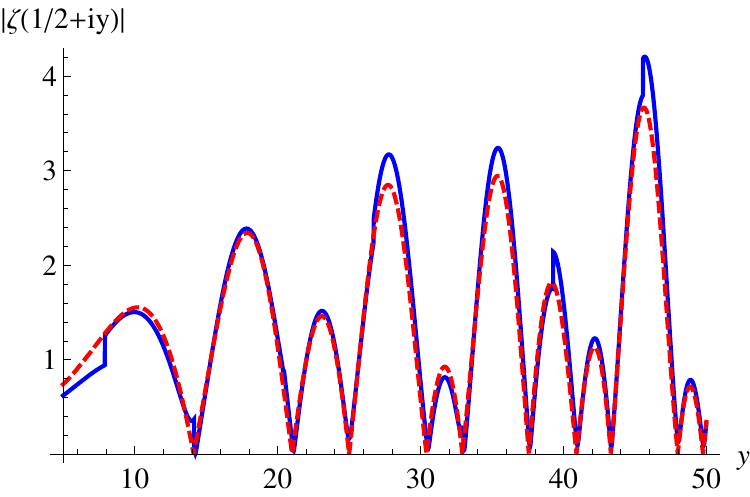}
\end{minipage}
\caption{
The dashed (red) line corresponds to the exact $\zeta(z)$ while
the solid (blue) line corresponds to the saddle point approximation
\eqref{zeta_saddle}.
\emph{ Left:} real values. 
\emph{ Right:} complex values on the critical line.}
\label{fig_complex_saddle}
\end{center}
\end{figure}

For real values of $z$ there is only one saddle-point corresponding
to $k=0$ in \eqref{zeta_saddle}. This approximation describes the
$\zeta(z)$ very  well over this range, as shown in 
Figure~\ref{fig_complex_saddle} (left).

For complex values in the approximation \eqref{zeta_saddle} we need
to consider the contribution for different branches $k$. On the
critical line $z=1/2+iy$, a plot of the absolute value of 
\eqref{zeta_saddle} is shown in Figure~\ref{fig_complex_saddle} (right). 
The number of saddle points varies with $y$. For instance, for
$y\sim 15$ we need $k=-2,-1,0$ and for $y \sim 30$ we need
$k=-4,-3,\dotsc,0$. Thus the sum in \eqref{zeta_saddle} ranges differently for
each value of $y$.  

Roughly,   the number of branches one needs 
is  $-|y|/2\pi < k \le 0$ 
which can be obtained  from the condition 
$\Re (s_k) > 0$ (the branches $k > 0$ can be neglected).
This approximation works surprisingly well,  as 
Figure~\ref{fig_complex_saddle} shows;   it clearly captures the zeros.   
However it remains difficult to characterize $\arg \zeta$ precisely,  because
of the sum over branches $k$.


\section{Discussion and open problems}
\label{sec:discussion}

In the second half of these lectures we have provided new 
insights into the RH.
We have actually proposed several concrete strategies towards 
its proof which are
all based on analysis.
Let us list some of the main open problems,   
namely what  would be needed to complete this program.   

\begin{enumerate}
\item Is there a more rigorous derivation of our main 
equation \eqref{Bzero}, 
i.e. $\theta + \theta' = (2n+1) \pi$?     
It is a powerful new equation that is valid whether
the RH is true or not and contains detailed information about the zeros since
they are  enumerated by an integer $n$.     
It is evidently correct and captures all zeros,  trivial or non-trivial, 
of the multitude of functions considered
in this paper,  including functions with zeros off of the critical line, 
where the RH fails. 
This seems to be a matter of carefully defining the limit as one 
approaches a zero.   

\item  Can one prove there is a unique solution to  the transcendental 
equation \eqref{exact_eq2} for every $n$?      
As explained,   this would imply the RH is true since then
$N_0 (T) = N(T)$.       
We provided arguments that the $\delta \to 0^+$  prescription 
smooths out the function $S(y)$ sufficiently so that there is  
indeed a unique solution.
We also showed that it is precisely this property that 
fails for the Davenport-Heilbronn function,
where the RH is false.

\item Can one prove that the real  
electric potential $\Phi (x,y)$  is a regular alternating function 
along the line $\Re (z) =1$?      
We showed this asymptotically,  but only  in a crude approximation.   
We argued that this would also establish the RH.   Can the latter be proven?

\item Can one prove directly that there are no zeros off of the critical
line by proving there are no 
solutions to $\theta (x,y) + \theta' (x, y) = (2n+1) \pi$   for $x\neq 1/2$?   
We suggested one approach,   based on the 
 observation  that the left hand side of this  equation  
 has very little dependence on $x$ for
$x>1/2$ (see Figure~\ref{trivial_zeros}).    
Since we know there are no solutions for $x\geq 1$,  i.e. there are no zeros
with $x \geq 1$,   this suggests that this would imply there are no zeros 
for $x > 1/2$, since the curve for $x>1/2$  is essentially the 
same as at $x =1$,
and there we know there are no solutions. 
 Can this argument be made more precise?   
This is an appealing
idea since the RH would  then be related to the fact that there are no 
zeros with $x \geq 1$,
which follows from the Euler product formula.    
Functions such as the Davenport-Heilbronn
function where the RH fails do not have an Euler product formula.
\end{enumerate}

\begin{added} 
A clearer  derivation of the equation \eqref{exact_eq2},  referred to as the Fran\c ca-LeClair equation
and proposed in \cite{FL1},   
was given in \cite{LecMussScattering}.
\end{added}
  
  
\begin{added}
Progress on open questions 2 and 4 above was recently presented  in \cite{LecMussScattering,SpectralFlow},  based on 
defining  a  quantum mechanical  model of a single particle scattering through impurities arranged on a circle,  where each impurity is associated with a prime number  through a scattering phase based on the Euler product.   The quantized energies $E_n (\sigma)$  are easily calculated from a Bethe-ansatz equation.    The quantized energies of this model are exactly the Riemann zeros on the critical line in the limit $\sigma \to \tfrac12$ from the right,  since they satisfy the Fran\c ca-LeClair equation.     In this model,   the pseudo-randomness of the primes implies  the particle is scattering through an essentially disordered potential, as in Anderson localization,  which implies an underlying  pseudo-random Hamiltonian.  The latter  property offers a new perspective on the origin of random matrix (GUE) statistics of the Riemann zeros.   Let us mention also that random matrix theory is used extensively in the study of disordered systems and Anderson localization.
\end{added}

\begin{added}  The spectral flow for the quantized energies  $E_n (\sigma)$ of the model in \cite{LecMussScattering}  was studied in \cite{SpectralFlow},   where it was argued that the scattering problem necessarily has real eigenvalues since the 
S-matrix is unitary,  and this would imply the Riemann Hypothesis is true.  This can be viewed as a scattering version of the 
 Hilbert-P\'olya conjecture since the S-matrix is $S=e^{-iH}$, where $H$ is a Hermitian Hamiltonian, and 
 $S^\dagger S =1$.\footnote{This is an over simplification:  in general, one  needs to consider a time-ordered exponential.}   For the RH, to crack this nut, seemingly made of steel, one
 must first create  a suitable tool before thinking 
about how to crack it precisely into two hemispheres.  
Cautiously speaking,  this recent work \cite{SpectralFlow} may establish the validity  of the Riemann Hypothesis based on quantum physics.    
This work shows that the Hilbert-P\'olya conjecture  is  true if one plucks  a  Hamiltonian of the right feather,  in particular one that leads to a \emph{scattering problem} rather than a bound state problem (as was presumably implied in the Hilbert-P\'olya conjecture).  
  However, it is probably not yet rigorous enough to constitute a  pure mathematical proof,  although 
it has all of the expected ingredients towards a proof of the Riemann Hypothesis:
  functional equation,  Euler product formula,  and pseudo-randomness of the primes that may induce GUE  random matrix statistics.
 Another nice feature of this approach  is that it easily extends to the General Riemann Hypothesis for $L$-functions based on Dirichlet characters and those based on cusp forms.   On the other hand,  the Davenport-Heilbronn counter-example has no such unitary S-matrix since the  linear combination of 
  Dirichlet $L$-functions does not admit an Euler product, showing why it fails within this framework \cite{SpectralFlow}.
\end{added}

\begin{added} 
Although the Euler product does not converge to the right of the critical line for Dirichlet $L$-functions based on principal characters,  such as Riemann $\zeta$ itself which corresponds to the trivial character 
$\chi (n)  = 1$ for all integer $n$,  a truncated Euler product can still provide a good 
  approximation.    We used this to estimate the $n=10^{100}$-th  zero to 103 decimals in \cite{Googleth}: 
Using  only  $10^{6}$ primes,  we found the
following $y_n$: 
\begin{equation*}
\begin{split}
n &= 10^{100}\mbox{th zero:}  \\ 
y_n &= 280690383842894069903195445838256400084548030162846 \\
&~~~~~ 045192360059224930922349073043060335653109252473.244\dotsm  
\end{split}
\end{equation*}
Obtaining this number took only a few minutes on a laptop using Mathematica.   The integer part is easily calculated using the Lambert approximation described above.    
 We are reasonably confident that the last $3$ digits $\sim .244$ are  correct  since   we checked that 
 they didn't 
 change between  using  $10^6$ and $5\times 10^6$ primes in the Euler product.   Although this needs more detailed investigation,  if correct  this google-th zero is far beyond what is currently known based on computational mathematics,  which are zeros around $n=10^{33}$.
\end{added}

\bigskip

\begin{acknowledgments}
We wish to thank the organizers Olaf Lechtenfeld and Ulrich Theis for the opportunity
to present these lectures.   
We also thank Denis Bernard and Timothy Healy  for useful discussions.
GF is supported by CNPq-Brazil.     For this updated version,  AL thanks Steve Gonek and Ghaith Hiary for discussions, as well as  Guilherme Fran\c ca and Giuseppe Mussardo for their subsequent 
collaboration on this topic.  
\end{acknowledgments}

\vfill\eject
\appendix

\section{The Perron formula} 
\label{sec:perron} 

Consider the sum
\beq 
A(x) = \sum_{n\leq x}{}^{'}  a(n)
\eeq
where $a(n)$ is an arithmetic function.    
Here, the prime on the summation indicates that the last term of the sum 
must be multiplied by $1/2$ when $x$ is an integer. 
One has the following integral: 
\beq
\label{Perron2}
\int_0^\infty A(x) x^{-z-1}  dx   =  
\int_0^\infty  \( \sum_{n\leq x} a(n) \)  x^{-z-1}  dx 
=  -\inv{z}  \int_0^\infty  \( \sum_{n\leq x} a(n) \)  \frac{d}{dx}  
\( x^{-z} \)   dx .
\eeq
Since $A(x)$ is a sequence of  step functions,  
and the derivative of a step function is
the Dirac delta-function,   
one has  $$\frac{d}{dx}  \sum_{n\leq x} a(n)  = a(n) \delta (x-n)$$
Integrating the above expression by parts one obtains 
\beq
\label{Perron3} 
\dfrac{g(z)}{z}  =  \int_0^\infty  A(x) x^{-z-1} dx .
\eeq
where 
\beq
g(z) = \sum_{n=1}^\infty \frac{a(n)}{n^z}.
\eeq

The standard definition of the Mellin transform $\CM f$  of the 
function $f$ is
\beq
\label{Mellin1}
(\CM f ) (s)  =  \vphi (s) =  \int_0^\infty x^{s-1} f(x)  dx  .
\eeq
Its inverse is well-known:
\beq
(\CM^{-1} \vphi ) (x)  =  f(x)  = 
\inv{2 \pi i}  \int_{c - i \infty}^{c+i \infty}  x^{-s} \vphi (s) ds .
\eeq
Identifying $s= -z$,  $\vphi (-z) = g(z)/z$,  and $f(x) = A(x)$  
one obtains  \eqref{Perron}.

\section{Counting formula on the entire critical strip}
\label{sec:number_zeros}

The argument principle enables us to count zeros and poles of a complex
function inside a simply connected bounded region. Let us consider two 
functions $f(z)$ and
$g(z)$, such that $f$ is analytic and $g$ is meromorphic inside
and on a given closed contour $\CC$. Thus we can write
\beq
g(z) = 
(z-w_1)^{n_1}
(z-w_2)^{n_2} \dotsm 
\dfrac{1}{(z-\bar{w}_1)^{\bar{n}_1}
(z-\bar{w}_2)^{\bar{n}_2}\dotsm} \, h(z)
\eeq
where $w_i$ and $\bar{w}_i$ denotes the zeros and poles, respectively, with
its respective multiplicities $n_i$ and order $\bar{n}_i$, and
$h(z)$ is analytic and without zeros or poles in and on $\CC$. 
Then we have
\beq
\dfrac{g'(z)}{g(z)} = \dfrac{d}{dz} \log g(z) = 
\dfrac{n_1}{z-w_1} + 
\dfrac{n_2}{z-w_2} + \dotsm -
\dfrac{\bar{n}_1}{z-\bar{w}_1} -
\dfrac{\bar{n}_2}{z-\bar{w}_2} - \dotsm + \dfrac{h'(z)}{h(z)}.
\eeq
Assuming that $g(z)$ has no zeros or poles on the contour $\CC$, we 
thus have
\beq
\oint_{\CC} f(z)\dfrac{g'(z)}{g(z)}dz = 
\sum_{i} n_i \oint_{\CC} \dfrac{f(z)}{z-w_i} dz -
\sum_{j} \bar{n}_j \oint_{\CC} \dfrac{f(z)}{z-\bar{w}_j} dz 
+\oint_{\CC}f(z)\dfrac{h'(z)}{h(z)}dz .
\eeq
Since $\tfrac{f h'}{h}$ is analytic, the last integral vanishes
by Cauchy's integral theorem. By Cauchy's integral formula
we therefore have
\beq
\dfrac{1}{2\pi i}\oint_{\CC} f(z)\dfrac{g'(z)}{g(z)} dz = 
\sum_{i} n_i f(w_i) - \sum_{j}\bar{n}_j f(\bar{w}_j).
\eeq
Now setting $f(z) = 1$ we have \emph{Cauchy's argument principle}
\beq\label{arg_principle}
\dfrac{1}{2\pi i} \oint_{\CC} \dfrac{g'(z)}{g(z)} dz
= N - \bar{N},
\eeq
where $N = \sum_{i} n_i$ is the number of zeros inside $\CC$, including
multiplicities, and $\bar{N} = \sum_{j} \bar{n}_j$ is the number
of poles, accounting its orders.

\begin{figure}
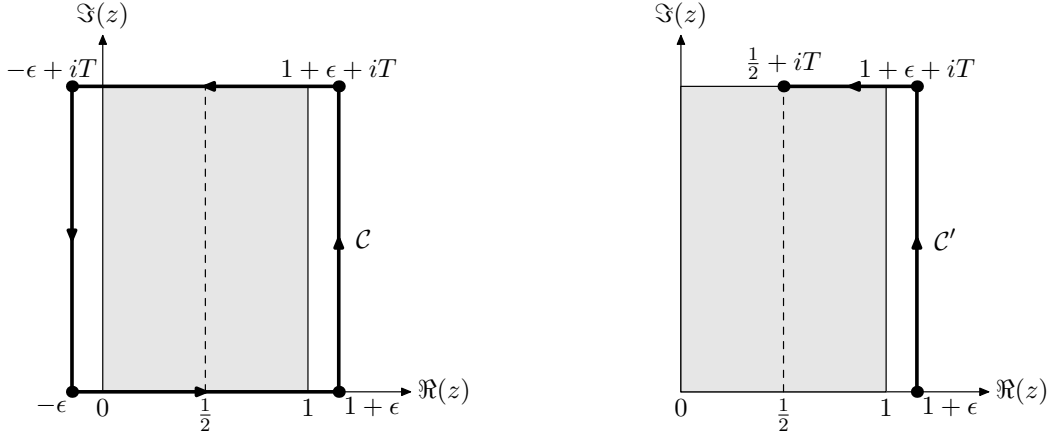

\begin{center}
\begin{minipage}{.49\textwidth}
\includegraphics[scale=0.9]{figs/fig5.mps}
\end{minipage}%
\begin{minipage}{.49\textwidth}
\includegraphics[scale=0.9]{figs/fig6.mps}
\end{minipage}
\caption{\emph{ Left:} contour $\CC$ for \eqref{arg_counting} in the 
counterclockwise direction. 
The contour must have no zeros on it. The critical line
$\Re(z) = 1/2$ splits the rectangle in half. \emph{ Right:} it is equivalent
to consider the contour $\CC'$.}
\label{fig:rectangle_counting}
\end{center}
\end{figure}

Now let us use \eqref{arg_principle} to count the zeros of $\zeta(z)$ inside
the region $0 \le \Re(z) \le 1$ and $0 \le \Im(z) \le T$.
Remember that $\zeta(z)$ and $\chi(z)$ have the same zeros in 
this region. Since $\zeta(z)$ has a simple pole at $z=1$, and so has
$\chi(z)$, we consider the function
\beq\label{xi_zeta}
\xi(z) \equiv \tfrac{1}{2}z(z-1)\pi^{-z/2}\,\Gamma\(z/2\) \zeta(z),
\eeq
which is entire and has no poles. Note that $\zeta(z)$, and thus $\xi(z)$, 
has no zeros along the positive real line, 
and according to \eqref{arg_principle} we must
assume that $T$ \emph{does not} correspond to the imaginary part of a zero.
Then from \eqref{arg_principle} we have
\beq\label{arg_counting}
N(T) = \dfrac{1}{2\pi}\Im\[
\oint_{\CC} \dfrac{\xi'(z)}{\xi(z)} dz \],
\eeq
where the contour $\CC$ is chosen to be the boundary of the rectangle
with vertices at the points $z\in \left\{-\epsilon, \, 1+\epsilon, \, 
1+\epsilon + iT, \, -\epsilon + iT\right\}$, as illustrated in
Figure~\ref{fig:rectangle_counting} (left). Since $\xi(z)$ is real on the real
line, the integration along the line segment $[-\epsilon, 1+\epsilon]$ 
does not contribute to \eqref{arg_counting}. 
Moreover, since $\xi(z)=\xi(1-z)$ we
can consider just the right half of the rectangle, whose contribution to
the integral \eqref{arg_counting} comes from line segments
$[1+\epsilon, 1+\epsilon + iT]$ and $[1+\epsilon+iT, 1/2+iT]$. This
is the contour $\CC'$ shown in Figure~\ref{fig:rectangle_counting} (right). 
Therefore we have
\begin{align}
N(T) &=  \dfrac{1}{\pi}\Im
\int_{\CC'} dz \, \dfrac{d }{dz} \log\xi(z)  \\
&= \dfrac{1}{\pi} \Im \int_{\CC'} dz \, \dfrac{d }{dz} 
\[ \log\(\tfrac{1}{2}z(z-1)\) +
\log\(\pi^{-z/2} \Gamma\( z/2 \) \) +
\log\zeta(z) \] . \label{n_integral}
\end{align}
The first term in \eqref{n_integral} gives
\beq
\dfrac{1}{\pi}\Im  
\left.\log\[\tfrac{1}{2}z(z-1)\] \right\vert_{1+\epsilon}^{1/2+iT} 
= \dfrac{1}{\pi} \Im  
\log\( - \tfrac{1}{2}\( T^2 + \tfrac{1}{4} \) \) = 1.
\eeq
For the second term in \eqref{n_integral}, since $\pi^{-z/2}\Gamma(z/2)$
is analytic and real on $z=1+\epsilon$, we have the contribution only
from the end point of the contour, which gives
\beq\label{second_term_integral}
\dfrac{1}{\pi} \Im  
\log \Gamma\(\tfrac{1}{4} + i\tfrac{T}{2}\)
- \dfrac{T}{2\pi}\log\pi \equiv \dfrac{1}{\pi}\vartheta\(T\),
\eeq
where we have the Riemann-Siegel $\vartheta$ function.
Therefore, the number of zeros in the critical strip up to height $T$ 
is given by
\beq\label{backlund_strip}
N(T) = \dfrac{1}{\pi}\vartheta(T) + 1 + S(T)
\eeq
where
\beq\label{S_variation}
S(T) \equiv \dfrac{1}{\pi} \Im
\int_{\CC'} d\( \log \zeta(z)\) = 
\dfrac{1}{\pi}\Delta_{\CC'}\arg\zeta(z).
\eeq
The above formula \eqref{backlund_strip} is known as the \emph{Backlund
counting formula}.
Note that $\Im\log \zeta  = \arg \zeta$. The
above result consists in computing the variation of $\arg\zeta(z)$ through
the line segments starting at $z=1+\epsilon$, 
where $\arg\zeta(1+\epsilon) = 0$, then up
to $z=1+\epsilon + iT$, then finally to $z=1/2+iT$. Since the integrand
in \eqref{S_variation} is a total derivative, it is tempting to
write $S(T) = 
\tfrac{1}{\pi}\Im\log\zeta\(1/2+iT\)=\tfrac{1}{\pi}\arg\zeta(1/2+iT)$, 
corresponding to the end point of the contour $\CC'$. However, this
must be carefully considered since $\log\zeta(z)$ can not be always
in the principal branch along $\CC'$. 
Nevertheless, if $\Re\(\zeta\)  > 0$ for all $z$ on $\CC'$, 
then $\log\zeta(z)$ is well defined and is always in 
the principal branch, i.e. $\log\zeta(z) = \log|\zeta(z)| + i\arg\zeta(z)$
with $-\pi/2 < \arg\zeta(z) \le \pi/2$.
In such a case it is valid
to set $S(T) = \tfrac{1}{\pi}\arg\zeta(1/2+iT)$, and we also have
$|S(T)| < 1/2$. 

Let us now consider $z=x+iy$ for $x \ge 2$. From the 
series $\zeta(z) = \sum_{n=1}^{\infty}n^{-z}$ we have
\beq
\Re\(\zeta\) = 1 + 
\dfrac{\cos\(y\log2\)}{2^x} + \dfrac{\cos\(y\log3\)}{3^x} + \dotsm
 \ge 1 - \( \dfrac{1}{2^x} + \dfrac{1}{3^x} + \dotsm \).
\eeq
Since $x\ge2$ we also have
\beq
\dfrac{1}{2^x} + \dfrac{1}{3^x} + \dfrac{1}{4^x} + \dotsm \le 
\dfrac{1}{2^2} + \dfrac{1}{3^2} + \dfrac{1}{4^2} + \dotsm  
= \zeta(2) - 1 = \dfrac{\pi^2}{6}-1
\eeq
Therefore
\beq\label{ReZeta2}
\Re\(\zeta\) > \dfrac{12-\pi^2}{6} > 0.
\eeq
For this reason it is often assumed $\epsilon = 1$ in
the contour $\CC'$ appearing in \eqref{S_variation}. In this way one only
needs to analyze the behaviour of $\Re \( \zeta \)$ along the 
horizontal line joining the points $2+iT$ and $\tfrac{1}{2}+iT$.

If one expands $\Gamma(z)$ in \eqref{second_term_integral} through 
Stirling's formula one obtains
\beq
\vartheta(T) = \dfrac{T}{2}\log\(\dfrac{T}{2\pi e}\) - \dfrac{\pi}{8}
+ O\(T^{-1}\).
\eeq
Then \eqref{backlund_strip} yields the \emph{Riemann-von Mangoldt counting
formula}
\beq
N(T) =  \dfrac{T}{2\pi}\log\(\dfrac{T}{2\pi}\) -\dfrac{T}{2\pi} + 
\dfrac{7}{8} + S(T) + O\(T^{-1}\).
\eeq

\section{Mathematica code}
\label{sec:mathematica}

Here we provide a simple Mathematica implementation
to compute the zeros of $L$-functions based on the previous
transcendental equations.
We will consider Dirichlet $L$-functions, since it involves more 
ingredients, like the
modulus $k$, the order $a$ and the Gauss sum $G(\tau)$. For the Riemann
$\zeta$-function the procedure is trivially adapted, as
well as for the Ramanujan $\tau$-function of section \ref{sec:ramanujan}.

The function \eqref{RSgen} is implemented as follows:
\begin{lstlisting}
RSTheta[y_, a_, k_] := Im[LogGamma[1/4 + a/2 + I*y/2]] - y/2*Log[Pi/k]
\end{lstlisting}
For the transcendental equation \eqref{exact} we have
\begin{lstlisting}
ExactEq[n_, y_, s_, a_, k_, j_, G_,  n0_] := 
    (RSTheta[y, a, k] + Arg[DirichletL[k, j, 1/2 + s + I*y]] -1/2*Arg[G])/Pi + a/4 + 1/2 - n + n0
\end{lstlisting}
Above, \lstinline{s} denotes $0<\delta\ll 1$, \lstinline{a} is the 
order \eqref{order}, \lstinline{k} is the
modulus, \lstinline{j} specify the Dirichlet character 
$\chi_{k,j}$ (as discussed
in section \ref{sec:dirichlet}), and \lstinline{G} is the Gauss
sum \eqref{tau}. Note that we also included \lstinline{n0}, discussed
after \eqref{almost_final}, but we always set $n_0=0$ for the cases
analysed in section \ref{sec:numerical_lfunc}. The implementation of the
approximate solution \eqref{approx} is
\begin{lstlisting}
Sgn[n_] := Which[n != 0, Sign[n], n == 0, -1]
A[n_, a_, G_, n0_] := Sgn[n]*(n - n0 + 1/2/Pi*Arg[G]) + (1 - 4*Sgn[n] - 2*a*(Sgn[n] + 1))/8
yApprox[n_, a_, G_, k_, n0_] := 2*Pi*Sgn[n]*A[n, a, G, n0]/LambertW[k*A[n, a, G, n0]/E]
\end{lstlisting}
One can then obtain the numerical solution of the transcendental 
equation \eqref{exact} as follows:
\begin{lstlisting}
FindZero[n_, s_, a_, k_, j_, G_, n0_, y0_, prec_] := 
 y /. FindRoot[ExactEq[n, y, s, a, k, j, G, n0], {y, y0}, PrecisionGoal -> prec/2, AccuracyGoal -> prec/2, WorkingPrecision -> prec]
\end{lstlisting}
Above, \lstinline{y0} will be given by the approximate solution \eqref{approx}.
The variable \lstinline{prec} will be adjusted iteratively.
Now the procedure described in section \ref{sec:numerical_exact}
can be implemented as follows:
\begin{lstlisting}
DirichletNZero[n_, order_, digits_, k_, j_, n0_] := (
    chi = DirichletCharacter[k, j, -1];
    a = Which[chi == -1, 1, chi == 1, 0];
    s = 10^(-3);
    prec = 15;
    G = Sum[DirichletCharacter[k, j, l]*Exp[2*Pi*I*l/k], {l, 1, k}];
    y = N[yApprox[n, a, G, k, n0], 20];
    absvalue = 1;
    While[absvalue > order,
        y = FindZero[n, s, a, k, j, G, n0, y, prec];
        Print[NumberForm[y, digits]];
        s = s/1000;
        prec = prec + 20;
        absvalue = Abs[DirichletL[k, j, 1/2 + I*y]];
    ]
    Print[ScientificForm[absvalue, 5]];
)
\end{lstlisting}
Above, the variable \lstinline{order} controls the accuracy of the solution. 
For instance, if \lstinline{order=10^(-50)}, it will iterate until the solution
is verified at least to $|L\(\tfrac{1}{2})+iy\)| \sim 10^{-50}$. The
variable \lstinline{digits} controls the number of decimal 
places shown in the output.

Let us compute the zero $n=1$, for the character \eqref{char73},
i.e. $k=7$ and $j=3$. We will verify the solution to $\sim 50$ decimal
places and print the results with $52$ digits. Thus executing
\begin{lstlisting}
DirichletNZero[1, 10^(-50), 52, 7, 3, 0]
\end{lstlisting}
the output will be
\begin{lstlisting}
4.35640188194945
...
4.356401624736284227279574790515519132971499551683092
4.356401624736284227279574790515519132971499294412496
4.356401624736284227279574790515519132971499294412239
4.1664*10^(-55)
\end{lstlisting}
Note how the decimal digits converge in each iteration. 
By decreasing \lstinline{order} and
increasing \lstinline{digits} it is possible to obtain  
highly accurate solutions.
Depending on the height of
the critical line under consideration, one should adapt the 
parameters \lstinline{s} and \lstinline{prec} appropriately. 
In Mathematica we were able to
compute solutions up to  $n\sim10^6$ for Dirchlet $L$-functions, and up to 
$n\sim10^9$ for the Riemann $\zeta$-function without problems.  
We were unable to go much higher only because Mathematica could not 
compute the $arg L$ term reliably.
To solve the transcendental
equations \eqref{exact_eq2} and \eqref{exact} for very high values
on the critical line is a challenging numerical problem. Nevertheless,
we believe that it can be done through a more specialized implementation.



%

\end{document}